\documentclass{ws-rmta}
\usepackage{psfrag,graphicx,epsfig,subfigure,hyperref,xcolor}
\newcommand*{\dif}{\mathop{}\!\mathrm{d}}
\newcommand{\re}{\textcolor{black}}
\begin{document}

\markboth{J. Ge et al.}
{LD-RMT and Its Applications in Deep Learning and Wireless Communications}

%
\catchline{}{}{}{}{}
%

\title{\re{Large-Dimensional Random Matrix Theory and Its Applications} in Deep Learning and Wireless Communications}

\author{Jungang Ge, and Ying-Chang Liang\footnote{
Corresponding author.}}

\address{Center for Intelligent Networking and Communications (CINC)\\
University of Electronic Science and Technology of China (UESTC)\\
Chengdu 611731, P. R. China\\
\email{gejungang@std.uestc.edu.cn, liangyc@ieee.org}}

\author{Zhidong Bai}

\address{Key Laboratory for Applied Statistics of MOE\\
School of Mathematics and Statistics\\
Northeast Normal University\\
Changchun 130024, P. R. China\\
\email{baizd@nenu.edu.cn}}

\author{Guangming Pan}

\address{School of Physical and Mathematical Sciences\\
Nanyang Technological University, Singapore 637371\\
\email{gmpan@ntu.edu.sg}}

\maketitle


\begin{abstract}
Large-dimensional random matrix theory, RMT for short, which originates from the research field of quantum physics, has shown tremendous capability in providing deep insights into large dimensional systems. With the fact that we have entered an unprecedented era full of massive amounts of data and large complex systems, RMT is expected to play more important roles in the analysis and design of modern systems. In this paper, we \re{review} the key results of RMT and its applications in two emerging fields: wireless communications and deep learning. In wireless communications, we show that RMT can be exploited to design the spectrum sensing algorithms for cognitive radio systems and to perform the design and asymptotic analysis for large communication systems. In deep learning, RMT can be utilized to analyze the Hessian, input-output Jacobian and data covariance matrix of the deep neural networks, thereby to understand and improve the convergence and the learning speed of the neural networks.  Finally, we highlight some challenges and opportunities in applying RMT to the practical large dimensional systems.
\end{abstract}



\section{Introduction} \label{sec:intro}

In the early $1940$'s, large dimensional random matrix theory, RMT for short, was first employed to study the complicated organizational structure of the heavy nuclei in the quantum mechanics. In particular, the $N\times N$ Hamiltonian matrices, whose elements are drawn from a probability distribution, are advocated to approximate the complex nuclei systems \cite{bai2014spectral, bai2010spectral}. Afterwards, the well-known {\it Wigner matrix} and {\it semicircular law} were proposed \cite{wigner1955characteristic}. From then on, RMT has been developed rapidly and many interesting results have emerged \cite{bai2014spectral, tulino2004random, tao2012topics, bai2010spectral, couillet2011random}, e.g., circular law, Mar\v{c}enko-Pastur law, etc. Nowadays, RMT has become an important research field in quantum physics and mathematics.

On the other hand, over the past decades, the computation speed and the storage capability of the computer has been increasing significantly due to the rapid development of the computer science, thus massive amounts of data can be collected and stored. RMT is regarded as a powerful tool to reveal the hidden patterns behind the large dimensional data. Besides, as the practical systems grow more and more complex, meaningful insights can be drawn through RMT. RMT has thus shown extraordinary capability in various research fields such as finance statistics, wireless communications, deep learning, etc. In this paper, we focus on the applications of RMT in the two emerging fields: wireless communications and deep learning.

In modern wireless communication systems, the number of users and the wireless traffic have been growing exponentially according to the report of {\it International Telecommunication Union} \cite{union2015imt}. As a consequence, the communication systems have to involve more degrees of freedom to support the communication demands. From one aspect, the degrees of freedom can be acquired by increasing the length of the spreading sequences in code-division-multiple-access (CDMA) systems. In another aspect, for the multiple-input-multiple-output (MIMO) systems, larger antenna arrays are employed to provide more degrees of freedom. It is worth noting that both the spreading sequences in CDMA systems and the channel matrices in MIMO systems can be modeled with random matrices and thus the systems can be analyzed by RMT \cite{tse1999linear,verdu1999spectral, tse2000linear}. For example, the capacity of the MIMO systems is related to the singular values of the channel matrix. With the knowledge of Wishart matrices and Mar\v{c}enko-Pastur law, the system capacity can be determined from the spectrum of the gram matrix of the channel matrix \cite{tulino2004random}. In addition, RMT is utilized to evaluate the asymptotic performance of the extremely large dimensional systems where both the number of users and that of the degrees of freedom go to infinity \cite{tse1999linear}. In return, the asymptotic results can provide us constructive instructions for the design of the large complex communication systems. On the other hand, improving the spectrum efficiency is also an effective method to accommodate the explosive wireless traffic. {\it Cognitive radio} (CR) technique provides us a novel way to further enhance the spectrum efficiency, i.e., allowing the so-called secondary users to use the licensed spectrum without disturbing the licensed primary users \cite{liang2008sensing,zeng2009eigenvalue,zeng2010review, liang2011cognitive, liang2020dynamic}. In the opportunistic CR, the secondary users have to determine whether the interested spectrum is occupied by the primary users through analyzing the signals sampled from the radio environment, and this is known as the spectrum sensing technique. In essence, the spectrum sensing problems are the conventional signal detection problems, i.e., identifying the existence of the primary users according to the signal samples from the radio environment. In the multi-antenna scenarios and the cooperative sensing scenarios, we can obtain sampled signal vectors (each vector is a signal sample) via multiple antennas or multiple sensors, respectively. Then the sample covariance matrix can be computed by the temporal signal samples acquired within the sensing duration. With RMT, it is observed that the sample covariance matrix when primary users are absent can be modeled with Wishart matrix and the sample covariance matrix when primary users are present can be modeled with the spiked model \cite{bianchi2009performance, bianchi2011performance, penna2009probability}. Consequently, many eigenvalue-based spectrum sensing algorithms have been developed upon this observation, i.e., by determining which random matrix model the sample covariance matrix should belong to.

Deep learning is regarded as the most significant breakthrough in the field of machine learning over the past two decades. It has shown that the state-of-the-art results in many areas such as computer vision, natural language processing, and human games are obtained by deep learning techniques\cite{lecun2015deep}. The strength of deep learning comes from the extremely complex deep neural networks, which are usually composed of millions or sometimes even billions of parameters \cite{adlam2019random}. The large complex neural networks are so powerful that they can approximate almost all possible functional relations between the inputs and the outputs. In addition, many advanced neural networks are proposed to extract the hidden patterns behind the large dimensional datasets, e.g., convolutional neural networks (CNNs), and recurrent neural networks (RNNs). However, the neural networks are often treated as black boxes with merely visible input-ports and output-ports since the neural networks and the datasets are too complex to understand due to their extremely large dimensions. This is quite similar to the dilemmas that are usually encountered in the quantum physics. With the fact that the large complex systems in quantum physics can be well approximated with random variables, we can also model the large complex neural networks with random variables. In addition, it is known that the neural networks are randomly initialized in general and the training stage may introduce only low-rank perturbations around the random configuration. This further justifies the assumption about the randomness in the neural networks. Therefore, RMT is expected to shed some light on understanding the neural networks. The recent research results have shown that the RMT-based analysis framework for the random neural networks can help us to understand and improve the deep learning technology. For example, it is observed that keeping all the singular values of the input-output Jacobian concentrate around $1$ can dramatically speed up the learning process \cite{schoenholz2016deep, pennington2017resurrecting, pennington2018emergence}. Moreover, the input-output Jacobian can be decomposed as a product of random matrices, and the characteristics of the singular values of the input-output Jacobian can be studied via RMT. The results can provide us constructive instructions to improve the performance of the deep neural networks by choosing the depth, the random weight initializations and the nonlinear activation functions. In addition, the Hessian of the neural networks contains a lot of information about the loss surface, and the spectrum of the Hessian at the critical points can be utilized to identify the saddle points or the local minima \cite{bray2007statistics, dauphin2014identifying, pennington2017geometry}. In the simplest case with several impractical assumptions, it is shown that the Hessian can be decomposed as a summation of a Wishart matrix and a Wigner matrix \cite{pennington2017geometry}. Thus, the spectrum of the Hessian can be analyzed using the results in RMT. In a recent work \cite{granziol2020beyond}, more complex random matrix models, such as random Wigner/Wishart ensemble products and percolated Wigner/Wishart ensembles, are proposed to approximate the Hessian more accurately. Besides, with the fact that highly skewed distributions means strong anisotropy in the embedded feature space which will derail the learning process, RMT is employed to study the spectra of data covariance matrices in the neural networks \cite{adlam2019random, pennington2017nonlinear}. The analytical results for the data covariance matrices help us identify a large series of activation functions that can preserve the spectra as the signal propagates through the neural networks. This also gives us some guidelines for designing new activation functions.
\re{
Last but not least, the limiting train error and generalization error of the overparametrized random neural networks can be analytically derived via RMT. The results exactly reveal the so-called {\it double descent phenomenon}, which explains the reason why the overparametrized neural networks with zero train error can generalize well without overfitting. Hence, this provides us deep insights into the outstanding performance of modern deep neural networks with millions or sometimes even billions of parameters. Furthermore, RMT can be also exploited to perform spectral analysis over the kernel matrices (e.g, conjugate kernel, neural tangent kernel) that are closely related to the training process of neural networks. The spectral behaviors of these kernel matrices in turn provide us possibly efficient ways to understand the training of neural networks.
}

Although there exist several classical books that review the basics of RMT and investigate the applications in wireless communications, e.g., \cite{bai2014spectral, tulino2004random, couillet2011random}, many new progresses have been made in recent years and not be included in the books. On the other hand, data science has become an important branch in modern digitalized society since the explosive data can be exploited via some advanced techniques to bring people great convenience. Machine learning, especially deep learning, is regarded as the most attractive technique that can extract a lot of beneficial knowledge from the big data. Intriguingly, many recent works show that RMT can also be utilized to help us understand and improve the deep learning technique. Therefore, in this paper, we try to provide a comprehensive sketch of the applications of RMT, including the latest applications in wireless communications and the recent progresses made in deep learning. We hope this article can establish a connection between engineering applications and mathematical field in which RMT will keep to be powerful.

The remainder of this paper is organized as follows. In Section \ref{sec:rmt}, we introduce the basic concepts and typical results in RMT. Section \ref{sec:spectrum_sensing} reviews the applications of RMT in designing spectrum sensing algorithms in cognitive radio systems. Section \ref{sec:LinearMURx} shows that RMT can be employed to analyze the asymptotic performance of the multiuser receivers in large communication systems. In Section \ref{sec:deep_learning}, we investigate some rudimentary explorations that apply RMT in understanding and improving the performance of neural networks. Important challenges and opportunities are discussed in Section \ref{sec:challenges}. Finally, Section \ref{sec:conclusions} concludes this paper.

\section{Basics of Large-Dimensional Random Matrix Theory}\label{sec:rmt}

In RMT, the results usually focus on the asymptotic regimes where the dimensions of the random matrices are extremely large or even infinite. The limiting results obtained in infinite-dimension cases can stunningly approximate the more practical finite-dimension scenarios very well, and this has been validated by many empirical results. Hence, it is quite significant to study the limiting behaviors of the random matrices. In this section, we introduce the basic concepts and celebrated results in RMT, which provide powerful theoretical support for analyzing the large dimensional communication systems and the emerging deep neural networks.

\subsection{Definitions and Notations}\label{subsec:definitions}

As the name suggests, a random matrix is a matrix whose entries are random variables. The behaviors of eigenvalues and eigenvectors of a random matrix are of main interest in RMT. In particular, most works focus on the characteristics of the eigenvalues (a.k.a. the spectrum) of the random matrices \cite{tulino2004random,bai2010spectral, couillet2011random, tao2012topics, bai2014spectral}. In addition, the spectra of Hermitian matrices are widely studied since their eigenvalues are real. Some definitions about the spectrum of a Hermitian matrix are as follows.

\begin{definition}
For an $N\times N$ (non-necessarily random) Hermitian (self-adjoint) matrix $\mathbf{T}_{N}$, its {\it empirical spectrum density} ({\it e.s.d.}) is defined as
\begin{equation}\label{eq:esd}
F^{\mathbf{T}_{N}}(x) = \frac{1}{N}\sum\limits_{j=1}^{N}1_{\{x,\lambda_{j}\leq x\}}(x).
\end{equation}
where $\lambda_{1},\cdots,\lambda_{N}$ are the eigenvalues of $\mathbf{T}_{N}$, $1_{\{x,\lambda_{j}\leq x\}}(x)$ is the indicator function, which equals $1$ when $\lambda_{j}\leq x$ or $0$ otherwise. When the dimension of $\mathbf{T}_{N}$ becomes large, or even goes to the infinity, i.e., $N\rightarrow \infty$, if its {\it e.s.d.}, namely, $F^{\mathbf{T}_{N}}$, converges to a non-random limit distribution $F^{\mathbf{T}}$, then $F^{\mathbf{T}}$ is defined as the {\it limit spectrum distribution} ({\it l.s.d.}) of $\mathbf{T}_{N}$. Most results in RMT are based on the weak convergence of $F^{\mathbf{T}_{N}}$ to $F^{\mathbf{T}}$, i.e., for all $x$ where $F^{\mathbf{T}}$ is continuous, $F^{\mathbf{T}_{N}}(x)-F^{\mathbf{T}}(x)\to 0$. The weak convergence is often denoted by
\begin{equation}\label{eq:lsd}
F^{\mathbf{T}_{N}}\Rightarrow F^{\mathbf{T}}.
\end{equation}
Although the weak convergence of $F^{\mathbf{T}_{N}}$ to $F^{\mathbf{T}}$ only holds for some specific random matrices in most cases, this will be described with the phrase $F^{\mathbf{T}_{N}}\Rightarrow F^{\mathbf{T}}$ almost surely (a.s.), which is also denoted by $F^{\mathbf{T}_{N}}\stackrel{a.s.}{\Rightarrow} F^{\mathbf{T}}$ in this article.
\end{definition}

\subsection{Semicircular Law and Mar\v{c}enko-Pastur Law}\label{subsec:MP_law}

The most well-known random matrices in RMT are Wishart matrices \cite{wishart1928generalised} and Wigner matrices \cite{wigner1955characteristic, wigner1958distribution}, which have been studied thoroughly since both of the two kinds of random matrices are Hermitian.

\begin{definition}
An $N\times N$ matrix $\mathbf{X}_N$ is a {\it Wigner matrix} if it is a Hermitian random matrix whose upper-triangular entries are independent zero-mean random variables with identical variance. $\mathbf{X}_N$ is referred to as a standard Wigner matrix when the identical variance is $\frac{1}{N}$.
\end{definition}

\begin{figure}[!t]
\begin{center}
\psfrag{x}[cc][cc][.7][0]{$x$}
\psfrag{y}[cc][cc][.7][0]{\rm probability density, $f^{\mathbf{X}}(x)$}
\epsfxsize=0.7\textwidth \leavevmode
\epsffile{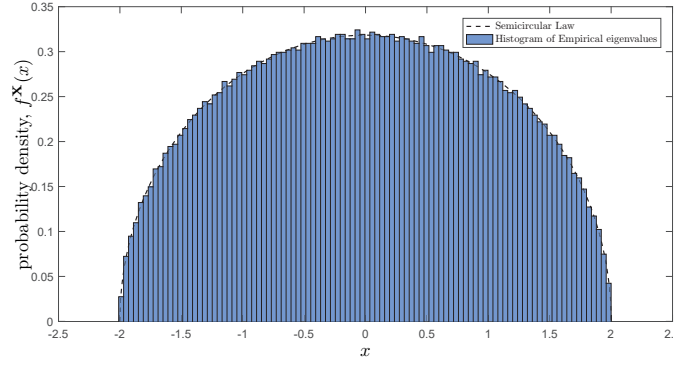}
\caption{Histogram of empirical eigenvalues and the semicircular law when $N=10000$.}\label{fig:SC_law}
\end{center}
\end{figure}

\begin{theorem}\label{th:SC_law}
Consider an $N\times N$ random Hermitian matrix $\mathbf{X}_N$ with independent entries $\mathbf{X}_{N,ij}$ such that $\mathbb{E}[\mathbf{X}_{N,ij}]=0$, $\mathbb{E}[|\mathbf{X}_{N,ij}|^2]=1/N$, and $\mathbf{X}_{N,ij}$ has a moment of order $2+\epsilon$ for an existing $\epsilon$, as $N\to\infty$, its {\it e.s.d.} converges weakly and almost surely towards a non-random distribution whose {\it probability density function} ({\it p.d.f.}), namely, $f^{\mathbf{X}}$, is given by \cite{bai2010spectral}
\begin{equation}\label{eq:SC_law}
f^{\mathbf{X}}(x) =
\left\{
\begin{aligned}
&\frac{1}{2\pi}\sqrt{4-x^2} &,{\rm if}\ |x|\leq2,\\
&0 &,{\rm otherwise}.
\end{aligned}
\right.
\end{equation}
\end{theorem}

As shown in Fig. \ref{fig:SC_law}, the graph of its {\it p.d.f.} looks like a semi-circle, and {\it Theorem} \ref{th:SC_law} is known as the {\it semicircular law}. In addition, the requirement of the moment of order $2+\epsilon$ can be discarded if the entries are {\it independent and identically distributed} ({\it i.i.d.}) \cite{couillet2011random}. Further, for a more general case where the identical variance of the entries becomes $\sigma^2/N$, the {\it e.s.d.} can be describe with the generalized semicircular law with an additional parameter $\sigma$. The semicircular law parameterized by $\sigma$ is given by

\begin{equation}\label{eq:SC_law_scaled}
f_{SC}(x; \sigma)=
\left\{
\begin{aligned}
&\frac{1}{2\pi\sigma^2}\sqrt{4\sigma^2-x^2} &,{\rm if}\ |x|\leq2\sigma,\\
&0 &,{\rm otherwise}.
\end{aligned}
\right.
\end{equation}

\begin{definition}
If the columns of the $N\times n$ random matrix $\mathbf{X}_N$ are zero-mean independent (real or complex) Gaussian vectors with covariance matrix $\mathbf{\Sigma}_N$, then the $N\times N$ random matrix $\mathbf{X}_N\mathbf{X}_N^H$ is a central (real or complex) Wishart matrix with $n$ degrees of freedom and covariance matrix $\mathbf{\Sigma}_N$. This is often denoted by $\mathbf{X}_N\mathbf{X}_N^{T}\sim\mathcal{W}_{N}(n, \mathbf{\Sigma}_N)$ for real Wishart matrices and $\mathbf{X}_N\mathbf{X}_N^{H}\sim\mathcal{CW}_{N}(n, \mathbf{\Sigma}_N)$ for complex Wishart matrices. Particularly, the Wishart matrix such that $\mathbf{\Sigma}_N=\mathbf{I}_{N}$ is also referred to as the {\it zero (or null) Wishart matrix}.
\end{definition}

\begin{remark}
Let $\mathbf{x}_1, \mathbf{x}_2, \cdots, \mathbf{x}_n\in\mathbb{C}^{N}$ be $n$ independent samples from a random process $\mathbf{x}\sim\mathcal{CN}(0, \mathbf{\Sigma}_N)$. Then we concatenate the $n$ samples to form a {\it sample matrix} $\mathbf{X}_N=[\mathbf{x}_1, \mathbf{x}_2, \cdots, \mathbf{x}_n]$. Hence, we have
\begin{equation}\label{eq:sample_matrix}
\sum\limits_{i=1}^{n}\mathbf{x}_{i}\mathbf{x}_i^H = \mathbf{X}_N\mathbf{X}_N^H.
\end{equation}
The term on the right side of \eqref{eq:sample_matrix}, namely, $\mathbf{X}_N\mathbf{X}_N^H$, is the Gram matrix of the random matrix $\mathbf{X}_N$ and the term on the left side, i.e., $\sum_{i=1}^{n}\mathbf{x}_{i}\mathbf{x}_i^H$, is related to the {\it sample covariance matrix} of the random process $\mathbf{x}$, which is defined as
\begin{equation}\label{eq:scm}
\hat{\mathbf{R}}_{\mathbf{x}\mathbf{x}}\triangleq\frac{1}{n}\sum\limits_{i=1}^{n}\mathbf{x}_{i}\mathbf{x}_i^H=\frac{1}{n}\mathbf{X}_N\mathbf{X}_N^H
=\left(\frac{1}{\sqrt{n}}\mathbf{X}_N\right)\left(\frac{1}{\sqrt{n}}\mathbf{X}_N\right)^H.
\end{equation}
Besides, $\mathbf{\Sigma}_N$ is referred to as the {\it population covariance matrix} of the random process $\mathbf{x}$. In signal detection problems, we will see that the sample covariance matrix under the pure noise case becomes a Wishart matrix. This is exactly the origin of the {\it null Wishart matrix} terminology.
\end{remark}

For a random process $\mathbf{x}$, we denote its {\it population covariance matrix} and {\it sample covariance matrix} by $\mathbf{R}_{\mathbf{x}\mathbf{x}}$ and $\hat{\mathbf{R}}_{\mathbf{x}\mathbf{x}}$, respectively. Moreover, $N$ and $n$ are referred to as the {\it population size} and {\it sample size}, respectively \cite{baik2006eigenvalues}. While the population size is fixed and the sample size goes to infinity, the sample covariance matrix is a good approximation of the population covariance matrix. However, as both the population size and the sample size become large with a constant ratio $N/n=c$, the sample covariance matrix does not approximate the population covariance matrix anymore. Fortunately, the {\it l.s.d.} of the sample covariance matrix is still related to the population covariance matrix. Considering an $N\times n$ sample matrix $\mathbf{X}_N\in\mathbb{C}^{N\times n}$ composed of $n$ {\it i.i.d.} samples with zero mean and covariance matrix $\mathbf{I}_N$, the corresponding sample covariance matrix can also be regarded as the Gram matrix of $\frac{1}{\sqrt{n}}\mathbf{X}_N$, in which $\mathbf{X}_N$ has {\it i.i.d.} entries of zero mean and unit variance. The convergence of the {\it e.s.d.} of the Gram matrix is proved by Mar\v{c}enko and Pastur, thus the limiting {\it e.s.d.}, namely, the {\it l.s.d.}, is known as the {\it Mar\v{c}enko-Pastur law} \cite{marvcenko1967distribution}, which unfolds as follows.

\begin{theorem}\label{th:MP_law}
Consider an $N\times n$ random matrix $\mathbf{X}_N\in\mathbb{C}^{N\times n}$ with independent entries of zero mean and unit variance. As $N, n\rightarrow\infty$ with a constant ratio $N/n=c$, the {\it e.s.d.} of $\mathbf{M}_N = \frac{1}{n}\mathbf{X}_N\mathbf{X}_N^{H}$ converges weakly and almost surely towards a non-random distribution whose {\it p.d.f.}, i.e., $f^{\mathbf{M}}(x;c)$, is given by
\begin{equation}\label{eq:MP_law}
f^{\mathbf{M}}(x; c) =
\left\{
\begin{aligned}
&\frac{1}{2\pi x c}\sqrt{(x-a)(b-x)} &,{\rm if}\ c<1,\\
&\left(1-\frac{1}{c}\right)\delta(x)+\frac{1}{2\pi x c}\sqrt{(x-a)(b-x)} &,{\rm otherwise},
\end{aligned}
\right.
\end{equation}
where $a=(1-\sqrt{c})^2$, $b=(1+\sqrt{c})^2$, and $\delta(x)$ is the {\it Dirac function} such that $\delta(x) = 1_{\{0\}}(x)$, which equals $1$ if $x=0$ or $0$ otherwise.
\end{theorem}

\begin{figure}[!t]
\begin{center}
\psfrag{x}[cc][cc][.5][0]{$x$}
\psfrag{y}[cc][cc][.5][0]{\rm probability density, $f^{\mathbf{M}}(x;c)$}
\epsfxsize=0.5\textwidth \leavevmode
\epsffile{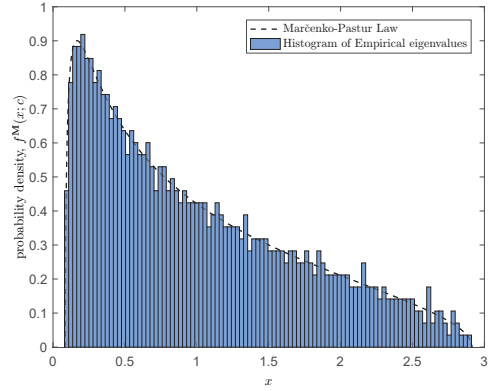}
\caption{Histogram of empirical eigenvalues and the Mar\v{c}enko-Pastur law when $c=0.5$, $N=1000$.}\label{fig:MP_law}
\end{center}
\end{figure}

\begin{figure}[!t]
\begin{center}
\psfrag{x}[cc][cc][.5][0]{$x$}
\psfrag{y}[cc][cc][.5][0]{\rm probability density, $f^{\mathbf{M}}(x;c)$}
\epsfxsize=0.5\textwidth \leavevmode
\epsffile{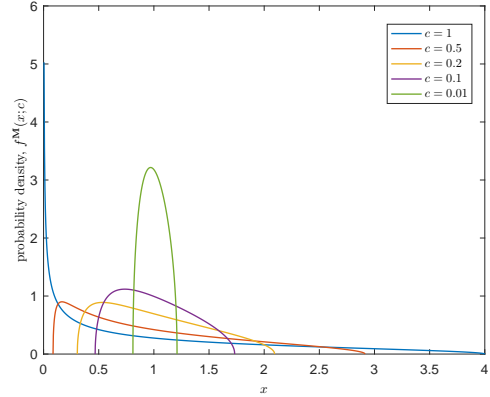}
\caption{Mar\v{c}enko-Pastur law for different $c$'s.}\label{fig:MP_law_c}
\end{center}
\end{figure}

The Mar\v{c}enko-Pastur distribution for $n = 2000$ and $N=1000$ is shown in Fig. \ref{fig:MP_law}. Besides, we also show the Mar\v{c}enko-Pastur distribution with different $c$'s in Fig. \ref{fig:MP_law_c}. It is worth noting that the entries of $\mathbf{X}$ are non-necessarily Gaussian in {\it Theorem} \ref{th:MP_law}. In addition, $\mathbf{M}_N$ is actually a sample covariance matrix $\hat{\mathbf{R}}_{\mathbf{x}\mathbf{x}}$ where $\mathbf{x}$ is a random process of zero mean with population covariance matrix $\mathbf{I}_N$. While $N$ is fixed, as $n\rightarrow\infty$, $c=N/n\to0$, the Mar\v{c}enko-Pastur distribution reduces to a single mass at $1$, this is consistent with the fact that the sample covariance matrix is an accurate approximation of the population covariance matrix in that case. Moreover, the Mar\v{c}enko-Pastur law also has a general form when the identical variance of the entries in $\mathbf{X}_N$ becomes $\sigma^2$, and the general Mar\v{c}enko-Pastur distribution with the additional parameter $\sigma$ is given by

\begin{equation}\label{eq:MP_law_scaled}
f_{MP}(x; c, \sigma)
\left\{
\begin{aligned}
&\frac{1}{2\pi x\sigma^2 c}\sqrt{(x-a_{\sigma})(b_{\sigma}-x)} &,{\rm if}\ c<1,\\
&\left(1-\frac{1}{c}\right)\delta(x)+\frac{1}{2\pi x\sigma^2 c}\sqrt{(x-a_{\sigma})(b_{\sigma}-x)}&,{\rm otherwise},
\end{aligned}
\right.
\end{equation}
\re{where $a_{\sigma}=\sigma^2(1-c)^2$ and $b_{\sigma}=\sigma^2(1+c)^2$}.

\subsection{\re{Stieltjes Transform and Free Probability Theory}}\label{subsec:Stieltjes_transform}
The Stieltjes transform is a powerful mathematical tool to prove many asymptotic results and conclusions in RMT. For example, the Mar\v{c}enko-Pastur law is exactly proved with the help of the Stieltjes transform \cite{couillet2011random}. To show the limiting results for more advanced random matrices, we first introduce the definition and some useful properties of the Stieltjes transform.
\begin{definition}
Let $F$ be a real-valued bounded measurable function over $\mathbb{R}$. The {\it Stieltjes transform} of $F$, denoted by $m_{F}(z)$, for $z\in Supp(F)^c$, is defined as
\begin{equation}\label{eq:Stieltjes_trans}
m_{F}(z)\triangleq\int_{-\infty}^{\infty}\frac{1}{\lambda-z}\dif F(\lambda),
\end{equation}
where $Supp(F)^c$ denotes the complex space complementary to the support of $F$; the support of $F$, i.e., $Supp(F)$, is the closure of the set $\{x\in\mathbb{R},f(x)>0\}$ and $f$ is the {\it p.d.f.} of $F$.
\end{definition}

Correspondingly, the {\it inverse Stieltjes transform} is defined as follows.
\begin{theorem}
If $x$ is a {\it continuity point} of $F$, then
\begin{equation}\label{eq:inverse_Stieltjes}
F(x)=\frac{1}{\pi}\lim\limits_{y\rightarrow 0^+}\int_{-\infty}^{x}\Im[m_F(x+iy)]\dif x,
\end{equation}
where the operator $\Im(\cdot)$ means to acquire the imaginary part.
\end{theorem}

The original intuition behind the Stieltjes transform is quite interesting and is illustrated as the following remark.
\begin{remark}
For a Hermitian random matrix $\mathbf{X}_N\in\mathbb{C}^{N\times N}$, the Stieltjes transform is given by

\begin{equation}\label{eq:intuition_Stieltjes}
\begin{aligned}
m_{F^{\mathbf{X}_N}}(z)&=\int\frac{1}{\lambda-z}\dif F^{\mathbf{X}_N}(\lambda)\\
&=\frac{1}{N}tr(\mathbf{\Lambda}-z\mathbf{I}_N)^{-1}\\
&=\frac{1}{N}tr(\mathbf{X}_N-z\mathbf{I}_N)^{-1},
\end{aligned}
\end{equation}
where $\mathbf{\Lambda}$ denotes the diagonal matrix consisting of the eigenvalues of $\mathbf{X}_N$. For notational simplicity, we also denote the Stieltjes transform of the {\it e.s.d.} of the Hermitian random matrix $\mathbf{X}_N$ by $m_{\mathbf{X}_N}\triangleq m_{F^{\mathbf{X}_N}}$ in the context. In \eqref{eq:intuition_Stieltjes}, it is observed that calculating the Stieltjes transform is equalent to working with the sum of diagonal entries of $(\mathbf{X}_N-z\mathbf{I}_N)^{-1}$. With the matrix inversion lemmas and some identities in matrix theory, it is quite simple to derive the limit of $tr(\mathbf{X}_N-z\mathbf{I}_N)^{-1}$. Thus, we can easily obtain a limit of {\it Stieltjes transform} of $F^{\mathbf{X}}$ as $N$ becomes large. The {\it l.s.d.} $F^{\mathbf{X}}$ such that $F^{\mathbf{X}_N}\stackrel{a.s.}{\Rightarrow} F^{\mathbf{X}}$ can be derived with the {\it inverse Stieltjes transform}. This is guaranteed by the following theorem \cite{bai2010spectral}.
\end{remark}

\begin{theorem}
Consider a set of bounded real functions $\{F_N\}$ satisfying $\lim\limits_{x\rightarrow-\infty}F_N(x)=0$. Then, $\forall z\in\mathbb{C}^{+}$
\begin{equation}\label{eq:limit_Stieltjes_transform}
\lim\limits_{N\rightarrow\infty}m_{F_{N}}(z) = m_{F}(z),
\end{equation}
{\it if and only if} there exists a function $F$ such that $\lim\limits_{x\rightarrow-\infty}F_N(x)=0$ and $|F_N(x)-F(x)|\to 0$ for all $x\in\mathbb{R}$.
\end{theorem}

An interesting identity between the Stieltjes transform of matrix $\mathbf{AB}$ and that of matrix $\mathbf{BA}$ when $\mathbf{AB}$ is Hermitian unfolds as follows.
\begin{corollary}
Let $\mathbf{A}\in\mathbb{C}^{N\times n}$, $\mathbf{B}\in\mathbb{C}^{n\times N}$, such that $\mathbf{AB}$ is Hermitian. For $z\in\mathbb{C}/\mathbb{R}$
\begin{equation}\label{eq:Stieltjes_indentity}
\frac{n}{N}m_{F^{\mathbf{BA}}}(z) = m_{F^{\mathbf{AB}}}(z)+\frac{N-n}{N}\frac{1}{z}.
\end{equation}
In particular, if $\mathbf{A}=\mathbf{B}^{H}=\mathbf{X}\in\mathbb{C}^{N\times n}$, for $z\in\mathbb{C}/\mathbb{R}^{+}$, \eqref{eq:Stieltjes_indentity} becomes
\begin{equation}\label{eq:particular_Stieltjes_indentity}
\frac{n}{N}m_{F^{\mathbf{X}^H\mathbf{X}}}(z) = m_{F^{\mathbf{X}\mathbf{X}^H}}(z)+\frac{N-n}{N}\frac{1}{z}.
\end{equation}
\end{corollary}

This identity is due to the fact that matrix $\mathbf{AB}$ and matrix $\mathbf{BA}$ have the same non-zero eigenvalues and different number of zero eigenvalues. Without loss of generality, assuming $n\geq N$, we denote the {\it p.d.f.} of the {\it e.s.d.} of $\mathbf{AB}$ and $\mathbf{BA}$ by $f^{\mathbf{AB}}(x)$ and $f^{\mathbf{BA}}(x)$, respectively, then we have
\begin{equation}\label{eq:pdf_identity}
\begin{aligned}
f^{\mathbf{BA}}(x)&=\frac{N}{n}f^{\mathbf{AB}}(\lambda)+\frac{n-N}{n}\delta(\lambda)\\
&=\left\{
\begin{aligned}
&\frac{N}{n}f^{\mathbf{AB}}(x)+\frac{n-N}{n}&, x=0;\\
&\frac{N}{n} f^{\mathbf{AB}}(x)&, x\neq0.
\end{aligned}
\right.
\end{aligned}
\end{equation}
With the Stieltjes transform, we finally obtain \eqref{eq:Stieltjes_indentity} via the following equation.
\begin{equation}\label{eq:proof_Stieltjes_transform}
\begin{aligned}
m_{F^{\mathbf{BA}}}(z) &= \int_{-\infty}^{\infty}\frac{1}{\lambda-z}dF^{\mathbf{BA}}(\lambda)\\
&=\int_{-\infty}^{\infty}\frac{1}{\lambda-z}f^{\mathbf{BA}}(\lambda)\dif \lambda\\
&=\int_{-\infty}^{\infty}\frac{1}{\lambda-z}\left(\frac{N}{n}f^{\mathbf{AB}}(\lambda)+\frac{n-N}{n}\delta(\lambda)\right)\dif \lambda\\
&=\frac{N}{n}m_{F^{\mathbf{AB}}}(z)-\frac{n-N}{n}\frac{1}{z}
\end{aligned}
\end{equation}

With the Stieltjes transform, we next introduce a kind of more complicated random matrices and the corresponding asymptotic results, which unfolds as the following theorem \cite{silverstein1995empirical}.
\begin{theorem}\label{th:advanced_model}
Let $\mathbf{B}_N=\mathbf{A}_N+\mathbf{X}_N^H\mathbf{T}_N\mathbf{X}_N$, where $\mathbf{X}_N\in \mathbb{C}^{N\times n}$ has independent entries with zero mean, variance $1/n$, and finite moment of order $2+\epsilon$ for some $\epsilon>0$ ($\epsilon$ is independent of $\mathbf{X}_N$), as $N$, $n$ grow large with a constant ratio $N/n=c$ $(0<c<\infty)$, $\mathbf{T}_N\in\mathbb{C}^{N\times N}$ is a diagonal matrix with real entries and its {\it e.s.d.} $F^{\mathbf{T}_N}$ converges weakly and almost surely to $F^{\mathbf{T}}$, $\mathbf{A}_N$ is a Hermitian matrix whose {\it e.s.d.} converges weakly and almost surely to $F^{\mathbf{A}}$. Then, the {\it e.s.d.} of $\mathbf{B}_N$, namely, $F^{\mathbf{B}_N}$ converges weakly and almost surely to a limit distribution $F^\mathbf{B}$ such that, for $z\in \mathbb{C}^{+}$, $m_{F^{\mathbf{B}}}(z)$ is the unique solution with positive imaginary part of
\begin{equation}\label{eq:Stieltjes_advanced}
m_{F^{\mathbf{B}}}(z)=m_{F^{\mathbf{A}}}\left(z-c\int\frac{t}{1+tm_{F^{\mathbf{B}}}(z)}\dif F^{\mathbf{T}}(t)\right).
\end{equation}
\end{theorem}
If the entries of $\mathbf{X}_N$ are identically distributed, \eqref{eq:Stieltjes_advanced} holds without requiring the finite moment of order $2+\epsilon$ \cite{couillet2011random}.

Under the particular case where $\mathbf{A}_N=0$, $\mathbf{B}_N$ reduces to a simpler form, i.e., $\mathbf{X}_N^H\mathbf{T}_N\mathbf{X}_N$. The matrix $\mathbf{T}_N^{\frac{1}{2}}\mathbf{X}_N\mathbf{X}_N^H\mathbf{T}_N^{\frac{1}{2}}$, where $\mathbf{T}_N^{\frac{1}{2}}$ denotes the Hermitian root of $\mathbf{T}_N$, can be regarded as the inverse Gram matrix of $\mathbf{X}_N^H\mathbf{T}_N\mathbf{X}_N$. To show the difference, we denote the {\it l.s.d.} of $\mathbf{T}_N^{\frac{1}{2}}\mathbf{X}_N\mathbf{X}_N^H\mathbf{T}_N^{\frac{1}{2}}$ and that of $\mathbf{X}_N^H\mathbf{T}_N\mathbf{X}_N$ by $F$ and $\underline{F}$, respectively. Besides, $\mathbf{T}_N^{\frac{1}{2}}\mathbf{X}_N\mathbf{X}_N^H\mathbf{T}_N^{\frac{1}{2}}$ is actually a general form of the sample covariance matrix while the population covariance matrix is $\mathbf{T}_N$. For example, the null Wishart matrix is a special case of $\mathbf{T}_N^{\frac{1}{2}}\mathbf{X}_N\mathbf{X}_N^H\mathbf{T}_N^{\frac{1}{2}}$ when $\mathbf{T}_N=\mathbf{I}_N$. With $\mathbf{A}_N=0$, \eqref{eq:Stieltjes_advanced} reduces to
\begin{equation}\label{eq:paticular_Stieltjes_advanced}
m_{\underline{F}}(z)=-\left(z-c\int\frac{t}{1+tm_{\underline{F}(z)}}\dif F^{\mathbf{T}}(t)\right)^{-1}.
\end{equation}
In addition, if we define $\mathbf{Y}_N = \mathbf{T}_N^{\frac{1}{2}}\mathbf{X}_N$, then we have $\mathbf{Y}_N\mathbf{Y}_N^H=\mathbf{T}_N^{\frac{1}{2}}\mathbf{X}_N\mathbf{X}_N^H\mathbf{T}_N^{\frac{1}{2}}$ and $\mathbf{Y}_N^H\mathbf{Y}_N=\mathbf{X}_N^H\mathbf{T}_N\mathbf{X}_N$. According to \eqref{eq:proof_Stieltjes_transform}, we can deduce the following equation:
\begin{equation}\label{eq:relation_sf}
m_{\underline{F}}(z) = cm_{F}(z)+(c-1)\frac{1}{z}.
\end{equation}
With \eqref{eq:relation_sf}, we can also obtain the Stieltjes transform of $F$ and therefore $F$ itself.


\re{
As we can see, the Stieltjes transform is a powerful tool to analyze the {\it l.s.d.} of complicated random matrix models. Besides, the free probability theory, which is closely related to the Stieltjes transform, is aimed to find an efficient approach to compute the spectrum of the products or the summations of the so-called {\it freely independent} matrices \cite{tao2012topics}. Here, we briefly introduce the key principles that play important roles in the free probability theory.
}

\re{
For a Hermitian random matrix $\mathbf{X}$, the Stieltjes transform $m_{F^{\mathbf{X}}}(z)$ can be obtained via \eqref{eq:intuition_Stieltjes}. Further, we define a moment generating function $M_{\mathbf{X}}$ as
\begin{equation}\label{eq:moment_generate_func}
M_{\mathbf{X}}\triangleq zm_{F^{\mathbf{X}}}(z)-1=\sum\limits_{k=1}^{\infty}\frac{m_k}{z^k},
\end{equation}
where $m_k=\int\lambda^k\dif F^{\mathbf{X}}(\lambda)$ is the $k$th moment of the {\it l.s.d.} of $\mathbf{X}$. Moreover, we denote the functional inverse of $M_{\mathbf{X}}$ by $M_{\mathbf{X}}^{-1}$ which obeys $M_{\mathbf{X}}(M_{\mathbf{X}}^{-1}(z))=M_{\mathbf{X}}^{-1}(M_{\mathbf{X}}(z))=z$. We finally define the {\it S-transform}, whose function is similar to that of the Stieltjes transform, as follows:
\begin{equation}\label{eq:stransform}
S_{\mathbf{X}}=\frac{1+z}{zM_{\mathbf{X}}^{-1}(z)}.
\end{equation}
The speciality of {\it S-transform} arises from its capability to deal with multiplications of random matrices. In particular, if two random matrix, e.g., $\mathbf{A}$ and $\mathbf{B}$, are {\it freely independent}, the {\it S-transform} of $\mathbf{AB}$ can be simply computed by
\begin{equation}\label{eq:stransformprod}
S_{\mathbf{AB}}=S_{\mathbf{A}}S_{\mathbf{B}}.
\end{equation}
}

\re{
Similarly, the R-transform is defined to compute the spectrum of the summation of {\it freely independent} matrices. For a Hermitian random matrix $\mathbf{X}$, the corresponding R-transform is given by
\begin{equation}\label{eq:R-transform}
R_{\mathbf{X}}(m_{F^{\mathbf{X}}}(z)) + \frac{1}{m_{F^{\mathbf{X}}}(z)} = z,
\end{equation}
where we recall that $m_{F^{\mathbf{X}}}(z)$ is the Stieltjes transform. For the {\it freely independent} random matrices, the R-transform of the summation of the random matrices is the summation of the R-transform of each random matrix. For example, if $\mathbf{A}$ and $\mathbf{B}$ are two freely independent random matrices, we have
\begin{equation}\label{eq:R-transform_sum}
R_{\mathbf{A}+\mathbf{B}} = R_{\mathbf{A}} + R_{\mathbf{B}}.
\end{equation}
}

\subsection{Characteristics of the Extreme Eigenvalues}\label{subsec:eigenvalues}
In the asymptotic regime, the {\it l.s.d.} of Wishart matrices and that of the Wigner matrices can be characterized by the Mar\v{c}enko-Pastur law and the semicircular law, respectively. It should be noted that the eigenvalues of a random matrix are actually a group of random variables. The limit spectrum distributions provide us the knowledge about the shapes of the spectra of the random matrices. However, the statistical characteristics of some specific eigenvalues are still unknown to us. For example, we may want to acquire the particular distributions of the extreme eigenvalues, i.e., the smallest and the largest eigenvalues. We may also want to know whether the extreme eigenvalues can be outside of the support of the limit spectra.
In \cite{bai1998no, yin1988limit}, it is shown that no eigenvalue can be found outside the support of the spectra for the general sample covariance matrices in terms of $\mathbf{T}_N^{\frac{1}{2}}\mathbf{X}_N\mathbf{X}_N^H\mathbf{T}_N^{\frac{1}{2}}$. This unfolds as the following theorem.

\begin{theorem}\label{th:no_eigenvaluesoutside}
Consider a matrix $\mathbf{X}_N\in\mathbf{C}^{N\times n}$ which has {\it i.i.d.} entries with zero mean, variance $ 1/n$, and finite fourth order moment, $\mathbf{T}_N\in\mathbb{C}^{N\times N}$ is a non-random matrix with uniformly bounded spectrum norm $\|\mathbf{T}_N\|$ and its {\it e.s.d.} $F^{\mathbf{T}_N}$ converges weakly and almost surely to a limit distribution function $H$. As shown in {\it Theorem} \ref{th:advanced_model}, the {\it e.s.d.} of $\mathbf{B}_N=\mathbf{T}_N^{\frac{1}{2}}\mathbf{X}_N\mathbf{X}_N^H\mathbf{T}_N^{\frac{1}{2}}\in\mathbb{C}^{N\times N}$ converges weakly and almost surely to a distribution function $F$ as $N, n\to\infty$ with $c_N=N/n\to c$ $(0<c<\infty)$. In addition, the {\it e.s.d.} of $\underline{\mathbf{B}}_N=\mathbf{X}_N^H\mathbf{T}_N\mathbf{X}_N$ converges weakly and almost surely towards $\underline{F}$ that satisfies
\begin{equation}\label{eq:relation_cdf}
\underline{F}(x) = cF(x)+(1-c)1_{[0,\infty]}(x)
\end{equation}
We denote $\underline{F}_N$ the distribution with Stieltjes transform $m_{\underline{F}_N}(z)$, which is the solution of the following equation of $m$ for $z\in \mathbb{C}^+$
\begin{equation}
m = -\left(z-\frac{N}{n}\int\frac{\tau}{1+\tau m}\dif F^{\mathbf{T}_N}(\tau)\right)^{-1}
\end{equation}
and define $F_N$ the distribution such that
\begin{equation}
\underline{F}_N=\frac{N}{n}F_N(x)+(1-\frac{N}{n})1_{[0, \infty)}(x).
\end{equation}
Let $N_0\in\mathbb{N}$, and choose an interval $[a, b]$ ($a,b\in(0, \infty])$ in an open interval outside the union of the supports of $F$ and $F_N$ for all $N\geq N_0$. For $\omega\in\Omega$, where $\Omega$ is the random space generating the series $\mathbf{X}_1, \mathbf{X}_2, \cdots$, denoting $\mathcal{L}_N(\omega)$ the set of eigenvalues of $\mathbf{B}_N(\omega)$, we have
\begin{equation}
P(\{\omega, \mathcal{L}_N(\omega)\cap[a,b]\neq\emptyset\  {\rm i.o.}\})=0,
\end{equation}
where ``i.o." means {\it infinitely often}.
\end{theorem}

{\it Theorem} \ref{th:no_eigenvaluesoutside} concretely means that, choosing an interval $[a, b]$ outside the union of the supports of $F$ and $F_N$ for all $N\geq N_0$, for all series $\mathbf{B}_1(\omega), \mathbf{B}_2(\omega), \cdots$, there exists $M(\omega)$ such that, for all $N\geq M(\omega)$, no eigenvalue of $\mathbf{B}_N(\omega)$ will appear in $[a, b]$. Besides, we define $F_K$ as the {\it l.s.d.} of $\mathbf{B}_N$ with $G=F^{\mathbf{T}_K}$. It is necessary to consider the supports of $F_N$ ($\forall N\geq N_0$) when only a few eigenvalues of $\mathbf{T}_N$ are isolated and finally contribute to $G$ with probability zero. Indeed, it is quite intuitive that, if the largest eigenvalue of $\mathbf{T}_N$ is much larger than the rest, at least one eigenvalue of $\mathbf{B}_N$ will also be larger than the rest (take $n\gg N$ to be convinced). This means that, if there no isolated eigenvalue in $\mathbf{T}_N$, no eigenvalue can be found outside the support of $F^{\mathbf{B}_N}$ as $N$ grows sufficiently large. The models in which $\mathbf{T}_N$ has isolated eigenvalues are referred to as the {\it spiked models}, which will be introduced later.

Now we consider the limiting statistical characteristics of the extreme eigenvalues, the main results on the limiting distributions of extreme eigenvalues originate from the work of Tracy and Widom \cite{tracy1996orthogonal}. The following results provide us the limit distributions of the extreme eigenvalues of Wigner matrices.
\begin{theorem}\label{th:TW_distribution}
Consider a Wigner matrix with independent Gaussian off-diagonal entries of zero mean and variance $\frac{1}{N}$ denoted by $\mathbf{X}_N\in\mathbb{C}^{N\times N}$, let $\lambda_{N}^{+}$, $\lambda_{N}^{-}$ denote the maximum eigenvalue and minimum eigenvalue of $\mathbf{X}_N$, respectively. Then, as $N\rightarrow\infty$, we have
\begin{eqnarray}
N^{\frac{2}{3}}(\lambda_{N}^{+}-2)\Rightarrow X^+\sim F_2 \label{eq:TW_distribution_max},\\
N^{\frac{2}{3}}(\lambda_{N}^{-}+2)\Rightarrow X^-\sim F^c_2 \label{eq:TW_distribution_min},
\end{eqnarray}
where $F_2$ is the {\it Tracy-Widom law} of order $2$ \cite{tracy2000distribution} given by
\begin{equation}\label{eq:TW_law}
F_2(t) = \exp\left(-\int_{t}^{\infty}(x-t)^2q^2(x)dx\right)
\end{equation}
with $q$ the {\it Painlev\'{e}} \uppercase\expandafter{\romannumeral2} function that solves the following differential equation
\begin{eqnarray}
q''(x) = xq(x)+2q^3(x)\label{eq:qx_1},\\
q(x)\sim{\rm Ai}(x)\  {\rm as}\  x\to+\infty,\label{eq:qx_2}
\end{eqnarray}
in which ${\rm Ai}(x)$ is the {\it Airy function} given by
\begin{equation}\label{eq:airy_func}
{\rm Ai}(x)=\frac{1}{2\pi}\int_{-\infty}^{\infty}e^{i\left(xt+\frac{t^3}{3}\right)}\dif t,
\end{equation}
and $F^c_2$ is defined as
\begin{equation}\label{eq:F_minus}
F^c_2(x)\triangleq 1-F_2(x).
\end{equation}
\end{theorem}

Besides, the random variables $\lambda_{N}^{+}$ and $\lambda_{N}^{-}$ are shown to be asymptotically independent \cite{bianchi2010asymptotic}. This thus provides us a way to study the asymptotic distribution of the condition number, i.e., $\lambda_{N}^{+}/\lambda_{N}^{-}$. The details unfold as the following theorem.
\begin{theorem}
With the assumptions in {\it Theorem} \ref{th:TW_distribution},
\begin{equation}\label{eq:asofconnum}
\left(N^{\frac{2}{3}}(\lambda_{N}^{+}-2), N^{\frac{2}{3}}(\lambda_{N}^{-}+2)\right)\Rightarrow (X^+,X^-)
\end{equation}
where $X^+$ and $X^-$ are independent random variables with distributions $F_2$, $F_2^c$, respectively. The random variable $\lambda_{N}^{+}/\lambda_{N}^{-}$ satisfies
\begin{equation}\label{eq:connumdistribution}
N^{\frac{2}{3}}\left(\frac{\lambda_{N}^{+}}{\lambda_{N}^{-}}+1\right)\Rightarrow-\frac{1}{2}(X^++X^-)
\end{equation}
\end{theorem}

The limiting distributions of extreme eigenvalues for the Wishart matrices in both real and complex cases are studied in \cite{johansson2000shape, feldheim2010universality}. The results are as follows.
\begin{theorem}\label{th:TW_law_wishart}
Let $\mathbf{X}_N\in\mathbb{C}^{N\times n}$ be a random matrix whose entries are {\it i.i.d.} zero-mean Gaussian variables with variance $1/n$. Denoting the largest and smallest eigenvalue of the Wishart matrix $\mathbf{X}_N\mathbf{X}_N^H$ by $\lambda_{N}^{+}$, $\lambda_{N}^{-}$, respectively. As $N,n\to\infty$ with $c=\lim N/n<1$, we have
\begin{eqnarray}
N^{\frac{2}{3}}\frac{\lambda_{N}^{+}-(1+\sqrt{c})^2}{(1+\sqrt{c})^{\frac{4}{3}}\sqrt{c}}\Rightarrow X\sim F_2,\label{eq:dfoflargestreigenvaluewishart}\\
N^{\frac{2}{3}}\frac{\lambda_{N}^{-}-(1-\sqrt{c})^2}{-(1-\sqrt{c})^{\frac{4}{3}}\sqrt{c}}\Rightarrow X\sim F_2,\label{eq:dfofsmallestreigenvaluewishart}
\end{eqnarray}
where $F_2$ is the {\it Tracy-Widom} distribution of order $2$ defined in \eqref{eq:TW_law}. In addition, the convergence result of $\lambda_{N}^{+}$ still holds for $c\geq 1$.
\end{theorem}

As we introduced in {\it Theorem} \ref{th:no_eigenvaluesoutside}, there are no eigenvalues outside the support of the {\it l.s.d.} of the Wishart matrix, i.e., the Mar\v{c}enko-Pastur distribution. With the assumptions and notations in {\it Theorem} \ref{th:TW_law_wishart}, the largest and the smallest eigenvalues converge to the edges of the support of the {\it l.s.d.} $F$ \cite{geman1980limit}. We recall that the edges of the Mar\v{c}enko-Pastur distribution are $a=(1-\sqrt{c})^2$, $b=(1+\sqrt{c})^2$. The limits of the two extreme eigenvalues are as follows \cite{geman1980limit, silverstein1985smallest}.
\begin{eqnarray}
\lambda_{N}^{+}\stackrel{a.s.}{\rightarrow}(1+\sqrt{c})^2, \label{eq:limitoflargesteigennullcase}\\
\lambda_{N}^{-}\stackrel{a.s.}{\rightarrow}(1-\sqrt{c})^2. \label{eq:limitofsmallesteigennullcase}
\end{eqnarray}

Note that \eqref{eq:dfoflargestreigenvaluewishart} in {\it Theorem} \ref{th:TW_law_wishart} has another form for real-valued random matrix $\mathbf{X}_N$, which is given as the following theorem \cite{johansson2000shape, johnstone2001distribution}.
\begin{theorem}\label{th:dfoflargestreigenvaluewishartreal}
Let $\mathbf{X}_N\in\mathbb{R}^{N\times n}$ be a random matrix whose entries are {\it i.i.d.} zero-mean Gaussian variables with variance $1/n$. Let $\mathbf{A}=n\mathbf{X}\mathbf{X}^H$, we denote the largest eigenvalue of $\mathbf{A}$ by $\lambda_{\max}(\mathbf{A})$.
Define two constant for centering and scaling as follows:
\begin{eqnarray}
  \mu_{n,N} &=& (\sqrt{n-1}+\sqrt{N})^2,\label{eq:centerconstant} \\
  \sigma_{n,N} &=& (\sqrt{n-1}+\sqrt{N})\left(\frac{1}{n-1}+\frac{1}{\sqrt{N}}\right)^{\frac{1}{3}}.\label{eq:scalingconstant}
\end{eqnarray}
As $N$, $n$ grow to infinity with $c=\lim_N\frac{N}{n}<1$, we have
\begin{equation}\label{eq:dfoflargestreigenvaluewishartreal}
\frac{\lambda_{\max}(\mathbf{A})-\mu_{n,N}}{\sigma_{n,N}}\to W_1\sim F_1,
\end{equation}
where $F_1$ is the {\it Tracy-Widom} law of order $1$ \cite{tracy2000distribution} given by
\begin{equation}\label{eq:TWorder1}
F_1(t)=\exp\left\{-\frac{1}{2}\int_{t}^{\infty}q(x)+(x-t)q^2(x)\dif x\right\}, t\in\mathbb{R},
\end{equation}
while $q(x)$ is the same with that defined in \eqref{eq:qx_1} and \eqref{eq:qx_2}.
\end{theorem}

\subsection{Spiked Models}\label{subsec:spiked_model}
We begin with a more detailed introduction to the aforementioned general sample covariance matrices. Let $\mathbf{T}_N$ be a fixed $N\times N$ non-negative definite Hermitian matrix. Let $\mathbf{X}_{N}\in\mathbb{C}^{N\times n}$ be a random matrix whose entries $\mathbf{X}_{N, ij}$ are {\it i.i.d.} complex random variables such that
\begin{equation}\label{eq:model_constraints}
\mathbb{E}(\mathbf{X}_{N,11})=0,\quad \mathbb{E}(|\mathbf{X}_{N,11}|^2)=1,\quad {\rm and}\quad \mathbb{E}(|\mathbf{X}_{N,11}|^4)<\infty.
\end{equation}

We use $\mathbf{B}_N=\frac{1}{n}\mathbf{T}_N^{\frac{1}{2}}\mathbf{X}_N\mathbf{X}_N^H\mathbf{T}_N^{\frac{1}{2}}$ to denote the sample covariance matrix, where $\mathbf{T}_N^{\frac{1}{2}}$ is a Hermitian square root of $\mathbf{T}_N$. Obviously, $\mathbf{T}_N$ is the population covariance matrix of the column vectors of $\mathbf{T}_N^{\frac{1}{2}}\mathbf{X}_N$. It is shown that this model covers various sample covariance matrices, since the population covariance matrices can be arbitrary. In addition, we denote the eigenvalues of $\mathbf{B}_N$ by $s_1^{(N)},s_2^{(N)},\cdots, s_N^{(N)}$. Thus, for some unitary matrix $\mathbf{U}_{\mathbf{B}}$, with the spectral decomposition (a.k.a. eigendecomposition) method, we have
\begin{equation}\label{eq:edofscm}
\mathbf{U}_{\mathbf{B}}\mathbf{B}_N\mathbf{U}_{\mathbf{B}}^{-1}=\left(
\begin{matrix}

  s_1^{(N)}&   &   & \\
  &   s_2^{(N)}&   & \\
  &   &   \ddots& \\
  &   &   &s_N^{(N)}
\end{matrix}\right)=diag(s_1^{(N)},s_2^{(N)},\cdots, s_N^{(N)}).
\end{equation}
For definiteness, we order the eigenvalues as $s_1^{(N)}\geq s_2^{(N)}\geq \cdots\geq s_N^{(N)}$.

Different from the typical null Wishart matrix, the so-called {\it spiked population model} proposed in \cite{johnstone2001distribution} allows some spikes, i.e., the eigenvalues not equal to 1, in the spectrum of the population covariance matrix $\mathbf{T}_N$. Without loss of generality, we assume that all the eigenvalues of $\mathbf{T}_N$ are $1$ except for the first $r$ eigenvalues. Let the first $r$ eigenvalues are $\alpha_1,\alpha_2,\cdots,\alpha_M$ with respective multiplicity $r_1, r_2,\cdots,r_M$, where $\alpha_1>\alpha_2>\cdots>\alpha_M$ are fixed real numbers for some $M\geq 0$ and $r_1, r_2,\cdots,r_M$ are fixed non-negative integers such that $r=r_1+r_2+\cdots+r_M$. Using the spectral decomposition again, for some unitary matrix $\mathbf{U}_{\mathbf{T}}$, we have
\begin{equation}\label{eq:edofpcm}
\mathbf{U}_{\mathbf{T}}\mathbf{T}_N\mathbf{U}_{\mathbf{T}}^{-1}=diag(\underbrace{\alpha_1,\cdots,\alpha_1}_{r_1},\underbrace{\alpha_2,\cdots,\alpha_2}_{r_2},\cdots,\underbrace{\alpha_M,\cdots,\alpha_M}_{r_M},\underbrace{1,\cdots,1}_{N-r}).
\end{equation}
Here, we set $r_0=0$ for definiteness. Obviously, the spiked model can be regarded as a finite-rank perturbation on the population covariance matrix of the null case \cite{bai2012sample}. In the context, we will use {\it sample eigenvalues} and {\it population eigenvalues} to represent the eigenvalues of the sample covariance matrix and that of the population covariance matrix, respectively.

The limiting laws of the sample eigenvalues of the spiked models unfold as the following theorem \cite{baik2006eigenvalues}.
\begin{theorem}\label{th:limitlawofspikedmodelcleq}
Assume $N, n\to\infty$ such that $N/n\to c$, where $c$ is a constant. Let $M_0$ be the number of $j$'s such that $\alpha_j>1+\sqrt{c}$, and let $M-M_1$ be the number of $j$'s such that $\alpha_j<1-\sqrt{c}$. Then we have the following results.
\begin{itemize}
  \item For $1\leq j\leq M_0$,
  \begin{equation}\label{eq:largelimitofspikecleq}
  s_{j,i}^{(N)}\triangleq s_{r_1+\cdots+r_{j-1}+i}^{(N)}\stackrel{a.s.}{\to}\phi(\alpha_j)=\alpha_j+\frac{c\alpha_j}{\alpha_j-1},\quad 1\leq i\leq r_j.
  \end{equation}
  \item The limits of the other sample eigenvalues depend on the value of $c$.

  -- If $c<1$, i.e., $N<n$, for $M_1+1\leq j\leq M$,
    \begin{equation}\label{eq:smalllimitofspikecleq}
    s_{j,i}^{(N)}\triangleq s_{N-r+r_1+\cdots+r_{j-1}+i}^{(N)}\stackrel{a.s.}{\to}\phi(\alpha_j)=\alpha_j+\frac{c\alpha_j}{\alpha_j-1},\quad 1\leq i\leq r_j.
    \end{equation}
    For the population eigenvalues inside $[1-\sqrt{c}, 1+\sqrt{c}]$, the following two sample eigenvalues satisfy
    \begin{equation}\label{eq:largelimitinsidecleq}
    s_{r_1+\cdots+r_{M_0}+1}^{(N)}\stackrel{a.s.}{\to}(1+\sqrt{c})^2,
    \end{equation}
    and
    \begin{equation}\label{eq:smalllimitinsidecleq}
    s_{N-r+r_1+\cdots+r_{M_1}}^{(N)}\stackrel{a.s.}{\to}(1-\sqrt{c})^2.
    \end{equation}
  -- If $c>1$, i.e., $N>n$, we have
    \begin{equation}\label{eq:largelimitinsidecgeq}
    s_{r_1+\cdots+r_{M_0}+1}^{(N)}\stackrel{a.s.}{\to}(1+\sqrt{c})^2,
    \end{equation}
    \begin{equation}\label{eq:smalllimitinsidecgeq}
    s_{n}^{(N)}\stackrel{a.s.}{\to}(1-\sqrt{c})^2,
    \end{equation}
    and
    \begin{equation}\label{eq:zeroeigenvaluescgeq}
    s_{n+1}^{(N)}=\cdots=s_{N}^{(N)}=0.
    \end{equation}
  -- If $c=1$, i.e., $N=n$, we have
    \begin{equation}\label{eq:largelimitinsideceq}
    s_{r_1+\cdots+r_{M_0}+1}^{(N)}\stackrel{a.s.}{\to}4,
    \end{equation}
    and
    \begin{equation}\label{eq:smalllimitinsideceq}
    s_{\min{\{n,N\}}}^{(N)}\stackrel{a.s.}{\to}0.
    \end{equation}
\end{itemize}
\end{theorem}

From {\it Theorem} \ref{th:limitlawofspikedmodelcleq}, we can see that, if all the non-unit population eigenvalues are {\it sufficiently close} to $1$ (i.e., $M_0=0$, $M_1=M$), the {\it l.s.d.} of the sample covariance matrix, namely, the Mar\v{c}enko-Pastur law is not disturbed and no sample eigenvalues have almost sure limits outside the support of the {\it l.s.d.}, i.e., $[(1-\sqrt{c})^2, (1+\sqrt{c})^2]$. Besides, the quantitative measure for evaluating whether the population eigenvalues are {\it sufficiently close} to $1$ turns to be whether the population eigenvalues are in the interval $[1-\sqrt{c}, 1+\sqrt{c}]$. More precisely, each population eigenvalue outside the interval $[1-\sqrt{c}, 1+\sqrt{c}]$ almost surely pulls one sample eigenvalue from the support $[(1-\sqrt{c})^2, (1+\sqrt{c})^2]$ of the {\it l.s.d.} of the null Wishart matrix and places it at $\alpha_j+\frac{c\alpha_j}{\alpha_j-1}$ in the limit.

In probability theory, there are two well-known theorems, namely, {\it law of large numbers} (LLR) and {\it central limit theorem} (CLT). The two theorems characterize a random variable by its limit and the fluctuation around the limit, respectively. {\it Theorem} \ref{th:limitlawofspikedmodelcleq} actually gives the limit of the extreme sample eigenvalues of the spiked models. In \cite{bai2008central}, the central limit theorems for the sample eigenvalues of the spiked models are studied. The conclusions unfold as follows.

We begin with a particular case of the spiked model in which $\mathbf{X}_N\in\mathbb{C}^{N\times n}$ has {\it i.i.d.} zero-mean entries with unit variance, and
\begin{equation}\label{eq:pcmclt}
\mathbf{T}_N = \left(
\begin{matrix}
  \mathbf{\Sigma} &  \\
   & \mathbf{I}_{p}
\end{matrix}\right),
\end{equation}
where $\mathbf{\Sigma}$ is a $r$ dimensional (non-necessarily diagonal) matrix and $p=N-r$. Hence, the $i$th column of $\mathbf{Y}=\mathbf{T}_N^{\frac{1}{2}}\mathbf{X}_N$ can be denoted by $\mathbf{y}_i=\mathbf{T}_N^{\frac{1}{2}}\mathbf{x}_i=(\mathbf{\xi}_i^T, \mathbf{\eta}_i^T)^T$ where $\mathbf{\xi}_i=[\mathbf{\xi}_i(1), \cdots, \mathbf{\xi}_i(r)]^T$, $\mathbf{\eta}_i=[\mathbf{\eta}_i(1), \cdots, \mathbf{\eta}_i(p)]^T$ are independent, of dimension $r$ and $p$, respectively. Obviously, $\xi_i$ is a random vector of zero mean and covariance matrix $\mathbf{\Sigma}$ while $\eta_i$ is a random vector of zero mean and covariance matrix $\mathbf{I}_p$. Thus, $\mathbf{S}_n=\frac{1}{n}\mathbf{Y}\mathbf{Y}^H$ is the sample covariance matrix of $\mathbf{y}_i$. Besides, we define $\mathbf{Y}_1 =\frac{1}{\sqrt{n}}\xi_{1:n}=\frac{1}{\sqrt{n}} [\xi_1, \cdots, \xi_n]$ and $\mathbf{Y}_2 =\frac{1}{\sqrt{n}}\eta_{1:n}=\frac{1}{\sqrt{n}} [\eta_1, \cdots, \eta_n]$, the sample covariance matrix is therefore
\begin{equation}\label{eq:scmbm}
\mathbf{S}_n=\frac{1}{n}\mathbf{Y}\mathbf{Y}^H=\left(
\begin{matrix}
  \mathbf{S}_{11} & \mathbf{S}_{12} \\
  \mathbf{S}_{21} & \mathbf{S}_{22}
\end{matrix}
\right)=\left(
\begin{matrix}
  \mathbf{Y}_1\mathbf{Y}_1^H & \mathbf{Y}_1\mathbf{Y}_2^H\\
  \mathbf{Y}_2\mathbf{Y}_1^H & \mathbf{Y}_2\mathbf{Y}_2^H
\end{matrix}
\right).
\end{equation}
Further, for $\lambda\notin[(1-\sqrt{c})^2, (1+\sqrt{c})^2]$, we define
\begin{equation}\label{eq:An}
\mathbf{A}_n=\mathbf{A}_n(\lambda)=\mathbf{Y}_2^H(\lambda\mathbf{I}-\mathbf{Y}_2\mathbf{Y}_2^H)^{-1}\mathbf{Y}_2,
\end{equation}
and
\begin{equation}\label{eq:Rn}
\mathbf{R}_n=\mathbf{R}_n(\lambda)=\frac{1}{\sqrt{n}}\{\xi_{1:n}(\mathbf{I}+\mathbf{A}_n)\xi_{1:n}^H-\mathbf{\Sigma} {\rm tr}(\mathbf{I}+\mathbf{A}_n)\}.
\end{equation}

Here, we consider the case where $c<1$ in {\it Theorem} \ref{th:limitlawofspikedmodelcleq}. Moreover, $K_j$ is used to denote the set of indexes of the sample eigenvalues outside the support of the Mar\v{c}enko-Pastur law due to the spike population eigenvalue $\alpha_j$. Obviously, for $j\in\{j|1\leq j\leq M_0 \ {\rm or}\ M_1+1\leq j\leq M\}$,
\begin{equation}\label{eq:kj}
K_j=\left\{
\begin{aligned}
  &\{r_1+\cdots+r_{j-1}+1, \cdots, r_1+\cdots+r_{j-1}+r_j\}&,\alpha_j>1+\sqrt{c}\\
  &\{N-r+r_1+\cdots+r_{j-1}+1,\cdots, N-r+r_1+\cdots+r_{j-1}+r_j\} &,\alpha_j<1-\sqrt{c}
\end{aligned}
\right.
\end{equation}
and the cardinality of $K_j$ is equal to $r_j$. Then, it is necessary to study the {\it central limit theorem} for the $r_j$-packed sample eigenvalues
\begin{equation}\label{eq:cltvariablesforalphaj}
\sqrt{n}[s_k^{(N)}-\phi(\alpha_j)],\quad k\in K_j,
\end{equation}
we recall that $\phi(\alpha_j)$ is the limit of $s_k^{(N)}$ ($k\in K_j$). Using the notations in \eqref{eq:largelimitofspikecleq} and \eqref{eq:smalllimitofspikecleq}, for each $\alpha_j$ outside $[1+\sqrt{c}, 1+\sqrt{c}]$, we consider the $r_j$ dimensional real vector $\sqrt{n}[s^{(N)}_{j, 1}-\phi(\alpha_j),\cdots, s^{(N)}_{j, r_j}-\phi(\alpha_j)]$. The {\it central limit theorem} for this vector is as follows.
\begin{theorem}\label{th:CLTspike}
For each $\alpha_j\notin[1+\sqrt{c}, 1+\sqrt{c}]$, the $r_j$ dimensional real vector $\sqrt{n}[s^{(N)}_{j, 1}-\phi(\alpha_j),\cdots, s^{(N)}_{j, r_j}-\phi(\alpha_j)]$ converges weakly to the distribution of the $r_j$ eigenvalues of the Gaussian random matrix
\begin{equation}\label{eq:CLTofsampleeigenvalues}
\frac{1}{1+cm_{3}[\phi(\alpha_j)]\alpha_j}\tilde{\mathbf{R}}_{jj}[\phi(\alpha_j)]
\end{equation}
where $\tilde{\mathbf{R}}_{jj}$ is the $j$-th diagonal block of $\tilde{\mathbf{R}}$ corresponding to indexes $\{u,v\in K_j\}$. $\tilde{\mathbf{R}}[\phi(\alpha_j)]=\mathbf{U}^H\mathbf{R}[\phi(\alpha_j)]\mathbf{U}$ where $\mathbf{U}$ is an unitary matrix such that
\begin{equation}\label{eq:decompositionofSigma}
\mathbf{\Sigma}=\mathbf{U}\left(
\begin{matrix}
  \alpha_1\mathbf{I}_{r_1}  &  &  \\
   & \ddots &  \\
   &  & \alpha_M\mathbf{I}_{r_M}
\end{matrix}\right)\mathbf{U}^H,
\end{equation}
and $m_3(\lambda)$ is defined as
\begin{equation}\label{eq:m3}
m_3(\lambda)=\int\frac{x}{(\lambda-x)^2}\dif F_{MP}(x;c),
\end{equation}
where $F_{MP}(x;c)$ is the {\it c.d.f.} of the Mar\v{c}enko-Pastur law parameterized by $c$.
\end{theorem}

{\it Theorem} \ref{th:CLTspike} shows that the limiting distribution of such $r_j$-packed sample eigenvalues are generally non-Gaussian and asymptotically dependent. However, if the multiplicity $r_j$ of the spike eigenvalue $\alpha_j$ equals to $1$, i.e., $\alpha_j$ is simple, then the corresponding sample eigenvalue is indeed Gaussian.

In \cite{bai2012sample}, the authors consider a generalized spiked model where the eigenvalues of the population covariance matrix of $\eta_i$ are non-necessarily equal to $1$ and derive the limiting laws of the sample eigenvalues and the {\it central limit theorem} for the packed sample eigenvalues. In addition, the block structure imposed in \eqref{eq:pcmclt} has been removed in \cite{zhang2020asymptotic}. Although the required mathematical tools are quite different, the obtained results and the conclusions are similar. Next, we consider the limiting behaviors of the extreme sample eigenvalues of the spiked models. The limiting distribution of the largest sample eigenvalue of the spiked model is given in the following theorem \cite{baik2005phase}.

\begin{theorem}\label{th:distributionoflargestspike}
Consider a particular spiked model where $\mathbf{X}_N\in\mathbb{C}^{N\times n}$ has {\it i.i.d.} Gaussian entries of zero mean and unit variance, and $\mathbf{T}_N=diag(\tau_1,\cdots, \tau_N)\in\mathbf{R}^{N\times N}$. Besides, for some fixed $r$ and $k$, $\tau_{r+1}=\cdots=\tau_{N}=1$ and $\tau_1=\cdots=\tau_k$ while $\tau_{k+1}, \cdots, \tau_r$ are in a compact subset of $(0, \tau_1)$. In the case $c=\lim N/n\to c<1$ as $N$, $n$ grow large, denoting the largest sample eigenvalue of $\frac{1}{n}\mathbf{T}_N^{\frac{1}{2}}\mathbf{X}_N\mathbf{X}_N^H\mathbf{T}_N^{\frac{1}{2}}$ by $\lambda_N^{+}$, we have:
\begin{itemize}
  \item if $\tau_1<1+\sqrt{c}$
  \begin{equation}\label{eq:limitofspikeleq}
    N^{\frac{2}{3}}\frac{\lambda_{N}^{+}-(1+\sqrt{c})^2}{(1+\sqrt{c})^{\frac{4}{3}}\sqrt{c}}\Rightarrow X\sim F_2,
  \end{equation}
  where $F_2$ is again the {\it Tracy-Widom} distribution defined in \eqref{eq:TW_law}.

  \item if $\tau_1>1+\sqrt{c}$
  \begin{equation}\label{eq:limitofspikegeq}
    \left(\tau_1^2-\frac{\tau_1^2c}{(\tau_1-1)^2}\right)^{\frac{1}{2}}n^{\frac{1}{2}}\left[\lambda_{N}^{+}-\left(\tau_1+\frac{\tau_1c}{\tau_1-1}\right)\right]\Rightarrow X_k\sim G_k,
  \end{equation}
  where $G_k$ is the distribution function of the largest eigenvalue of the $k\times k$ {\it Gaussian unitary ensemble} (GUE) \cite{couillet2011random}:
  \begin{equation}\label{eq:dfofleGUE}
    G_k(x)=\frac{1}{Z_k}\int_{-\infty}^{x}\cdots\int_{-\infty}^{x}\prod\limits_{1\leq i<j\leq k}|\xi_i-\xi_j|^2\prod\limits_{i=1}^{k}e^{-\frac{1}{2}\xi_i^2}\dif \xi_1\cdots\dif \xi_k,
  \end{equation}
  and $Z_k$ is a normalization constant; $\xi_1, \cdots, \xi_k$ denote the corresponding $k$ eigenvalues. In particular, $G_1(x)$ is the Gaussian distribution function, and this is consistent with the conclusion from {\it Theorem} \ref{th:CLTspike}.
\end{itemize}
\end{theorem}

With {\it Theorem} \ref{th:CLTspike} and\ref{th:distributionoflargestspike}, we can see that, if the largest population eigenvalue is not large enough to pull out a sample eigenvalue from the support of the Mar\v{c}enko-Pastur distribution, then the distribution of the largest sample eigenvalue is same with that in {\it Theorem} \ref{th:TW_law_wishart}. On the contrary, if the largest population eigenvalue exceeds the critical threshold, i.e., $1+\sqrt{c}$, the corresponding $k$-packed eigenvalues have a central limit. In particular, if $k=1$, the largest sample eigenvalue satisfies a Gaussian distribution which is provided in \cite{feral2009largest}. If we define
\begin{equation}\label{eq:expectationlargestspikedmodel}
    \mu(\lambda_N^{+})=\tau_1+\frac{c\tau_1}{\tau_1-1},
\end{equation}
\begin{equation}\label{eq:stddeviationlargestspikedmodel}
    v(\lambda_N^{+})=\tau_1\sqrt{1-\frac{c}{(\tau_1-1)^2}},
\end{equation}
Then the distribution of $\lambda_N^{+}$ can be described as
\begin{equation}\label{eq:distributionlargestspikedmodel}
n^{\frac{1}{2}}\frac{\lambda_N^{+}-\mu(\lambda_N^{+})}{v(\lambda_N^{+})}\sim \mathcal{N}(0, 1)
\end{equation}

\section{Large-Dimensional Random Matrix Theory in Cognitive Radio}\label{sec:spectrum_sensing}

Cognitive radio (CR) has been a hot topic in wireless communications in recent years since it substantially improves the spectrum efficiency via allowing secondary users to use spectrum that is licensed to the primary users. One of the basic principles in cognitive radio is that the secondary users should not affect the transmission of primary users. In the opportunistic CR, the secondary users are supposed to sense the state of the spectrum before launching data transmission. If the spectrum is detected to be occupied by the primary users, the secondary users should not start their transmission. On the contrary, the secondary users can exploit the vacant spectrum to transmit. In some sense, the performance of the designed spectrum sensing algorithms determines how much improvement can a CR system realize in terms of the overall spectrum efficiency. In this section, we will introduce the applications of RMT in designing the spectrum sensing methods.

\subsection{Basics of Spectrum Sensing}\label{subsec:sysmodelofsensing}

We consider a general scenario in cognitive radio, in which the secondary user (SU) is equipped with $N$ antennas and tries to sense the radio spectrum of its interest. The signal model here is actually same with that in cooperative sensing scenarios \cite{penna2009cooperative, cardoso2008cooperative}. The SU can obtain $n$ samples within the sensing interval, then make a decision on whether there exist active primary users (PUs). Thus, the sensing samples may come from one of the following two hypotheses:
\begin{itemize}
  \item {\bf Hypothesis 0} ($\mathcal{H} 0$): No active primary users exist in the vicinity of the SU. Hence, the sensing samples are actually drawn from additive white Gaussian noise (AWGN) process, the $i$-th sample is given by
     \begin{eqnarray}
     \mathbf{x}_i = \mathbf{u}_i, \label{eq:H0}
     \end{eqnarray}
     where $\mathbf{x}_i = [x_i(1),x_i(2),\cdots, x_i(N)]^{T}$, $i= 1, 2, \cdots, n$, $\mathbf{u}_i$ is the AWGN vector with zero mean and covariance matrix $\sigma_{u}^{2}\mathbf{I}_N$.

  \item {\bf Hypothesis 1} ($\mathcal{H} 1$): Without loss of generality, we assume that there are $K$ active PUs in the vicinity of the SU. Hence, the sensing samples which are composed of received signals from primary users and the noise vector, are denoted by
      \begin{eqnarray}
      \mathbf{x}_i = \mathbf{Hs}_i + \mathbf{u}_i. \label{eq:H1}
      \end{eqnarray}
      where the $N\times K$ matrix $\mathbf{H}$ denotes the channel from the $K$ primary users to the SU. $\mathbf{s}_i = [s_i(1),s_i(2),\cdots, s_i(K)]^{T}$ denotes the transmitted signals from the $K$ PUs (also can be regarded as a primary transmitter with $K$ antennas).
\end{itemize}

In addition, there are some mild assumptions which are usually considered in the literatures as follows:
\begin{description}
  \item[$\mathcal{AS} 1$:] Both the signal vector $\mathbf{s}_i$ and the noise vector $\mathbf{u}_i$ are independent temporally, and $\mathbf{s}_i$ is independent of $\mathbf{u}_i$.
  \item[$\mathcal{AS} 2$:] $\mathbf{s}_i$ is composed of {\it i.i.d.} Gaussian random variables of mean zero and variance $\sigma_{s}^{2}$.
\end{description}

Then, we concatenate the sensing samples to form a $N\times n$ dimensional observation matrix $\mathbf{X} = [\mathbf{x}_1, \mathbf{x}_2, \cdots, \mathbf{x}_n]$. Similarly, we define $\mathbf{S} = [\mathbf{s}_1, \mathbf{s}_2, \cdots, \mathbf{s}_n]$, $\mathbf{U} = [\mathbf{u}_1, \mathbf{u}_2, \cdots, \mathbf{u}_n]$. Hence, under the two hypotheses, we respectively have
\begin{eqnarray}
&\mathcal{H} 0:& \mathbf{X} = \mathbf{U}, \label{eq:obmatrixH0}\\
&\mathcal{H} 1:& \mathbf{X} = \mathbf{H}\mathbf{S} + \mathbf{U}. \label{eq:obmatrixH1}
\end{eqnarray}

Based on the observation matrix, we can construct various test statistics to solve this conventional signal detection problem \cite{liang2008sensing}. Obviously, there are two possible sensing results, namely, {\it PUs are absent} or {\it PUs are present}, which are usually denoted by $\mathcal{D}_0$ and $\mathcal{D}_1$, respectively. In particular, we are often interested in two probabilities, namely, {\it probability of detection} $P_d$ and {\it probability of false alarm} $P_{fa}$, which are respectively given by
\begin{eqnarray}
  P_d &=& P(\mathcal{D}_1|\mathcal{H}_1), \label{eq:pd} \\
  P_{fa} &=& P(\mathcal{D}_1|\mathcal{H}_0). \label{eq:pfa}
\end{eqnarray}
Higher $P_d$ can better protect the data transmissions of the PUs. However, this will cause higher $P_{fa}$ and higher $P_{fa}$ reduces chances of the SUs to access the idle channels, therefore degrades the throughput of the SUs. An obvious conclusion is that, if the PUs are perfectly protected, the SUs will be not permitted to access the channels that are allocated to the PUs. This is contrary to the original motivation of cognitive radio. Thus, in cognitive radio, we often consider how to maximize $P_d$ while keeping $P_{fa}$ below a certain threshold $\alpha_f$, i.e., $P_{fa}\leq \alpha_f$.

Since the primary signals and the additive noise are independent, it is observed that the expectation of the received power under the two hypotheses are quite different. Motivated by this, {\it Energy detection} (ED) is proposed and then becomes the most popular scheme in spectrum sensing. The test statistic of energy detection is give by
\begin{equation}\label{eq:tsofed}
T^{(ED)}(\mathbf{X})=\frac{1}{nN}\sum\limits_{i=1}^{n}\|\mathbf{x}_i\|^2.
\end{equation}
In general, there are two main steps when we design a sensing algorithm: the first step is to construct a test statistic, which is denoted by $T(\mathbf{X})$ in this paper, to test $\mathcal{H}_0$ against $\mathcal{H}_1$; The second step is to determine the detection threshold $\gamma$ for the designed test statistic, then we can declare that PUs are present when $T(\mathbf{X})>\gamma$ or say that PUs are absent otherwise \cite{zeng2010review}. For example, in energy detection, with accurate noise power $\sigma_u^2$, we can simply set
\begin{equation}\label{eq:EDgamma}
\gamma^{(ED)}=\sigma_u^2.
\end{equation}

According to the law of the large numbers, we can imagine that the threshold will work well in the regime where the number of samples is sufficiently large. However, due to the practically limited sensing time, the number of samples is therefore limited. We then need to set $\gamma$ with the knowledge of the distribution of $T(\mathbf{X})$ under $\mathcal{H}_0$ and a given tolerable false alarm probability. Without loss of generality, we here consider the real-valued case, i.e., both the noise and the signal are real random variables, according to the {\it central limit theorem}, $T(\mathbf{X})$ under $\mathcal{H}_0$ can be approximated by a Gaussian distribution given by
\begin{equation}\label{eq:EDdistributionofTA}
T^{(ED)}(\mathbf{X})\sim\mathcal{N}\left(\sigma_u^2, \frac{2\sigma_u^4}{nN}\right).
\end{equation}
Hence, for some given $P_{fa}$ and $n$, $\gamma$ is set as
\begin{equation}\label{eq:EDgamma1}
\gamma^{(ED)}=\sqrt{\frac{2}{nN}}Q^{-1}(P_{fa})+1,
\end{equation}
where
\begin{equation}\label{eq:Qfunction}
Q(t)=\frac{1}{\sqrt{2\pi}}\int_{t}^{+\infty}e^{-\frac{u^2}{2}}\dif u.
\end{equation}

It has been proved that energy detection is optimal for {\it i.i.d.} signal, i.e., under $\mathcal{AS} 2$ \cite{kay1993fundamentals}. The correlation of the signals will degrade its detection performance. In addition, the energy detection requires accurate noise power, namely, $\sigma_u^2$, to realize a good detection performance. In practice, the noise uncertainty problems usually exist due to the estimation errors of the noise power, further to incur a dramatic degradation of the detection performance \cite{cabric2006spectrum, zeng2009eigenvalue, sonnenschein1992radiometric, tandra2005fundamental}. Hence, several so-called blind spectrum sensing methods, which do not require the accurate estimate of the noise power, are proposed to overcome the noise uncertainty. Among them, the approaches based on eigenvalues of the sample covariance matrix achieve a notable performance.

\subsection{Sample Covariance Matrix under the Two Hypotheses}\label{subsec:SCMmodel}

The {\it sample covariance matrix} intrinsically indicates the existence of active PUs \cite{zhang2010multi}. Intuitively, this can be verified by the difference between the {\it population covariance matrices} under the two hypotheses: $\mathcal{H}_0$ and $\mathcal{H}_1$.
\begin{enumerate}
  \item {\it Pure Noise Case}: Under $\mathcal{H}_0$, the samples are actually {\it i.i.d.} Gaussian noise vectors. Thus, the sample covariance matrix can be expressed by a null Wishart matrix \cite{tulino2004random} with $n$ degrees of freedom and covariance matrix $\sigma_{u}^{2}{\bf I}_N$. In general, we denote the sample covariance matrix by $\hat{R}_{\mathbf{x}}$, thus, the sample covariance matrix under $\mathcal{H}_0$ is given as
    \begin{eqnarray}
    \hat{\mathbf{R}}_{\mathbf{x}\mathbf{x}} = \frac{1}{n}\mathbf{X}\mathbf{X}^{H} = \frac{1}{n}\mathbf{U}\mathbf{U}^{H} = \frac{1}{n}\sum\limits_{i=1}^{n}{\mathbf{u}_i\mathbf{u}^{H}_i}=\hat{\mathbf{R}}_{\mathbf{u}\mathbf{u}}. \label{eq:SCMH0}
    \end{eqnarray}
    We recall that when the number of the samples is sufficiently large, i.e., $n\rightarrow\infty$, the {\it sample covariance matrix} is a good approximation of the {\it population covariance matrix}. Thus, we have
    \begin{eqnarray}
    \hat{\mathbf{R}}_{\mathbf{u}\mathbf{u}}\to\mathbf{R}_{\mathbf{u}\mathbf{u}}=\mathbb{E}[\mathbf{u}_i\mathbf{u}^{H}_i]=\sigma_{u}^{2}{\bf I}_N. \label{eq:PCMH0}
    \end{eqnarray}
  \item {\it Signal-plus-Noise Case}: Under $\mathcal{H}_1$, the samples are composed of PUs' signals and the additive noise. With $\mathcal{AS} 1$ and $\mathcal{AS} 2$, the {\it population covariance matrix} can be written as
\begin{equation}\label{eq:PCMH1}
\mathbf{R}_{\mathbf{x}\mathbf{x}}=\mathbb{E}[\mathbf{x}_i\mathbf{x}^{H}_i]=\mathbb{E}[(\mathbf{Hs}_i + \mathbf{u}_i)(\mathbf{Hs}_i + \mathbf{u}_i)^H]=\sigma_s^2\mathbf{H}\mathbf{H}^H+\sigma_{u}^{2}{\bf I}_N.
\end{equation}
The corresponding {\it sample covariance matrix} is given by \cite{zeng2009eigenvalue}
\begin{equation}\label{eq:SCMH1}
\begin{aligned}
\hat{\mathbf{R}}_{\mathbf{x}\mathbf{x}}=\frac{1}{n}\sum\limits_{i=1}^{n}\mathbf{x}_i\mathbf{x}^{H}_i&=\frac{1}{n}\sum\limits_{i=1}^{n}[(\mathbf{Hs}_i + \mathbf{u}_i)(\mathbf{Hs}_i + \mathbf{u}_i)^H]\\
&\approx \mathbf{H}\frac{1}{n}\sum\limits_{i=1}^{n}\mathbf{s}_i\mathbf{s}^{H}_i\mathbf{H}^H+\frac{1}{n}\sum\limits_{i=1}^{n}\mathbf{u}_i\mathbf{u}^{H}_i\\
&= \mathbf{H}\hat{\mathbf{R}}_{\mathbf{s}\mathbf{s}}\mathbf{H}^H+\hat{\mathbf{R}}_{\mathbf{u}\mathbf{u}}.
\end{aligned}
\end{equation}
\end{enumerate}

Obviously, the $N$ population eigenvalues of $\mathbf{R}_{\mathbf{u}\mathbf{u}}$ are identical and equal to $\sigma_u^2$. The $N$ population eigenvalue of $\mathbf{R}_{\mathbf{x}\mathbf{x}}$ under $\mathcal{H}_0$ are thus $\sigma_u^2, \cdots, \sigma_u^2$. On the contrary, denoting the $N$ population eigenvalues of $\sigma_s^2\mathbf{H}\mathbf{H}^H$ by $\rho_1>\cdots>\rho_N$, the $N$ population eigenvalues of $\mathbf{R}_{\mathbf{x}\mathbf{x}}$ under $\mathcal{H}_1$ are respectively $\sigma_u^2+\rho_1, \cdots, \sigma_u^2+\rho_N$, which are obviously different from that under $\mathcal{H}_0$. Therefore, the active primary users can be detected by computing the eigenvalues of the population covariance matrix. This is exactly the original motivation to develop the eigenvalue-based methods. However, we only have access to the {\it sample covariance matrix} in practice. Similar to the previous analysis of the energy detection method, the {\it sample covariance matrix} can not approximate the {\it population covariance matrix} well due to the limited amount of sensing samples. This can also be verified by the conclusions about the Wishart matrix, e.g., the Mar\v{c}enko-Pastur law. Thus, how to acquire the distribution of the test statistic becomes the main obstacle in designing the eigenvalue-based approaches. Fortunately, the results in Section \ref{sec:rmt} provide the relations between the sample eigenvalues and the population eigenvalues. Besides, the Tracy-Widom law for the extreme eigenvalues of the sample covariance matrices under the signal-plus-noise case is derived in \cite{zhang2020tracy}. Next, we will see that the aforementioned results from RMT can greatly help us develop the eigenvalue-based spectrum sensing approaches.

\subsection{Eigenvalue-based Spectrum Sensing}\label{subsec:eigensensing}

In \cite{zeng2008maximum}, a {\it maximum eigenvalue detection} (MED) method is proposed. Here, we use the notation $\lambda_N^{+}$ to denote the largest eigenvalue of the sample covariance matrix, the test statistic is given by
\begin{equation}\label{eq:tsofMED}
T^{(MED)}(\mathbf{X})=\frac{\lambda_N^{+}}{\sigma_u^2}.
\end{equation}
For the real-valued case under $\mathcal{H}_0$, the limit distribution of $\lambda_N^{+}$ can be obtained via {\it Theorem} \ref{th:dfoflargestreigenvaluewishartreal}. With a given $\alpha_f$, the detection threshold $\gamma$ of MED is given by
\begin{equation}\label{eq:MEDgamma1}
\gamma^{(MED)}=\frac{(\sqrt{n}+\sqrt{N})^2}{n}\left(1+\frac{(\sqrt{n}+\sqrt{N})^{-\frac{2}{3}}}{(nN)^{\frac{1}{6}}}F_1^{-1}(1-\alpha_f)\right).
\end{equation}
For the complex-valued case, we just need to modify \eqref{eq:MEDgamma1} by replacing $F_1$ with $F_2$.

It is worth noting that the MED method is not a blind sensing approach since it also relies a lot on the accuracy of the estimate of the noise power. To solve this problem, a proper substitute for the noise power is obviously required to design a fully blind sensing approach. Note that in \eqref{eq:PCMH1}, if the rank of $\mathbf{H}$ is less than $N$, the smallest eigenvalue of the population covariance matrix is exactly equal to $\sigma_u^2$. Using the smallest eigenvalue of the sample covariance matrix, which is denoted by $\lambda_N^-$, as the estimation of the noise power, we get the {\it condition number detection} (CND) (a.k.a. {\it maximum-minimum eigenvalue} (MME) detection) method proposed in \cite{zeng2009eigenvalue}. Furthermore, in energy detection, if we replace the noise power with $\lambda_N^-$, we get the {\it energy with minimum eigenvalue} (EME) detection method. Similarly, for the real-valued case, we describe the CND and EME methods as follows.
\begin{itemize}
  \item CND:
  \begin{eqnarray}
  &T^{(CND)}(\mathbf{X})=\frac{\lambda_N^{+}}{\lambda_N^{-}},\label{eq:tsofCND}\\
  &\gamma^{(CND)}=\frac{(\sqrt{n}+\sqrt{N})^2}{(\sqrt{n}-\sqrt{N})^2}\left(1+\frac{(\sqrt{n}+\sqrt{N})^{-\frac{2}{3}}}{(nN)^{\frac{1}{6}}}F_1^{-1}(1-\alpha_f)\right).\label{eq:CNDgamma1}
  \end{eqnarray}
  \item EME:
  \begin{eqnarray}
  &T^{(EME)}(\mathbf{X})=\frac{T^{(ED)}(\mathbf{X})}{\lambda_N^{-}},\label{eq:tsofEME}\\
  &\gamma^{(EME)}=\left(\sqrt{\frac{2}{nN}}Q^{-1}(\alpha_f)+1\right)\frac{n}{(\sqrt{n}-\sqrt{N})^2}\label{eq:EMEgamma1}.
  \end{eqnarray}
\end{itemize}
Again, to obtain the CND method for complex-valued case, we just need to modify \eqref{eq:CNDgamma1} by replacing $F_1$ with $F_2$. It should also be pointed out that the two methods above are obtained by substituting $\sigma_u^2$ with the limit of $\lambda_{N}^{-}$, i.e., \eqref{eq:limitofsmallesteigennullcase}, straightforwardly. The distribution of the test statistics are obtained through only the limiting distribution of the largest eigenvalue while the limiting distribution of $\lambda_{N}^{-}$ is actually not considered \cite{zeng2009eigenvalue}. As a consequence, this approximation inevitably induces inaccuracy in the above two methods. In the later research, the inaccuracy problem is solved and the exact distribution of the condition number of the sample covariance matrix is computed\cite{penna2009cooperative}. Specifically, the limiting distribution of the condition number is derived from the limiting distributions of the two extreme eigenvalues in {\it Theorem} \ref{th:TW_law_wishart} and a general method to compute the distributions of quotients of independent random variables in \cite{curtiss1941distribution}. The distribution of the condition number is referred as to the {\it Tracy-Widom-Curtiss} distribution \cite{zhang2011spectrum, zhang2016exact}. With the exact distribution of the condition number, the detection threshold can be calculated more accurately and the performance of the CND method is further improved.

Besides, \cite{zhang2010multi} proposes to perform the blind spectrum sensing with the ratio of the {\it arithmetic mean} (AM) to the {\it geometric mean} (GM) of the eigenvalues, which is derived from the {\it generalized likelihood ratio test} (GLRT) paradigm \cite{kay1993fundamentals}. This detection method is thus known as the {\it arithmetic to geometric mean} (AGM) method. Denoting the eigenvalues of the sample covariance matrix with $s_1^{(N)}\geq s_2^{(N)}\geq \cdots\geq s_N^{(N)}$, the test statistic of the AGM detection method is given by
\begin{eqnarray}
T^{(AGM)}(\mathbf{X})=\frac{\frac{1}{N}\sum\limits_{k=1}^{N}s_k^{(N)}}{\left(\prod\limits_{k=1}^{N}s_k^{(N)}\right)^{\frac{1}{N}}}.\label{eq:tsofAGM}
\end{eqnarray}
Since \eqref{eq:tsofAGM} is quite complex, it is intractable to compute the detection threshold analytically. Alternatively, the threshold can be computed by the Monte-Carlo method with a given $\alpha_f$. Furthermore, \cite{bouallegue2018blind} propose a new detection method that performs AGM detection with only extreme eigenvalues, i.e., {\it mean-to-square extreme eigenvalue} (MSEE), whose test statistic is described as follows:
\begin{equation}\label{eq:tsMSEE}
T^{(MSEE)}(\mathbf{X})=\frac{\frac{1}{2}(\lambda_N^++\lambda_N^-)}{\sqrt{\lambda_N^+\lambda_N^-}}.
\end{equation}
Note that the test statistic in \eqref{eq:tsMSEE} can be regarded as a function of that of the CND method, the detection threshold $\gamma^{(MSEE)}$ can therefore be obtained analytically via $\gamma^{(CND)}$:
\begin{equation}\label{eq:gammaMSEE}
\gamma^{(MSEE)}=G^{-1}\left(\frac{(\sqrt{n}+\sqrt{N})^2}{nN}\left[1+\frac{(\sqrt{n}+\sqrt{N})^{-\frac{2}{3}}}{(nN)^{\frac{1}{6}}} F_1^{-1}(1-\alpha_f)\right]\right),
\end{equation}
where $G(x)=2x^2-1+2x\sqrt{x^2-1}$. Moreover, there exist some other similar eigenvalue-based spectrum sensing algorithms, such as the methods based on simplified predicted eigenvalue threshold (SPET) \cite{hassan2013multiple}, maximum-eigenvalue-to-the-geometric-mean (MEGM) \cite{pillay2012blind}, etc.

The underlying mechanism of the eigenvalue-based methods can also be explained by the results from RMT \cite{cardoso2008cooperative}. When the primary users are absent, the sample covariance matrix is actually a Wishart matrix. With {\it Theorem} \ref{th:no_eigenvaluesoutside}, we know that the no eigenvalue can be found outside the support of the Mar\v{c}enko-Pastur distribution. On the contrary, when the primary users are present, the sample covariance matrix can be described with the spiked model where the primary signals perform a low-rank perturbation on the null Wishart matrix. There may exist eigenvalues outside the distribution of the Mar\v{c}enko-Pastur distribution due to the large spikes. Therefore, the eigenvalue-based methods are able to distinguish which kind of random matrices the sample covariance matrix belongs to. However, as shown in {\it Theorem} \ref{th:distributionoflargestspike}, if the power of the perturbation (proportional to the SNR of the primary signals) is not large enough, there will be no eigenvalues outside the support of the Mar\v{c}enko-Pastur distribution. As a consequence, the eigenvalue-based methods will fail to detect the primary signals in the low SNR regime. Besides, {\it Theorem} \ref{th:distributionoflargestspike} also provides us a way to deal with the low SNR case. We can increase the number of samples, i.e, $n$, to reduce the limit ratio $c = N/n$, further to separate the spiked sample eigenvalues outside the support of the Mar\v{c}enko-Pastur distribution. This conclusion can be verified by the simulation results in the literatures about the eigenvalue-based spectrum sensing methods.

One can notice that, the detection thresholds in the above spectrum sensing methods are obtained with the distribution of the test statistics under $\mathcal{H}_0$ and the given $\alpha_f$. However, the detection performance is rarely analyzed since the sample covariance matrix under $\mathcal{H}_1$ is usually intractable. Thanks to the advanced results of the spiked model, the detection performance (in terms of probability of detection, probability of miss detection, or the error exponent) of some sensing methods under the single primary user case can be evaluated analytically \cite{bianchi2009performance, bianchi2011performance, penna2009probability}. For example, the detection performance of the CND method is evaluated in \cite{penna2009probability}. With the general method to derive the distribution of the quotient of independent random variables from \cite{curtiss1941distribution}, the authors propose to exploit the asymptotic independence between the largest and the smallest sample eigenvalue to derive the distribution of the test statistic. The limiting distribution of the largest sample eigenvalue in the spiked model, i.e., \eqref{eq:distributionlargestspikedmodel}, and the limiting distribution of the smallest sample eigenvalue i.e., \eqref{eq:dfofsmallestreigenvaluewishart}, are used to compute the distribution of the test statistic under $\mathcal{H}_1$. With the detection threshold calculated before, the probability of miss detection can be computed accurately.

\section{Large-Dimensional Random Matrix Theory in Large Communication Systems}\label{sec:LinearMURx}

In this section, we focus on the large multiuser systems in wireless communications. To support the communications of multiple users simultaneously, the resource for each user must be orthogonal or almost-orthogonal in some domain that can usually be the frequency domain, the space domain, or the code domain. As a consequence, the  corresponding methods to realize multiple access are thus respectively known as frequency-division multiple access (FDMA), space-division multiple access (SDMA), and code-division multiple access (CDMA). In the context, we mainly consider the uplink multiuser communications under the SDMA case and CDMA case.

\subsection{A Brief Overview of Multiuser Receivers}\label{subsec:conventionalRx}

In the direct-sequence code-division multiple access (DS-CDMA) systems, the information symbols of different users are transmitted via different spreading codes (a.k.a. signature sequences). The degrees of freedom are thus provided in the code domain to support the multiuser communications. We consider a general scenario where the spreading codes of different users are {\it randomly} and {\it independently} chosen \cite{tse1999linear, verdu1999spectral, tse2000linear}. Assuming that the length of the spreading code is $N$ and the total number of users is $K$, the received signal at the base station (BS) in a symbol-synchronous CDMA system can be modeled as
\begin{equation}\label{eq:Rxsignal_multiuser}
\begin{aligned}
\mathbf{x} &= \sum\limits_{k=1}^{K} \mathbf{h}_ks_k + \mathbf{u}\\
& = \mathbf{Hs} + \mathbf{u},
\end{aligned}
\end{equation}
where $\mathbf{h}_k$ and $s_k$ respectively denote the spreading code and the transmitted symbol of user $k$; $\mathbf{u}$ denotes the additive Gaussian noise vector; $\mathbf{H} = [\mathbf{h}_1, \cdots, \mathbf{h}_K]$ denotes the concatenated spreading code matrix; $\mathbf{s} = [s_1, \cdots, s_K]^{T}$ denotes the symbol vector consisting of the transmitted symbols of all users. For the SDMA case, the degrees of freedom are provided in the space domain, i.e., via multiple antennas. Consider the scenario where the channels from different users to the base station are of independent Rayleigh fading, the received signal at the base station can be still modeled with \eqref{eq:Rxsignal_multiuser} \cite{liang2006block}. The only difference is that $\mathbf{h}_k$ here represents the single-input-multiple-output (SIMO) channel from user $k$ to the base station. In addition, $K$ and $N$ are referred to as the signal dimension and observation dimension \cite{bai2014spectral}, respectively.

One can imagine that, since the spreading codes (or channels) of different users are random and thus not perfectly orthogonal, the users' transmitted symbols are inevitably interfering with each other at the receiver.
\re{
Therefore, the multiuser receivers are proposed to recover the transmitted symbols of each user as accurate as possible. In particular, the multiuser receivers can be divided into two main categories, namely, linear multiuser receivers and non-linear multiuser receivers.} Before detailed descriptions for the multiuser receivers, we make the following mild assumptions.
\begin{description}
  \item[$\mathcal{AS} 1$:] The transmitted symbols of different users are independent zero-mean random variables. The average trasmit power of user $k$ is $\mathbb{E}[|s_k^2|]=p_k$, for $k = 1, \cdots, K$.
  \item[$\mathcal{AS} 2$:] The additive Gaussian noise vector $\mathbf{u}$ is zero mean with covariance matrix $\mathbb{E}[\mathbf{u}\mathbf{u}^H]=\sigma_u^2\mathbf{I}_N$. In addition, it is independent of the transmitted symbols of users.
\end{description}

For linear multiuser receivers, the signal recovery process can be expressed as
\begin{equation}\label{eq:linear_multiuser_receiver}
\hat{\mathbf{s}}=\mathbf{W}^H\mathbf{x}=\mathbf{W}^H\mathbf{Hs} + \mathbf{W}^H\mathbf{u},
\end{equation}
where $\hat{\mathbf{s}}$ denotes the estimate of the users' symbols; $\mathbf{W}$ is exactly the matrix form of the linear receivers. Note that $\mathbf{W}=[\mathbf{w}_1, \cdots, \mathbf{w}_K]$, $\mathbf{w}_k$ can be considered as an extractor for the transmitted symbol of user $k$. Thus, we have
\begin{equation}\label{eq:signalextractor1}
\hat{s}_k = \mathbf{w}_k^H\mathbf{x}.
\end{equation}
Substitute \eqref{eq:Rxsignal_multiuser} to \eqref{eq:signalextractor1}, the formula of $\hat{s}_k$ can be written as
\begin{equation}\label{eq:signalextractor2}
\hat{s}_k = \mathbf{w}_k^H\mathbf{h}_ks_k + \sum\limits_{j\neq k}\mathbf{w}_k^H\mathbf{h}_js_j + \mathbf{w}_k^H\mathbf{u}.
\end{equation}
The signal-to-interference-plus-noise ratio (SINR) of user $k$, namely, $\gamma_k$, is thus given by
\begin{equation}\label{eq:SINR}
\begin{aligned}
\gamma_k &= \frac{p_k|\mathbf{w}_k^H\mathbf{h}_k|^2}{\sum_{j\neq k}p_j|\mathbf{w}_k^H\mathbf{h}_j| + \sigma_u^2\|\mathbf{w}_k\|^2}\\
&=\frac{p_k|\mathbf{w}_k^H\mathbf{h}_k|^2}{\mathbf{w}_k^H(\mathbf{H}_k\mathbf{D}_k\mathbf{H}_k^H+\sigma_u^2\mathbf{I}_N)\mathbf{w}_k},
\end{aligned}
\end{equation}
where
\begin{equation}\label{eq:Hi}
\mathbf{H}_k = [\mathbf{h}_1, \cdots, \mathbf{h}_{k-1}, \mathbf{h}_{k+1}, \cdots, \mathbf{h}_{K}],
\end{equation}
and
\begin{equation}\label{eq:Hi}
\mathbf{D}_k = diag([p_1, \cdots, p_{k-1}, p_{k+1}, \cdots, p_{K}]).
\end{equation}

The most well-known linear multiuser receivers are the zero-forcing (ZF) receiver (a.k.a. the decorrelator), the maximum-ratio combining (MRC) receiver (a.k.a. the matched-filter receiver), and minimum mean-square-error (MMSE) receiver. We first introduce the basic principles of the three linear multiuser receivers.

\begin{itemize}
  \item MRC receiver: The MRC receiver aims to extract its intended signal without considering the interference from the other users. The signal extractor for user $k$ is designed as
      \begin{equation}\label{eq:MRCRx}
        \mathbf{w}_k^{(MRC)} = \frac{\mathbf{h}_k}{\|\mathbf{h}_k\|^2}.
      \end{equation}
      The denominator is to ensure that the signal estimate of user $k$ is unbiased. It can be observed that the MRC receiver is optimal in single-user system but suffers from the interference severely. Therefore, it can achieve near-optimal performance when the interference power from other users is negligible.
  \item ZF receiver: The ZF receiver is designed to null out the interference from the other users. In other words, $\mathbf{W}$ is supposed to make the matrix product, i.e., $\mathbf{W}^H\mathbf{H}$, be an identity matrix. Thus, the matrix form of the ZF receiver is exactly the Moore–Penrose pseudo-inverse of the channel matrix. When $N \geq K$, the ZF receiver can be expressed as
      \begin{equation}\label{eq:ZFRx}
        \mathbf{W}^{(ZF)} = \mathbf{H}(\mathbf{H}^H\mathbf{H})^{-1}.
      \end{equation}
       The ZF receiver performs well when the interference power from the other users are very strong, and this often happens in the near-far resistance scenario of the conventional CDMA systems. However, the ZF receiver often suffers from the noise enhancement.
  \item MMSE receiver: The MMSE receiver in the context, is actually the linear MMSE (LMMSE) receiver, which is the optimal linear receiver maximizing the output SINR since it is aimed to minimize the mean-square-error (MSE) between the extracted symbols and the transmitted symbols.
      \begin{align}\label{eq:MMSERx}
      \mathbf{w}_k^{(MMSE)}&=\arg\min\limits_{\mathbf{w}_k}\mathbb{E}[|\mathbf{w}_k^H\mathbf{x}-s_k|^2]\\
        \label{eq:MMSERx1}&= p_k(\mathbf{H}\mathbf{D}\mathbf{H}^H+\sigma_u^2\mathbf{I}_N)^{-1}\mathbf{h}_k\\ \label{eq:MMSERx2}&=\frac{p_k(\mathbf{H}_k\mathbf{D}_k\mathbf{H}_k^H+\sigma_u^2\mathbf{I}_N)^{-1}\mathbf{h}_k}{1+p_k\mathbf{h}_k^H(\mathbf{H}_k\mathbf{D}_k\mathbf{H}_k^H+\sigma_u^2\mathbf{I}_N)^{-1}\mathbf{h}_k}
      \end{align}
      where $\mathbf{D} = diag([p_1, \cdots, p_{K}])$. \eqref{eq:MMSERx2} is obtained via the matrix inversion lemma. Substitute \eqref{eq:MMSERx2} into \eqref{eq:SINR}, we can get the output SINR of the MMSE receiver for user $k$ as follows
      \begin{equation}\label{eq:SINR_MMSE}
        \gamma_k^{(MMSE)}=p_k\mathbf{h}_k^H(\mathbf{H}_k\mathbf{D}_k\mathbf{H}_k^H+\sigma_u^2\mathbf{I}_N)^{-1}\mathbf{h}_k.
      \end{equation}

\end{itemize}

\re{
There are also many nonlinear signal detection methods for the multiuser systems, namely, nonlinear multiuser receivers, such as sphere decoding \cite{vikalo2006sphere}, or the generalized decision feedback equalizer/receiver (GDFE) which is proved to be equivalent to the VBLAST receiver \cite{ginis2001relation}. In particular, sphere decoding, which belongs to the lattice search techniques \cite{damen2003maximum}, realizes a near maximum likelihood (ML) detection performance with a lower complexity than ML detector. However, the computational complexity and memory demand increase dramatically as the signal dimension grows. On the other hand, a block-iterative GDFE (BI-GDFE) was proposed for the signal detection in large multiple-input-multiple-output (MIMO) communications systems that are also known as massive MIMO systems nowadays \cite{liang2006block}. The underlying mechanism of BI-GDFE is to detect the transmitted symbols in an iterative manner. To be specific, we obtain $\hat{\mathbf{s}}$ by performing MMSE detection over $\mathbf{x}$ and then make hard decisions over $\hat{\mathbf{s}}$ to get $\bar{\mathbf{s}}$ at one iteration. The decisions made in the previous iteration are then utilized to cancel the multiuser interference and finally we get the signal estimate of the current iteration. The signal detection process of BI-GDFE at $l$-th iteration can be described as
\begin{equation}\label{eq:BIGDFE}
\hat{\mathbf{s}_l}=\mathbf{F}_l^H\mathbf{x}+\mathbf{D}_l\bar{\mathbf{s}}_{l-1},
\end{equation}
where $\hat{\mathbf{s}_l}$ is the signal estimate at $l$-th iteration, $\bar{\mathbf{s}}_{l-1}$ denotes the hard decision of signal estimate at $(l-1)$-th iteration, $\mathbf{F}_l$ and $\mathbf{D}_l$ respectively denote the feed-forward equalizer (FFE) and the feedback equalizer (FBE) at $l$-th iteration. Without loss of generality, we here assume that the average transmit power of each user is the same, i.e., $p_k=p, \forall k=1, \cdots, K$. Besides, the elements of $\bar{\mathbf{s}}_{l}$ are assumed to be {\it i.i.d.} variables with zero mean and variance $p$ \cite{chan2001class}. Moreover, the input-decision-correlation (IDC) coefficient at $l$-th iteration, namely, $\rho_l$, is defined with
\begin{equation}\label{eq:IDC}
\mathbb{E}[\mathbf{s}\bar{\mathbf{s}}_{l}^H]=\rho_lp\mathbf{I}_K.
\end{equation}
With $\mathbf{A}_l=diag(\mathbf{F}_l^H\mathbf{H})$\footnote{\re{$\mathbf{A}_l=diag(\mathbf{F}_l^H\mathbf{H})$ denotes the diagonal matrix whose diagonal elements are the same with that of $\mathbf{F}_l^H\mathbf{H}$.}}, the optimal FFE and FBE that maximize the output SINR at $l$-th iteration are given by
\begin{equation}\label{eq:FFE}
\mathbf{F}_l=\left[(1-\rho_{l-1}^2)\mathbf{H}\mathbf{H}^H+\frac{\sigma_u^2}{p}\mathbf{I}_N\right]^{-1}\mathbf{H},
\end{equation}
\begin{equation}\label{eq:FBE}
\mathbf{D}_l=\rho_{l-1}(\mathbf{A}_l-\mathbf{F}_l^H\mathbf{H}).
\end{equation}
The maximum SINR for the $k$-th symbol of $\mathbf{s}$ at $l$-th iteration is given by
\begin{equation}\label{eq:SINRBIGDFE}
\gamma_k^{(BI-GDFE)}(l) = \frac{|\mathbf{A}_l(k,k)|^2p}{\mathbf{R}_{\tilde{\mathbf{u}}_l}(k,k)},
\end{equation}
where $\mathbf{A}_l(k,k)$ denotes the $k$-th element of $k$-th column of $\mathbf{A}_l$, $\mathbf{R}_{\tilde{\mathbf{u}}_l}$ denotes the covariance matrix of the equivalent additive noise and is given by
\begin{equation}\label{eq:Rul}
\mathbf{R}_{\tilde{\mathbf{u}}_l} = \frac{p(1-\rho_{l-1}^2)}{\rho_{l-1}^2}\mathbf{D}_l\mathbf{D}_l^H+\sigma_u^2\mathbf{F}_l^H\mathbf{F}_l.
\end{equation}
With properly selected IDC coefficients \cite{liang2006block}, the detection performance of BI-GDFE in the high SNR regime can approach the single user matched filter bound after a few iterations in large MIMO systems.
}

\subsection{Asymptotic Performance Analysis via RMT}\label{subsec:Asymptotic_SINR}
In general, the performance of a receiver can be evaluated by its output SINR since higher SINR means higher achievable data rate or lower bit-error-rate (BER) in communication systems. Here, we mainly focus on the output SINR of the multiuser receivers in the asymptotic regime, where both the signal dimension and the observation dimension become infinitely large with a constant ratio, i.e., $K\to\infty$, $N\to\infty$ with $K/N\to c$. With the two assumptions in \ref{subsec:conventionalRx}, we here make another mild assumption on the random channels (or the random spreading codes).
\begin{description}
  \item[$\mathcal{AS} 3$:] The channel (or the spreading code) of user $k$ is expressed as
  \begin{equation}\label{eq:randomchannel}
    \mathbf{h}_k = \frac{1}{\sqrt{N}}[v_{1k}, \cdots, v_{Nk}]^{T}, \forall k=1,\cdots, K,
  \end{equation}
  where $v_{nk}$'s ($\forall n=1,\cdots, N$) are {\it i.i.d} random variables that satisfy $\mathbb{E}[v_{nk}]=0$ and $\mathbb{E}[|v_{nk}|^2]=1$.
\end{description}

We first give a sketch of ideas to analyze the limit SINR of the MRC receiver. Substitute \eqref{eq:MRCRx} into \eqref{eq:SINR}, we get output SINR of user $k$ as follows
\begin{equation}\label{eq:SINR_MRC}
  \gamma_k^{(MRC)}
  =\frac{p_k\|\mathbf{h}_k\|^4}{\sum_{j\neq k}p_j|\mathbf{h}_k^H\mathbf{h}_j|^2+\|\mathbf{h}_k\|^2\sigma_u^2}
\end{equation}
Since $N\to\infty$, with the law of large numbers, $\|\mathbf{h}_k\|^2$ and $\|\mathbf{h}_k\|^4$ almost surely converge to $1$. Thus, the limit SINR is mainly determined by the interference from the other users, i.e., the first term in the denominator. Note that $K\to\infty$ and
\begin{equation}\label{eq:spread_gain}
\mathbb{E}[|\mathbf{h}_k^H\mathbf{h}_j|^2]=\mathbb{E}[(\sum\limits_{n=1}^{N}v_{nk}v_{nj})(\sum\limits_{m=1}^{N}v_{mk}v_{mj})]
=\mathbb{E}[\sum\limits_{n=1}^{N}v_{nk}^2v_{nj}^2]=\frac{1}{N},
\end{equation}
using the law of large number law again, $\sum_{j\neq k}p_j|\mathbf{h}_k^H\mathbf{h}_j|^2$ almost surely converges to $\frac{1}{N}\sum_{j\neq k}p_j$. In addition, $p_1, p_2, \cdots, p_K$ also can be regarded as a series of samples from a distribution whose {\it c.d.f.} is denoted by $F(p)$, $\sum_{j\neq k}p_j$ almost converges to $(K-1)\int_{0}^{\infty}p\dif F(p)$, which is almost surely equivalent to $K\int_{0}^{\infty}p\dif F(p)$. Thus, we have
\begin{equation}\label{eq:limititerference}
  \sum_{j\neq k}p_j|\mathbf{h}_k^H\mathbf{h}_j|^2\to c\int_{0}^{\infty}p\dif F(p).
\end{equation}
Finally, we obtain the limit SINR of the MRC receiver for user $k$ as
\begin{equation}\label{eq:limitSINRMRC}
  \gamma_k^{(MRC)}\to\bar{\gamma}_k^{(MRC)}=\frac{p_k}{c\int_{0}^{\infty}p\dif F(p)+\sigma_u^2}.
\end{equation}
The above analysis for the limit SINR of user $k$ under MRC receiver case is quite intuitive. The rigorous proof can be found in \cite{tse1999linear}. Besides, the conclusion in \eqref{eq:limitSINRMRC} gives us some enlightenments about the asymptotic results when both the signal dimension and the observation dimension go to infinity with a constant ratio. \eqref{eq:spread_gain} can be seen as a processing gain in suppressing the interference from the other users, and the MRC receiver can reduce the interference power to $1/N$ of the original averagely. On the other hand, the total interference power grows when the total number of users increases. As a consequence, the SINR converges to a constant value as the number of signal dimension and the observation dimension go to the infinity simultaneously with a constant ratio.

Next, we describe the limit SINR of user $k$ under ZF receiver. The results are given in \cite{tse1999linear} as follows.
\begin{equation}\label{eq:limitSINRZF}
\gamma_k^{(ZF)}\to\bar{\gamma}_k^{(ZF)}=\left\{
\begin{aligned}
&\frac{p_k}{\sigma_u^2}(1-c), & c<1,\\
&0, & c\geq 1.
\end{aligned}
\right.
\end{equation}
The conclusion can be explained from a geometric perspective. The ZF receiver tries to extract the symbol of user $k$ by projecting $\mathbf{h}_k$ onto a subspace which is orthogonal to all the columns in $\mathbf{H}_k$. If $c\geq 1$, we can not find a such subspace due to $K\geq N$. Thus, the interference can not be nulled out and the SINR tends to zero. With $\mathcal{V}_k\triangleq(span(\{\mathbf{h}_1, \cdots, \mathbf{h}_{k-1}, \mathbf{h}_{k+1}, \mathbf{h}_K\}))^{\perp}$, we denote the projection of $\mathbf{h}_k$ onto $\mathcal{V}_k$ by $\mathbf{r}_k$, the SINR of user $k$ can be expressed as
\begin{equation}\label{eq:SINR_ZF}
\gamma_k^{(ZF)} = \frac{p_k}{\sigma_u^2}\mathbf{r}_k^H\mathbf{r}_k=\frac{p_k}{\sigma_u^2}\|\mathbf{r}_k\|^2.
\end{equation}
Using {\it Lemma} $4.2$ in \cite{tse1999linear}, we have $\|\mathbf{r}_k\|^2\to 1-c$, thus \eqref{eq:limitSINRZF} is obtained.

The method to obtain the limit SINR for the MMSE receiver is quite delicate. Denoting the spectrum decomposition of $\mathbf{H}_k\mathbf{D}_k\mathbf{H}_k^H+\sigma_u^2\mathbf{I}_N$ by $\mathbf{Q}_k^H\mathbf{\Lambda}_k\mathbf{Q}_k$, \eqref{eq:SINR_MMSE} can be further expressed as
\begin{equation}\label{eq:SINR_MMSE_decom}
\begin{aligned}
    \gamma_k^{(MMSE)}&=p_k\mathbf{h}_k^H(\mathbf{Q}_k^H\mathbf{\Lambda}_k\mathbf{Q}_k+\sigma_u^2\mathbf{I}_N)^{-1}\mathbf{h}_k\\
    &=p_k\mathbf{h}_k^H\mathbf{Q}_k^H(\mathbf{\Lambda}_k+\sigma_u^2\mathbf{I}_N)^{-1}\mathbf{Q}_k\mathbf{h}_k\\
    &=p_k(\mathbf{Q}_k\mathbf{h}_k)^H(\mathbf{\Lambda}_k+\sigma_u^2\mathbf{I}_N)^{-1}(\mathbf{Q}_k\mathbf{h}_k)
\end{aligned}
\end{equation}
where $\mathbf{\Lambda}=diag([\lambda_1+\sigma_u^2, \cdots, \lambda_N+\sigma_u^2])$ and $\lambda_1, \cdots, \lambda_N$ are the eigenvalues of $\mathbf{H}_k\mathbf{D}_k\mathbf{H}_k^H$; and $\mathbf{Q}_k$ is the corresponding unitary matrix whose columns are the eigenvectors of $\mathbf{H}_k\mathbf{D}_k\mathbf{H}_k^H+\sigma_u^2\mathbf{I}_N$. Let $\mathbf{y}_k\triangleq \mathbf{Q}_k\mathbf{h}_k$, then we obtain
\begin{equation}\label{eq:SINR_MMSE1}
\gamma_k^{(MMSE)}=\sum\limits_{n=1}^{N}\frac{[{y}_k(n)]^2p_k}{\lambda_n+\sigma_u^2}.
\end{equation}
Using {\it Lemma} $4.3$ in \cite{tse1999linear}, we have
\begin{equation}\label{eq:limitSINR_MMSE}
\gamma_k^{(MMSE)}\to\bar{\gamma}_k^{(MMSE)}=\int_{0}^{\infty}\frac{p_k}{\lambda+\sigma_u^2}\dif G(\lambda)=p_k\int_{0}^{\infty}\frac{1}{\lambda+\sigma_u^2}\dif G(\lambda),
\end{equation}
where $G(\lambda)$ denotes the {\it l.s.d.} of $\mathbf{H}_k\mathbf{D}_k\mathbf{H}_k^H$. We recall that the Stieltjes transform of a limit spectrum density function, e.g., $G(\lambda)$, is defined as
\begin{equation}\label{eq:mG}
  m_{G}(z)=\int_{0}^{\infty}\frac{1}{\lambda-z}\dif G(\lambda).
\end{equation}
Thus, we now know $\gamma_k^{(MMSE)}$ almost surely converges to $p_km_{G}(-\sigma_u^2)$. On the other hand, the Stieltjes transform of $G(\lambda)$ has been studied in {\it Theorem} \ref{th:advanced_model}, the conclusion in \eqref{eq:paticular_Stieltjes_advanced} can be straightforwardly exploited to obtain $m_{G}$ as follows
\begin{equation}\label{eq:mGadvancedmodel}
m_{G}(z) = -\left(z-c\int\frac{p}{1+pm_{G}(z)}\dif F(p)\right)^{-1}.
\end{equation}
Substitute $m_{G}(-\sigma_u^2)=\frac{\bar{\gamma}_k^{(MMSE)}}{p_k}$ into \eqref{eq:mGadvancedmodel}, we finally obtain a equivalent for the limit output SINR of the MMSE receiver as follows:
\begin{equation}\label{eq:eqforlimitSINR_MMSE}
\bar{\gamma}_k^{(MMSE)}=\frac{p_k}{\sigma_u^2+c\int_{0}^{\infty}\frac{p_kp\dif F(p)}{p_k+p\bar{\gamma}_k^{(MMSE)}}}, {\rm when}\ N,K\to\infty\ {\rm with}\ \frac{K}{N}\to c.
\end{equation}
Obviously, there does not exist explicit formula for the limit SINR under the MMSE receiver case. However, when the received power of each user is equal to $p$, a simpler form can be obtained by solving a quadratic equation. The solution is given by
\begin{equation}\label{eq:eqpowerlimitSINR_MMSE}
\gamma_k^{(MMSE)}\to\bar{\gamma}_k^{(MMSE)}=\frac{(1-c)p}{2\sigma_u^2}-\frac{1}{2}+\sqrt{\frac{(1-c)^2p^2}{4\sigma_u^2}+\frac{(1+c)p}{2\sigma_u^2}+\frac{1}{4}}, \forall k=1,\cdots,K.
\end{equation}

\re{
The limit SINR of MMSE receiver also provides an efficient way to study the limit SINR of BI-GDFE in the asymptotic regime. With \eqref{eq:SINRBIGDFE}, the output SINR of user $k$ can be rewritten as
\begin{align}\label{eq:SINRBIGDFEre}
\gamma_k^{(BI-GDFE)}(l) &= \frac{|\mathbf{f}_k^H\mathbf{h}_k|^2}{(1-\rho_l^2)\sum_{j\neq k}|\mathbf{f}_k^H\mathbf{h}_k|^2+\frac{\sigma_u^2}{p}\|\mathbf{f}_k\|^2}\\
\label{eq:SINRBIGDFEMMSE}&=\frac{1}{1-\rho_l^2}\frac{|\mathbf{f}_k^H\mathbf{h}_k|^2}{\sum_{j\neq k}|\mathbf{f}_k^H\mathbf{h}_k|^2+\frac{\sigma_u^2}{p(1-\rho_l^2)}\|\mathbf{f}_k\|^2}.
\end{align}
Compare \eqref{eq:FFE} and \eqref{eq:MMSERx1}, \eqref{eq:SINRBIGDFEMMSE} and \eqref{eq:SINR_MMSE}, we can easily observe that the second multiplication component of \eqref{eq:SINRBIGDFEMMSE} is equivalent to the output SINR of the linear MMSE receiver operating under where the identical receiver power of each user is $(1-\rho_l^2)p$. Therefore, considering that the limit SINR of the MMSE receiver, namely, $\bar{\gamma}_k^{(MMSE)}$ in \eqref{eq:eqpowerlimitSINR_MMSE}, is a function with respect to the receive power of each user, we can denote the function by
\begin{equation}\label{eq:eqpowerlimitSINR_MMSEasfunction}
\bar{\gamma}_k^{(MMSE)} = g(p)=\frac{(1-c)p}{2\sigma_u^2}-\frac{1}{2}+\sqrt{\frac{(1-c)^2p^2}{4\sigma_u^2}+\frac{(1+c)p}{2\sigma_u^2}+\frac{1}{4}}.
\end{equation}
Finally, the limit output SINR of BI-GDFE at $l$-th iteration is given by \cite{liang2006block}
\begin{equation}\label{eq:limitSINRBIGDFE}
\begin{aligned}
\bar{\gamma}_k^{(BI-GDFE)}(l)&=\frac{1}{1-\rho_l^2}g[(1-\rho_l^2)p]\\
&=\frac{1}{1-\rho_l^2}\left[\frac{(1-c)(1-\rho_l^2)p}{2\sigma_u^2}-\frac{1}{2}+\sqrt{\frac{(1-c)^2(1-\rho_l^2)^2p^2}{4\sigma_u^2}+\frac{(1+c)(1-\rho_l^2)p}{2\sigma_u^2}+\frac{1}{4}}\right].
\end{aligned}
\end{equation}
}

In essence, the output SINRs of the multiuser receivers are random variables with some specific distributions. The aforementioned analysis actually gives the limits of the random variables as $N,K\to\infty$ with $K/N\to c$. However, for finite $N$ and $K$, the details about the distributions of output SINRs of the multiuser receivers are not clarified. In \cite{tse2000linear, liang2007asymptotic}, the limit distributions of the output SINRs for multiuser receivers are studied. As a consequence, the output SINR of each particular user is asymptotically Gaussian for large $N$. The obtained results about the limit output SINR actually only give the mean of the Gaussian distributions. To show the details about the Gaussian distributions, we here need a further assumption that the random channels (or spreading codes) satisfy $\mathbb{E}[|v_{nk}|^8]<\infty$. This assumption can be relaxed to finite fourth-order moment, but the stronger assumption is made to simplify the proofs in \cite{tse2000linear}. Without loss of generality, we consider the asymptotic SINR distribution of user $1$.

For the ZF receiver when $N>K$, the SINR of user $1$ is
\begin{equation}\label{eq:SINR_ZF1}
\gamma_1 = \frac{p_1}{\sigma_u^2[\mathbf{H}^{H}\mathbf{H}]_{1,1}}.
\end{equation}
The analysis of the fluctuations around the limit SINR starts from finding a equivalent but more useful form of \eqref{eq:SINR_ZF1}. According the introduction of the multiuser receivers in Section \ref{subsec:conventionalRx}, the ZF receiver and MMSE receiver are identical in the large SNR regime, i.e.,
\begin{equation}\label{eq:ZF_SINR_approx}
\sigma_u^2\lim\limits_{\sigma_u^2\to 0}\gamma^{(ZF)}=\sigma_u^2\lim\limits_{\sigma_u^2\to 0}\gamma^{(MMSE)}
\end{equation}
Since the interference from the other users are fully nulled out in the ZF receiver, we can assume the received power of each of the other users is equal to $p$. Then, with \eqref{eq:SINR_MMSE} and \eqref{eq:SINR_ZF1}, we have
\begin{equation}\label{eq:ZF_SINR_approx1}
\frac{1}{[\mathbf{H}^{H}\mathbf{H}]_{1,1}}=\lim\limits_{\sigma_u^2\to 0}\sigma_u^2\mathbf{h}_1^H(p\mathbf{H}_1\mathbf{H}_1^H+\sigma_u^2\mathbf{I}_N)^{-1}\mathbf{h}_1.
\end{equation}
Denoting the spectrum of $p\mathbf{H}_1\mathbf{H}_1^H$ as $\mathbf{O}^H\mathbf{F}\mathbf{O}$, where $\mathbf{F}=diag([\lambda_1, \cdots, \lambda_N])$, we have
\begin{equation}\label{eq:ZF_SINR_approx3}
\begin{aligned}
\lim\limits_{\sigma_u^2\to 0}\sigma_u^2\mathbf{h}_1^H(p\mathbf{H}_1\mathbf{H}_1^H+\sigma_u^2\mathbf{I}_N)^{-1}\mathbf{h}_1
&=\lim\limits_{\sigma_u^2\to 0}\sigma_u^2\mathbf{h}_1^H\mathbf{O}^H(\mathbf{F}+\sigma_u^2\mathbf{I}_N)^{-1}\mathbf{O}\mathbf{h}_1\\
&=\mathbf{h}_1^H\mathbf{O}^H\mathbf{A}\mathbf{O}\mathbf{h}_1,
\end{aligned}
\end{equation}
where $\mathbf{A}=diag([0, \cdots, 0, 1,\cdots, 1])=\lim\limits_{\sigma_u^2\to 0}\sigma_u^2(\mathbf{F}+\sigma_u^2\mathbf{I}_N)^{-1}$ since
\begin{equation}\label{eq:limitA}
\lim\limits_{\sigma_u^2\to 0}\frac{\sigma_u^2}{\lambda_i+\sigma_u^2}=\left\{
\begin{aligned}
&1, & \lambda_i=0,\\
&0, & \lambda_i\neq0.
\end{aligned}
\right.
\end{equation}
The number of $1$'s in the diagonal of $\mathbf{A}$ is the number of zero eigenvalues of $\mathbf{H}_1\mathbf{H}_1^H$, i.e., $N-K+1$. Moreover, the {\it l.s.d} of $\mathbf{A}$ is given by
\begin{equation}\label{eq:lsdA}
f^{\mathbf{A}}(\lambda)=c\delta(\lambda) + (1-c)\delta(\lambda).
\end{equation}
Under the real-valued case where the transmitted symbols and random channels are real, the distribution of $\mathbf{h}_1^H\mathbf{O}^H\mathbf{A}\mathbf{O}\mathbf{h}_1$ can be obtained via the following result from RMT, which is proved in \cite{marvcenko1967distribution, tse2000linear}. We have
\begin{equation}\label{eq:lemmaGaussian}
\sqrt{N}\left[\mathbf{h}_1^H\mathbf{O}^H\mathbf{A}\mathbf{O}\mathbf{h}_1-\frac{1}{N}tr(\mathbf{A})\right]\overset{\mathcal{D}}{\to}\mathcal{N}(0, a)
\end{equation}
where
\begin{equation*}
\begin{aligned}
a &= 2\int\lambda^2f^{\mathbf{A}}(\lambda)\dif\lambda+(\mathbb{E}[|v_{11}]|^4-3)\left(\int\lambda f^{\mathbf{A}}(\lambda)\dif\lambda\right)^2\\
&=2(1-c)+(\mathbb{E}[|v_{11}]|^4-3)(1-c)^2,\ {\rm when}\ f^{\mathbf{A}}(\lambda)\ {\rm is\ given\ as}\ \eqref{eq:lsdA}.
\end{aligned}
\end{equation*}
Following this train of thought, in the large SNR regime, we have
\begin{equation}\label{eq:distributionZFSINR}
\lim\limits_{\sigma_u^2\to 0}\gamma_1=\lim\limits_{\sigma_u^2\to 0}\frac{p_1}{\sigma_u^2[\mathbf{H}^{H}\mathbf{H}]_{1,1}}=\frac{p_1}{\sigma_u^2}\mathbf{h}_1^H\mathbf{O}^H\mathbf{A}\mathbf{O}\mathbf{h}_1.
\end{equation}
Substitute \eqref{eq:limitSINRZF} into \eqref{eq:distributionZFSINR}, we then obtain the asymptotic Gaussian distribution under the case where $N\to\infty$ with $c<1$ as follows
\begin{equation}\label{eq:limitdistributionZFSINR}
\sqrt{N}\left(\gamma_1^{(ZF)}-\frac{p_1}{\sigma_u^2}(1-c)\right)\overset{\mathcal{D}}{\to}\mathcal{N}\left(0, \left(\frac{p_1}{\sigma_u^2}\right)^2a\right).
\end{equation}

The distribution of the output SINR for the MMSE receiver can be analyzed in a similar manner. In \cite{tse2000linear}, the special case where the received power of all the users are the same is considered and \eqref{eq:eqpowerlimitSINR_MMSE} gives the convergence point of the output SINR. Assuming that the received powers of all the users are equal to $p$, for user $1$, the output SINR of the MMSE receiver, i.e., \eqref{eq:SINR_MMSE}, becomes
\begin{equation}\label{eq:SINR_MMSE_eqpower}
\gamma_1^{(MMSE)}=p\mathbf{h}_1^H(p\mathbf{H}_1\mathbf{H}_1^H+\sigma_u^2\mathbf{I}_N)^{-1}\mathbf{h}_1.
\end{equation}
Denote the spectrum decomposition of $p\mathbf{H}_1\mathbf{H}_1^H$ by $\mathbf{O}^H\mathbf{F}\mathbf{O}$,
\begin{equation}\label{eq:SINR_MMSE_eqpower}
\gamma_1^{(MMSE)}=p\mathbf{h}_1^H\mathbf{O}^H(\mathbf{F}+\sigma_u^2\mathbf{I}_N)^{-1}\mathbf{O}\mathbf{h}_1
\end{equation}
Using {\it Lemma} $3.2$ in \cite{tse2000linear} and denoting the {\it l.s.d.} of $p\mathbf{H}_1\mathbf{H}_1^H$ by $G(\lambda)$, we have
\begin{equation}\label{eq:SINR_MMSE_eqpower_approx}
\gamma_1^{(MMSE)}\approx\frac{p}{N}tr(\mathbf{F}+\sigma_u^2\mathbf{I}_N)^{-1}=p\int\frac{1}{\lambda+\sigma_u^2}\dif G(\lambda)
\end{equation}
On the other hand, using the result in \eqref{eq:lemmaGaussian}, we have
\begin{equation}\label{eq:distributionMMSESINR}
\sqrt{N}\left(\gamma_1^{(MMSE)}-\frac{p}{N}tr(\mathbf{F}+\sigma_u^2\mathbf{I}_N)^{-1}\right)\overset{\mathcal{D}}{\to}\mathcal{N}(0, b),
\end{equation}
where
\begin{equation*}
b = 2\int\left(\frac{p}{\lambda+\sigma_u^2}\right)^2\dif G(\lambda) + (\mathbb{E}[|v_{11}]|^4-3)\left[\int\frac{p}{(\lambda+\sigma_u^2)}\dif G(\lambda)\right]^2.
\end{equation*}
Note that $$\int\frac{p}{(\lambda+\sigma_u^2)}\dif G(\lambda)=\bar{\gamma}_1^{(MMSE)}$$ is actually the limit SINR in \eqref{eq:eqpowerlimitSINR_MMSE} and $$\int\frac{p}{(\lambda+\sigma_u^2)^2}\dif G(\lambda)=-\frac{\dif \bar{\gamma}_1^{(MMSE)}}{\dif(\sigma_u^2)},$$ we finally get
\begin{equation*}
b = \frac{2\bar{\gamma}_1^{(MMSE)}(1+\bar{\gamma}_1^{(MMSE)})^2}{\frac{\sigma_u^2}{p}(1+\bar{\gamma}_1^{(MMSE)})^2+c} + (\mathbb{E}[|v_{11}]|^4-3)(\bar{\gamma}_1^{(MMSE)})^2.
\end{equation*}

It should be noted that the result in \eqref{eq:distributionMMSESINR} is also obtained under the real-valued case. For the complex-valued case, \cite{liang2007asymptotic} proves: the variance of the output SINR under the complex-valued case is half of that under the real-valued case. The proof exploits the fact that the suboptimal MMSE receiver becomes optimal when the users have the same received power and the results about the asymptotic SINR distribution for the suboptimal MMSE receiver.

\subsection{Massive Connectivity Scenario}\label{subsec:massive_connectivity}

In recent years, the massive machine type communication (mMTC, a.k.a. massive connectivity or massive access) has been regarded as a significant scenario in the future communication networks \cite{zhu2010exploiting, liu2018massive1, liu2018massive2, chen2020massive}. A representative application of massive connectivity is the cellular Internet of Things (IoT), which can be regarded as an extension of the conventional multiuser system. In massive connectivity, the data traffic of the devices is sporadic and only a quite small number of the devices are active in a coherence interval, and thus we just need to decode the messages of the active devices.

The biggest difference of massive connectivity from the conventional multiuser systems is that the number of potential devices is much more than the available degrees of freedom while the number of active devices is usually less than the available degrees of freedom. Similarly, the degrees of freedom in massive connectivity can be provided by either code domain or the space domain \cite{zhu2010exploiting}. As a promising technique in 5G and beyond, massive multiple-input-multiple-output (MIMO) is expected to be capable of supporting massive devices. Moreover, massive MIMO is found to be especially suitable for massive connectivity \cite{liu2018massive1}. Therefore, it is preferred that the degrees of freedom are provided by the large number of antennas at the BS. Considering a general massive connectivity scenario where the BS with $M$ antennas serves $N$ potential single-antenna devices, the received signal model is given by
\begin{equation}\label{eq:signal_model_massive_connect}
\mathbf{x} = \sum\limits_{n=1}^{N}\alpha_n\mathbf{h}_ns_n + \mathbf{u}=\sum\limits_{k\in\mathcal{K}}\mathbf{h}_ks_k + \mathbf{u},
\end{equation}
where $\alpha_n\in\{0, 1\}$ is a binary indicator to represent the activity of device $n$, i.e., $\alpha_n=1$ for device $n$ is active; $\mathbf{h}_n\in\mathbb{C}^{M\times 1}\sim\mathcal{CN}(0, \beta_n\mathbf{I}_M)$ denotes the channel of device $n$ and $\beta_n$ denotes the path loss of device $n$; $\mathcal{K}$ is the set of active devices and $K=|\mathcal{K}|$ denotes the cardinality of $\mathcal{K}$.

The signal detection in massive connectivity usually performs in a two-phase manner. In the first phase, the BS detects the activities of all the potential devices and estimates the channels from the active devices. In the second phase, the BS decodes the transmitted symbols of each active device using the channel state information (CSI) acquired in the previous phase. Here, we assume that the channels are estimated via the pilot sequences of length $L$ in the first phase. In massive connectivity, $L$ is usually much smaller than $N$ due to the limited pilot length. Hence, it is impossible to allocate orthogonal pilot sequences to all potential devices. In the context, we consider a non-orthogonal pilot allocation scheme where each device $n$ is allocated to a random pilot sequence consisting of {\it i.i.d.} random variables with zero mean and variance $1/L$. Besides, each device sends its pilot sequence synchronously in the first phase. Denoting the identical transmit power of the active devices by $\rho^{pilot}$, the total transmit energy of each active device is denoted by $\xi=L\rho^{pilot}$.

Here, we are also interested in the asymptotic regime where the $L,K,N\to\infty$ with their ratios converge to some fixed constants, i.e., $N/L\to \omega$ and $K/N\to \epsilon$ with $\omega, \epsilon\in(0,\infty)$ while the total transmit power remains unchanged. In the first phase, an MMSE-based approximate message passing (AMP) algorithm is proposed to detect the activities of the potential devices and estimate the channels in \cite{liu2018massive1}. Besides, it is shown that the activity detection is nearly perfect when the number of the antennas goes to infinity. However, the channel estimation can not be perfect due to the non-orthogonal pilot sequences. The estimated channel and the channel estimation error of an active device $k$ are denoted by $\hat{\mathbf{h}}_k$ and $\Delta\mathbf{h}_k=\mathbf{h}_k-\hat{\mathbf{h}}_k$, respectively. After the MMSE-based AMP algorithm converges, the covariance matrices of $\hat{\mathbf{h}}_k$ and $\Delta\mathbf{h}_k$, are respectively given by
\begin{equation}\label{eq:CM_estimated_channel}
\mathbf{R}_{\hat{\mathbf{h}}_k\hat{\mathbf{h}}_k} = v_k(M)\mathbf{I},
\end{equation}
\begin{equation}\label{eq:CM_estimated_channel_error}
\mathbf{R}_{\Delta\mathbf{h}_k\Delta\mathbf{h}_k} = \Delta v_k(M)\mathbf{I},
\end{equation}
where $v_k(M)$ and $\Delta v_k(M)$ converge to as $M\to\infty$
\begin{equation}\label{eq:estimated_channel_variance_limitM}
\lim\limits_{M\to\infty}v_k(M)=\frac{\beta_k^2}{\beta_k+\tau_{\infty}^2},
\end{equation}
\begin{equation}\label{eq:estimated_channel_error_variance_limitM}
\lim\limits_{M\to\infty}\Delta v_k(M)=\frac{\beta_k\tau_{\infty}^2}{\beta_k+\tau_{\infty}^2}.
\end{equation}
In \eqref{eq:estimated_channel_variance_limitM} and \eqref{eq:estimated_channel_error_variance_limitM}, $\tau_{\infty}^2$ is the fixed-point solution to the following state evolution of the AMP algorithm as $M\to\infty$:
\begin{equation}\label{eq:tau0}
\tau_0^2=\frac{\sigma_u^2}{\xi}+\omega\epsilon\mathbb{E}_{\beta}[\beta],
\end{equation}
\begin{equation}\label{eq:taut1}
\tau_{t+1}^2=\frac{\sigma_u^2}{\xi}+\omega\epsilon\mathbb{E}_{\beta}\left[\frac{\beta\tau_{t}^2}{\beta+\tau_{t}^2}\right], t\geq 0.
\end{equation}
According to in \cite[Theorem 1]{liu2018massive2}, in the high SNR regime where $\omega\epsilon<1$, i.e., $L>K$, the fixed-point solution to \eqref{eq:taut1} is unique and converges as follows
\begin{equation}\label{eq:tauinftyconvergence}
\tau_{\infty}^2\to \frac{\sigma_u^2}{\xi(1-\omega\epsilon)}.
\end{equation}
Then $v_k$ and $\Delta v_k$ can be approximated by
\begin{equation}\label{eq:estimated_channel_variance_approx}
v_k=\frac{\beta_k^2}{\beta_k+\frac{\sigma_u^2}{\xi(1-\omega\epsilon)}},
\end{equation}
and
\begin{equation}\label{eq:estimated_channel_error_variance_approx}
\Delta v_k=\frac{\beta_k\frac{\sigma_u^2}{\xi(1-\omega\epsilon)}}{\beta_k+\frac{\sigma_u^2}{\xi(1-\omega\epsilon)}}.
\end{equation}

In the second phase, the received signal at BS is given by
\begin{equation}\label{eq:received_signal_massive}
\mathbf{x} = \sum\limits_{n\in\mathcal{K}}\mathbf{h}_n\sqrt{\rho^{data}}s_n + \mathbf{u},
\end{equation}
where $s_n\sim\mathcal{CN}(0,1)$ denotes the transmit symbol of device $n$; $\rho^{data}$ denotes the identical transmit power of all the active devices; $\mathbf{u}\sim\mathcal{CN}(0, \sigma_u^2\mathbf{I})$ is the AWGN at BS. With the estimated CSI in the first phase, the multiuser receivers can be employed to decode the messages of the active devices. As introduced in Section \ref{subsec:conventionalRx}, we denote the linear signal extractor to recover the signal of device $k\in\mathcal{K}$ by $\mathbf{w}_k$, the estimate of $s_k$ is given by
\begin{equation}\label{eq:sk_hat_massive}
\begin{aligned}
\hat{s}_k &= \mathbf{w}_k^H\left(\sum\limits_{n\in\mathcal{K}}\mathbf{h}_n\sqrt{\rho^{data}}s_n + \mathbf{u}\right)\\
&=\mathbf{w}_k^H\hat{\mathbf{h}}_k\sqrt{\rho^{data}}s_k+\mathbf{w}_k^H\sum\limits_{j\in\mathcal{K}, j\neq k}\hat{\mathbf{h}}_j\sqrt{\rho^{data}}s_j + \mathbf{w}_k^H\sum\limits_{n\in\mathcal{K}}\Delta\mathbf{h}_n\sqrt{\rho^{data}}s_n + \mathbf{w}_k^H\mathbf{u}.
\end{aligned}
\end{equation}
In \eqref{eq:sk_hat_massive}, the BS regards the estimated channel $\hat{\mathbf{h}}_k$ as the real channel $\mathbf{h}_k$ and treats the term $\mathbf{w}_k^H\sum\limits_{n\in\mathcal{K}}\Delta\mathbf{h}_n\sqrt{\rho^{data}}s_n$ as another additional noise. The SINR for decoding $s_k$ is therefore
\begin{equation}\label{eq:outputSINRmassive}
\gamma_k = \frac{\rho^{data}|\mathbf{w}_k^H\hat{\mathbf{h}}_k|^2}{\rho^{data}\sum\limits_{j\in\mathcal{K},j\neq k}|\mathbf{w}_k^H\hat{\mathbf{h}}_j|^2+
\rho^{data}\|\mathbf{w}_k\|^2\sum\limits_{n\in\mathcal{K}}\frac{\beta_n\tau_{\infty}^2}{\beta_n+\tau_{\infty}^2}+\sigma_u^2\|\mathbf{w}_k\|^2
}.
\end{equation}

The statistics of the estimated channels and the errors have been shown in \eqref{eq:CM_estimated_channel} -- \eqref{eq:estimated_channel_error_variance_limitM}. Besides, the estimated channels are nearly Gaussian in the massive MIMO limit. Two multiuser receivers are considered here: the MRC receiver and the MMSE receiver, which are respectively given as
\begin{equation}\label{eq:MRC_Rx_massive}
\mathbf{w}_k^{MRC}=\hat{\mathbf{h}}_k,
\end{equation}
and
\begin{equation}\label{eq:MMSE_Rx_massive}
\mathbf{w}_k^{MMSE}=\left(\sum\limits_{n\in\mathcal{K}}\rho^{data}\hat{\mathbf{h}}_n\hat{\mathbf{h}}_n^H+
\sum\limits_{n\in\mathcal{K}}\frac{\rho^{data}\beta_n\tau_{\infty}^2}{\beta_n+\tau_{\infty}^2}\mathbf{I}+\sigma_u^2\mathbf{I}\right)^{-1}\hat{\mathbf{h}}_k.
\end{equation}
Now we return to the asymptotic regime where $K,L,M,N$ go to infinity with the constant ratios, i.e., $\omega$, $\epsilon$ and an additional ratio $c=K/M$ ($c\in(0,\infty)$), the limit output SINR of the two receivers are respectively given by \cite{liu2018massive2}
\begin{equation}\label{eq:MRCoutputSINRmassive_limit}
\gamma_k^{MRC}\to\bar{\gamma}_k^{MRC}=\frac{\beta_k^2}{c\mathbb{E}[\beta](\beta_k+\tau_{\infty}^2)}, \forall k,
\end{equation}
and
\begin{equation}\label{eq:MMSEoutputSINRmassive_limit}
\gamma_k^{MMSE}\to\bar{\gamma}_k^{MMSE}=\frac{\beta_k^2}{\beta_k+\tau_{\infty}^2}\Gamma, \forall k,
\end{equation}
where $\Gamma$ is the unique fixed-point solution of the following equation:
\begin{equation}\label{eq:Gamma}
\Gamma = \frac{1}{c\mathbb{E}\left[\frac{\beta^2}{\beta+\tau_{\infty}^2+\beta^2\Gamma}\right] + c\mathbb{E}\left[\frac{\beta\tau_{\infty}^2}{\beta+\tau_{\infty}^2}\right]}.
\end{equation}

The formulas in \eqref{eq:MRCoutputSINRmassive_limit} and \eqref{eq:MMSEoutputSINRmassive_limit} are more involved compared to that in Section \ref{subsec:Asymptotic_SINR} due to the considerations of the channel estimation errors. The proofs mainly exploit the mathematical methods in \cite{tse1999linear, wagner2012large} and the statistics of $\hat{\mathbf{h}}_k$ and $\Delta\mathbf{h}_k$. It is worth noting that the results in \eqref{eq:MRCoutputSINRmassive_limit} and \eqref{eq:MMSEoutputSINRmassive_limit} are the same with that in \eqref{eq:limitSINRMRC} and \eqref{eq:eqforlimitSINR_MMSE} if the channel estimation had been perfect, i.e., $\hat{\mathbf{h}}_k=\mathbf{h}_k$. In other words, \eqref{eq:MRCoutputSINRmassive_limit} and \eqref{eq:MMSEoutputSINRmassive_limit} extend the conclusions in \eqref{eq:limitSINRMRC} and \eqref{eq:eqforlimitSINR_MMSE} to a more general case where the channel estimation error for each active device is considered.

\section{Large-Dimensional Random Matrix Theory in Deep Learning}\label{sec:deep_learning}

Deep learning has shown its state-of-the-art performance in many fields such as computer vision, natural language processing, human games, etc \cite{bengio2013representation,lecun2015deep,buduma2017fundamentals,caterini2018deep}. In deep learning, the deep neural networks empower the machines to be capable of human-like behaviors \cite{krizhevsky2012imagenet, schmidhuber2015deep}. More and more advanced neural network architectures are proposed to improve the performance of deep learning in some particular learning tasks. However, the neural networks are usually regarded as black boxes with merely visible input-ports and output-ports since the neural networks and the datasets are too complex to understand due to their extremely large dimensions. It is therefore hard to answer the questions such as why the deep neural networks perform so well, and how to improve the learning speed of the neural networks. Despite that some empirical tricks can be exploited to tune the neural networks, rigorous theories from the mathematics are needed to further promote the development of deep learning. In this section, we introduce some preliminary explorations that try to explain the properties of the neural networks from the perspective of RMT.

\subsection{Preliminaries and Background of Neural Networks}\label{subsec:basicsdeeplearning}

The phrase, {\it neural networks}, is actually a generic term for the various neural networks that are designed for different specific learning tasks. The popular ones among them, such as the convolutional neural networks (CNNs) popularly used in computer vision \cite{kalchbrenner2014convolutional} and the recurrent neural networks (RNNs) widely used in time series prediction \cite{rumelhart1986learning, elman1990finding, mikolov2010recurrent, zaremba2014recurrent}, have attracted a lot of attention for their extraordinary performance in solving specific problems. In this section, we introduce the basics of the most fundamental neural networks composed of only fully-connected layers, i.e., deep neural networks (DNNs), which are also known as the multi-layer perceptrons (MLPs) \cite{gardner1998artificial}.

In general, the deep fully-connected neural networks are used to approximate the extremely complex nonlinear functions that represent the hidden relations between the inputs and outputs of the networks. Obviously, only employing the linear operations to construct the neural networks is not enough to realize the complex functions. There are also nonlinear operations in the neural networks, i.e., the activation functions. Here, we mainly focus on the feed-forward neural networks. In particular, we consider an $L$-layer feed-forward neural network of synaptic weights $\mathbf{W}^1, \cdots, \mathbf{W}^L$ with $L+1$ neural activity vectors $\mathbf{x}^0, \cdots, \mathbf{x}^L$. Denoting the number of neurons in layer $l$ by $N_l$, we have $\mathbf{x}^l\in\mathbb{R}^l$ and $\mathbf{W}^l\in \mathbb{R}^{N_l\times N_{l-1}}$, the feed-forward dynamics elicited by the input $\mathbf{x}^0$ is given by \cite{pennington2017resurrecting,poole2016exponential}
\begin{equation}\label{eq:activation}
\mathbf{x}^l=\phi(\mathbf{h}^l),
\end{equation}
\begin{equation}\label{eq:feedforward}
\mathbf{h}^l=\mathbf{W}^l\mathbf{x}^{l-1}+\mathbf{b}^l,\ l = 1, \cdots, L,
\end{equation}
where $\mathbf{b}^l$ is the bias vector and $\mathbf{h}^l$ denotes the inputs to neurons at layer $l$; $\phi(\cdot)$ is the component-wise nonlinear activation function that transforms the pre-activations $\mathbf{h}^l$ to the post-activations $\mathbf{x}^l$.

In the applications of RMT for physics, approximating the constituents with random variables has made vital progresses in understanding large complex systems. Analogously, we may gain some insights via approximating the large complex modern neural networks with random variables in the similar way. In addition, the random configurations are related to random feature and kernel methods and define the initial loss surface  \cite{pennington2017nonlinear}, which is the geometric representation of the loss function with respect to the weights. Hence, the literatures are usually interested in the general ensembles of random neural networks where both the synaptic weights and the biases are {\it i.i.d.} Gaussian random variables. The explorations for understanding the neural networks start from an abundance of relevant matrices that are of theoretical and practical interest. The most attractive matrices are the input-output Jacobian \cite{pennington2017resurrecting, pennington2018emergence, pastur2020random, pastur2020random1, ling2019spectrum}, the Hessian of the loss function with respect to the weights \cite{kawaguchi2016deep, granziol2020beyond, pennington2017geometry, granziol2019towards, choromanska2015loss, dauphin2014identifying}, and the data covariance matrices of each layer in the neural networks \cite{pennington2017nonlinear, adlam2019random, benigni2019eigenvalue, liao2018spectrum}. For example, the knowledge of the input-output Jacobian can help us improve the learning speed by properly setting weight initialization and choosing the nonlinear activation functions. The Hessian contains the information about the loss surface, thus, studying the Hessian may give us a explanation about why the deep learning performs so well in spite of the non-convex loss functions. The data covariance matrices provide us a insight about how spectra of the data covariance matrices propagate through the neural networks.
\re{Moreover, RMT can also be exploited to understand the training and generalization performance of neural networks by deriving the limit training error and generalization error \cite{couillet2016asymptotic,louart2018random,hastie2019surprises,mei2019generalization}, or performing spectral analysis over the relevant kernel matrices, e.g. conjugate kernel (CK) \cite{fan2020spectra}, neural tangent kernel (NTK) \cite{jacot2018neural}.}

\re{
Before introducing the numerous works on the random feed-forward neural networks, we here stress that there also exist a few researches which are related to some advanced neural networks, i.e., CNNs \cite{xiao2018dynamical}, RNNs \cite{couillet2016asymptotic,chen2018dynamical}, generative adversarial networks (GANs) \cite{seddik2020random}, etc. For example, the input-output Jacobian spectra of CNNs and RNNs are analyzed in \cite{xiao2018dynamical} and  \cite{couillet2016asymptotic}, respectively. Besides, \cite{couillet2016asymptotic} derives the limiting train error and generalization error of {\it linear echo state neural networks}, which are actually a class of RNNs. These works will be discussed detailedly in the following sections. Another notable work studying the GAN-data, i.e., \cite{seddik2020random}, proves that the deep learning representations of the data produced by GAN (a.k.a. GAN-data) behaves as Gaussian mixtures. In particular, GAN is composed of two neural networks, namely, generative network and discriminative network. The generative network tries to learn the mapping from a latent space to the true data distribution of interest, while the discriminative network distinguishes data produced by the generator from the true data distribution. \cite{seddik2020random} proposes to describe the deep learning representations of GAN-data with {\it concentrated vectors} \cite{ledoux2001concentration}, which can be obtained by applying successive {\it Lipschitz} operations \cite{tao2012topics} to Gaussian random vectors. The spectral behaviors (e.g., spectral distribution and dominant eigenvectors) of the covariance matrix of the deep learning representations of GAN-data can be analyzed via RMT, and are shown to be the same with that of Gaussian mixture model (GMM) with the same means and covariances in the asymptotic regime.
}

\subsection{Achieving Dynamical Isometry with the Knowledge of the Input-Output Jacobian}\label{subsec:Jacobian}

It is well-known that the weight initialization has a strong impact on the learning speed in the training stage of deep learning. For example, making the mean squared singular value of the network's input-output Jacobian be $\mathcal{O}(1)$, i.e., stay constant for different depths of the neural network, \re{can prevent the gradients from vanishing or exploding exponentially}. In addition, keeping the mean squared singular value of the network's input-output Jacobian close to $1$ means that the norm of a randomly chosen error vector can be preserved {\it on average} in the back-propagation process \cite{pennington2017resurrecting}. Further, ensuring that all the singular values of the input-output Jacobian are concentrated near $1$ can approximately preserve the norm of every error single error vector and dramatically speed up the learning process \cite{saxe2013exact}. This phenomenon is known as a property called {\it dynamical isometry}. However, how to achieve dynamical isometry in neural networks is still a problem that has attracted a lot of attention. It is preliminarily shown that the distribution of the singular values of the input-output Jacobian depends on the depth of the network, the weight initialization, and the choice of nonlinear activation functions \cite{pennington2017resurrecting, pennington2018emergence}. Hence, it is quite essential to study how to control the entire distribution of the singular values of input-output Jacobian in deep learning.

Without loss of generality, we consider an $L$-layer network of width $N$ where $N_l=N$ ($l=1,\cdots, L$) and $\mathbf{W}\in\mathbb{R}^{N\times N}$. Based on the model described in \eqref{eq:activation} and \eqref{eq:feedforward}, the network's input-output Jacobian $\mathbf{J}\in\mathbb{R}^{N\times N}$ is given by
\begin{equation}\label{eq:Jacobian}
\mathbf{J}=\frac{\partial \mathbf{x}^L}{\partial \mathbf{x}^0}=\prod\limits_{l=1}^{L}\mathbf{D}^l\mathbf{W}^l,
\end{equation}
where $\mathbf{D}^l$ is a diagonal matrix whose entries $D_{ij}^l=\phi'(h_i^l)\delta_{ij}$, $\delta_{ij}$ is the Kronecker delta function, which equals $1$ when $i=j$ or $0$ otherwise. The input-output Jacobian is closely related to the back-propagation process in which the output errors are propagated backward to update the weight matrix layer by layer. If the input-output Jacobian is well-conditioned, then all the weight layers are expected to be well-conditioned.

Here, we consider the random neural networks with randomly initialized weights and biases. The biases $b_i^{l}$ are {\it i.i.d.} Gaussian random variables with zero mean and variance $\sigma_b^2$. For the weight initialization, two random ensembles are assumed: i) {\it Gaussian weights} whose entries $W_{ij}^l$ are {\it i.i.d.} Gaussian random variables with zero mean and variance $\sigma_w^2/N$; ii) {\it orthogonal weights} that are drawn from a uniform distribution over the scaled orthogonal matrices satisfying $(\mathbf{W}^l)^T\mathbf{W}^l=\sigma_w^2\mathbf{I}$. While the mean squared singular value of the input-output Jacobian is set to $1$ by proper rescaling, two metrics of our main interest are the largest singular value $s_{max}$ of the input-output Jacobian $\mathbf{J}$ (or the largest eigenvalue $\lambda_{max}$ of $\mathbf{J}\mathbf{J}^T$) and the variance $\sigma_{\mathbf{J}\mathbf{J}^T}^2$ of the eigenvalues of $\mathbf{J}\mathbf{J}^T$. They quantify the behaviors of the squared singular values around $1$, and thus the conditioning of the input-output Jacobian. If $\lambda_{max}\gg 1$ and $\sigma_{\mathbf{J}\mathbf{J}^T}^2\gg 1$, the input-output Jacobian is ill-conditioned and this will yield the slow learning speed \cite{pennington2017resurrecting}.

We start from reviewing the signal propagation process in the neural networks. The random matrices $\mathbf{D}^l$ in \eqref{eq:Jacobian} depend on the empirical distributions of the pre-activations $h_i^l$ ($i=1,\cdots, N$) entering the nonlinear activation function $\phi(\cdot)$. The propagation of the empirical distributions of the pre-activations among different layers are studied in \cite{poole2016exponential, schoenholz2016deep}. In the large $N$ regime, \re{it is shown that the empirical distributions of the pre-activations} converge to independent Gaussian distributions with zero mean and indentical variance $q^l$, where $q^l$ \re{is independent over the index $i$} and is given in a recursion way as follows
\begin{equation}\label{eq:ql}
q^l = \sigma_w^2\int\mathcal{D}h\phi(\sqrt{q^{l-1}}h)^2+\sigma_b^2,
\end{equation}
where $q^0=\frac{1}{N}\sum_{i=1}^{N}(h_i^0)^2$, and $\mathcal{D}h=\frac{\dif h}{\sqrt{2\pi}}\exp(-\frac{h^2}{2})$. Besides, there exists a fixed point for \eqref{eq:ql} as follows
\begin{equation}\label{eq:ql_fp}
q^*=\sigma_w^2\int\mathcal{D}h\phi(\sqrt{q^{*}}h)^2+\sigma_b^2.
\end{equation}
Obviously, if we let $q^0=q^*$ by choosing a proper $\mathbf{h}^0$, the propagation actually starts from the fixed point, thus the distribution of $\mathbf{D}^l$ is independent of $l$. Intriguingly, \cite{poole2016exponential} shows that even if the propagation is not started from the fixed point, the empirical distribution will reach the fixed point after a few layers. Thus, we can reasonably assume $q^l=q^*$ in the deep networks.

In addition, there is another quantity, namely, the mean squared singular values of the matrix $\mathbf{D}\mathbf{W}$, which determines whether the gradients exponentially explode or vanish in the deep networks. It is defined as
\begin{equation}\label{eq:chi}
\chi = \frac{1}{N}tr[(\mathbf{D}\mathbf{W})^T\mathbf{D}\mathbf{W}]=\sigma_w^2\int\mathcal{D}h[\phi'(\sqrt{q^{*}}h)]^2.
\end{equation}
In particular, when $\chi > 1$, the back-propagated gradients to update the weights will explode exponentially. On the contrary, the gradients will vanish exponentially when $\chi < 1$. Thus, the so-called the criticality condition, i.e., $\chi=1$, ensures proper initializations without exploding or vanishing gradients.
\re{Either vanishing gradients or exploding gradients will result in the failure of training of deep neural networks. Thus, the analysis for the behaviors of $\lambda_{max}$ and $\sigma_{\mathbf{J}\mathbf{J}^T}^2$ only makes sense under the criticality condition.}
In the following, the behaviors of $\lambda_{max}$ and $\sigma_{\mathbf{J}\mathbf{J}^T}^2$ are studied while $\chi$ is kept to be around $1$.

Note that the input-output Jacobian is a product term of $\mathbf{D}^l$ and $\mathbf{W}^l$ ($l=1,\cdots, L$) in \eqref{eq:Jacobian}, free probability theory can possibly be utilized to compute the spectrum of the input-output Jacobian.
In \cite{pennington2017resurrecting, pennington2018emergence}, it is shown that the {\it S-transform} of $\mathbf{J}\mathbf{J}^T$ can be rewritten as using the free probability theory
\begin{equation}\label{eq:sjj}
S_{\mathbf{J}\mathbf{J}^T}=\prod\limits_{l=1}^{L}S_{\mathbf{W}^l(\mathbf{W}^l)^T}S_{\mathbf{D}^l(\mathbf{D}^l)^T}
=\prod\limits_{l=1}^{L}S_{\mathbf{W}^l(\mathbf{W}^l)^T}S_{(\mathbf{D}^l)^2}
=S_{\mathbf{W}\mathbf{W}^T}^LS_{\mathbf{D}^2}^L,
\end{equation}
which is derived using the fact that the weights $\mathbf{W}^l$ ($l=1,\cdots, L$) have identical distribution and $\mathbf{D}^l$ ($l=1,\cdots, L$) are of independently identical distribution due to the reasonable assumption, namely, $q^l=q^*$. Hence, \eqref{eq:sjj} provides us a useful method to compute the {\it l.s.d.} of $\mathbf{J}\mathbf{J}^T$ in the large $N$ regime: i) calculate the {\it l.s.d} of $\mathbf{W}\mathbf{W}^T$ and $\mathbf{D}^2$; ii) compute the corresponding Stieltjes transforms and S-transforms of $\mathbf{W}\mathbf{W}^T$ and $\mathbf{D}^2$ according to \eqref{eq:intuition_Stieltjes}, \eqref{eq:moment_generate_func}, and \eqref{eq:stransform}; iii) compute the S-transform of $\mathbf{J}\mathbf{J}^T$ via \eqref{eq:sjj}; iv) Convert the S-transform to the corresponding Stieltjes transform and finally obtain $f^{\mathbf{J}\mathbf{J}^T}(\lambda)$ using the inverse Stieltjes transform.

\begin{figure*}[!t]
\begin{center}
\psfrag{f}[cc][cc][.6][0]{$L$}
\psfrag{b}[cc][cc][.6][0]{$\lambda_{max}$}
\psfrag{c}[cc][cc][.6][0]{$L$}
\psfrag{d}[cc][cc][.6][0]{$\sigma_{\mathbf{J}\mathbf{J}^T}^2$}

\subfigure[$\lambda_{max}$]{%
\epsfxsize=0.4\textwidth \leavevmode
\epsffile{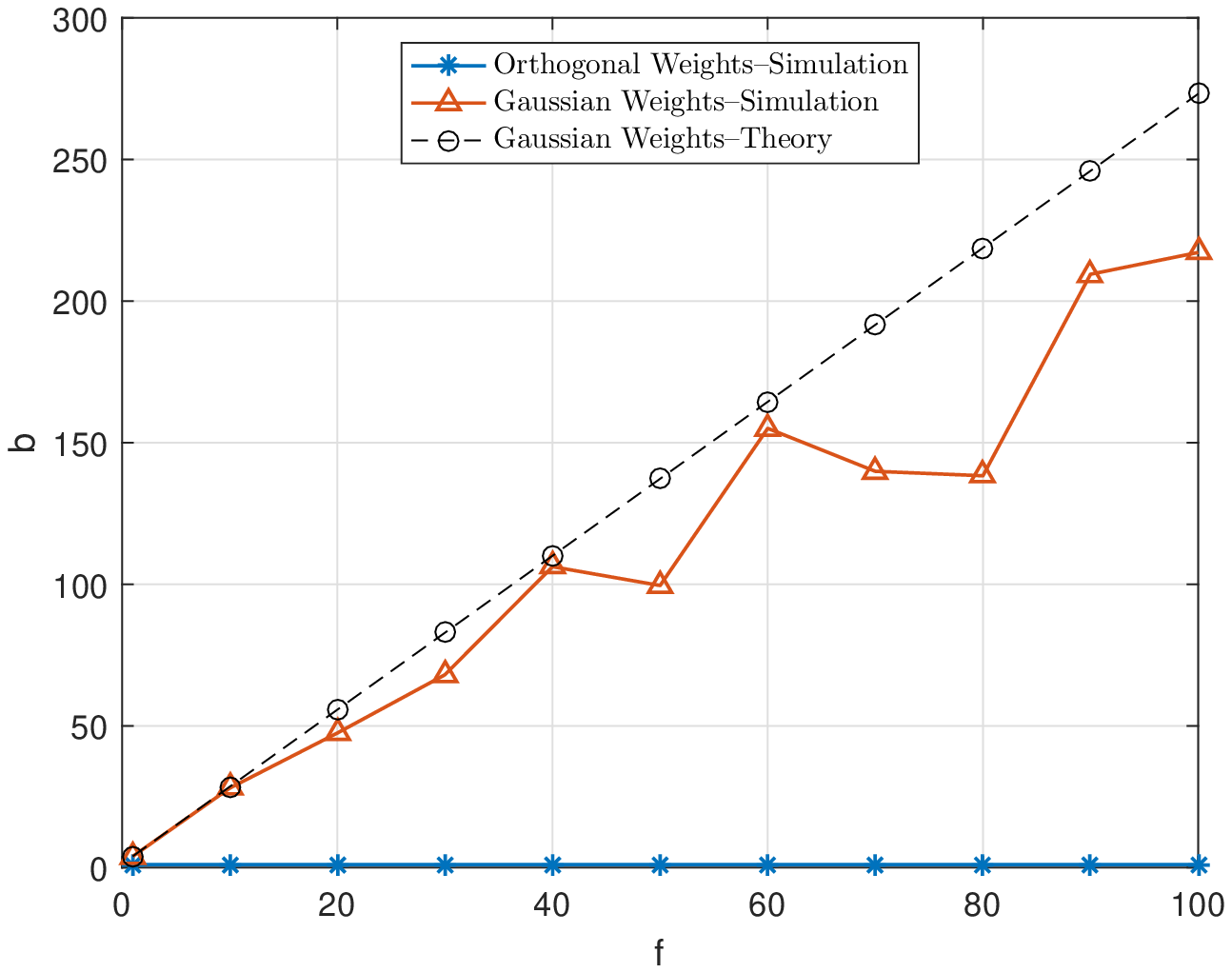}}\quad
\subfigure[$\sigma_{\mathbf{J}\mathbf{J}^T}^2$]{%
\epsfxsize=0.4\textwidth \leavevmode
\epsffile{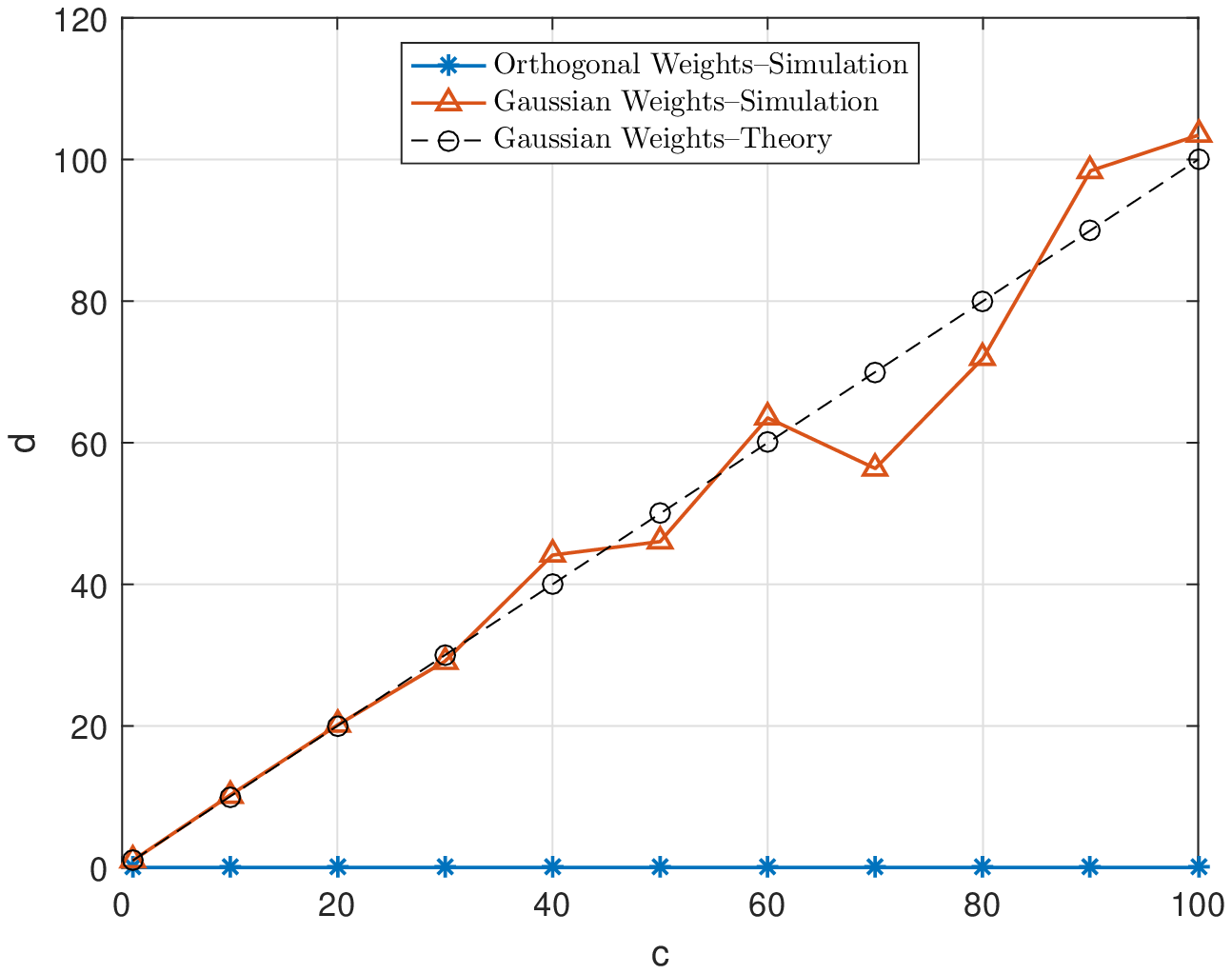}}

\caption{Linear growths of $\lambda_{max}$ and $\sigma_{\mathbf{J}\mathbf{J}^T}^2$ with respect to $L$ for the Gaussian random weights in linear neural networks. \re{The results are obtained on a single realization.}}\label{fig:valsta_linear_activation}
\end{center}
\end{figure*}

The computation of $f^{\mathbf{J}\mathbf{J}^T}(\lambda)$ is quite complex, we here omit the details and only present the results and corresponding conclusions. For the linear networks that have no nonlinear activation functions, the Jacobian $\mathbf{J}$ reduces to $\prod_{l=1}^{L}\mathbf{W}^l$. When the network is initialized with random orthogonal weights, all the singular values are $1$, and therefore realizing perfect dynamical isometry. For the Gaussian random weights, $\mathbf{J}\mathbf{J}^T=\prod_{l=1}^{L}\mathbf{W}^l(\mathbf{W}^l)^T$ becomes a product of Wishart matrices, whose {\it l.s.d.} is studied in
\cite{neuschel2014plancherel}. The variance of the eigenvalues of $\mathbf{J}\mathbf{J}^T$ is thus given by $\sigma_{\mathbf{J}\mathbf{J}^T}^2=L$. The largest eigenvalue of $\mathbf{J}\mathbf{J}^T$ is $\lambda_{max}=s_{max}^2=L^{-L}(L+1)^{L+1}$. For large $L$, it is observed that $\lambda_{max}$ scales as $\lambda_{max}\sim eL$. The linear growths of $\lambda_{max}$ and $\sigma_{\mathbf{J}\mathbf{J}^T}^2$ are validated in Fig. \ref{fig:valsta_linear_activation}. This means the breakdown of dynamical isometry and the poor conditioning in deep linear Gaussian networks.

\begin{figure*}[!t]
\begin{center}
\psfrag{f}[cc][cc][.6][0]{$L$}
\psfrag{b}[cc][cc][.6][0]{$\lambda_{max}$}
\psfrag{c}[cc][cc][.5][0]{$L$}
\psfrag{d}[cc][cc][.6][0]{$\sigma_{\mathbf{J}\mathbf{J}^T}^2$}

\subfigure[$\lambda_{max}$, \rm ReLU]{%
\epsfxsize=0.31\textwidth \leavevmode
\epsffile{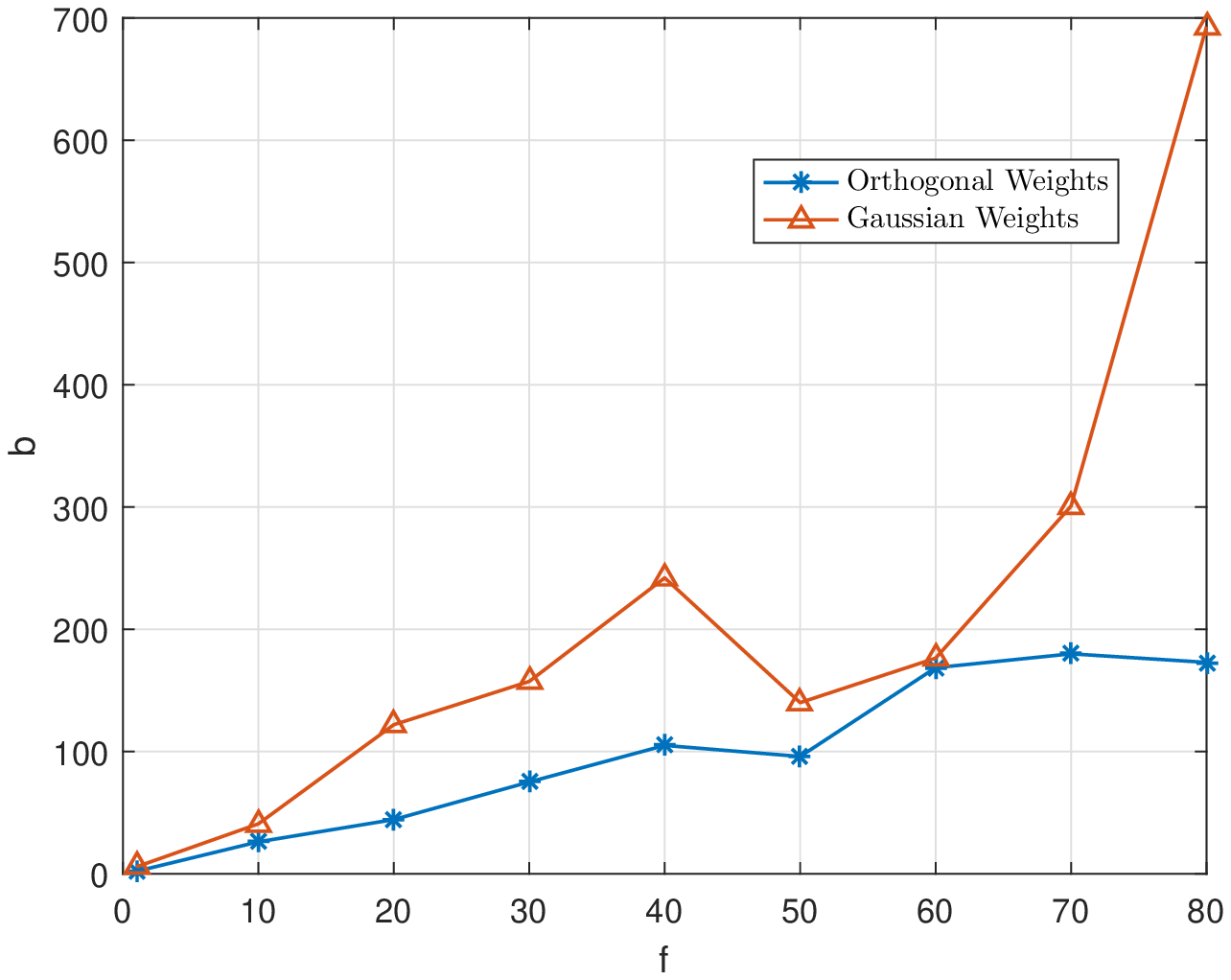}}\quad
\subfigure[$\lambda_{max}$, \rm hard-tanh]{%
\epsfxsize=0.31\textwidth \leavevmode
\epsffile{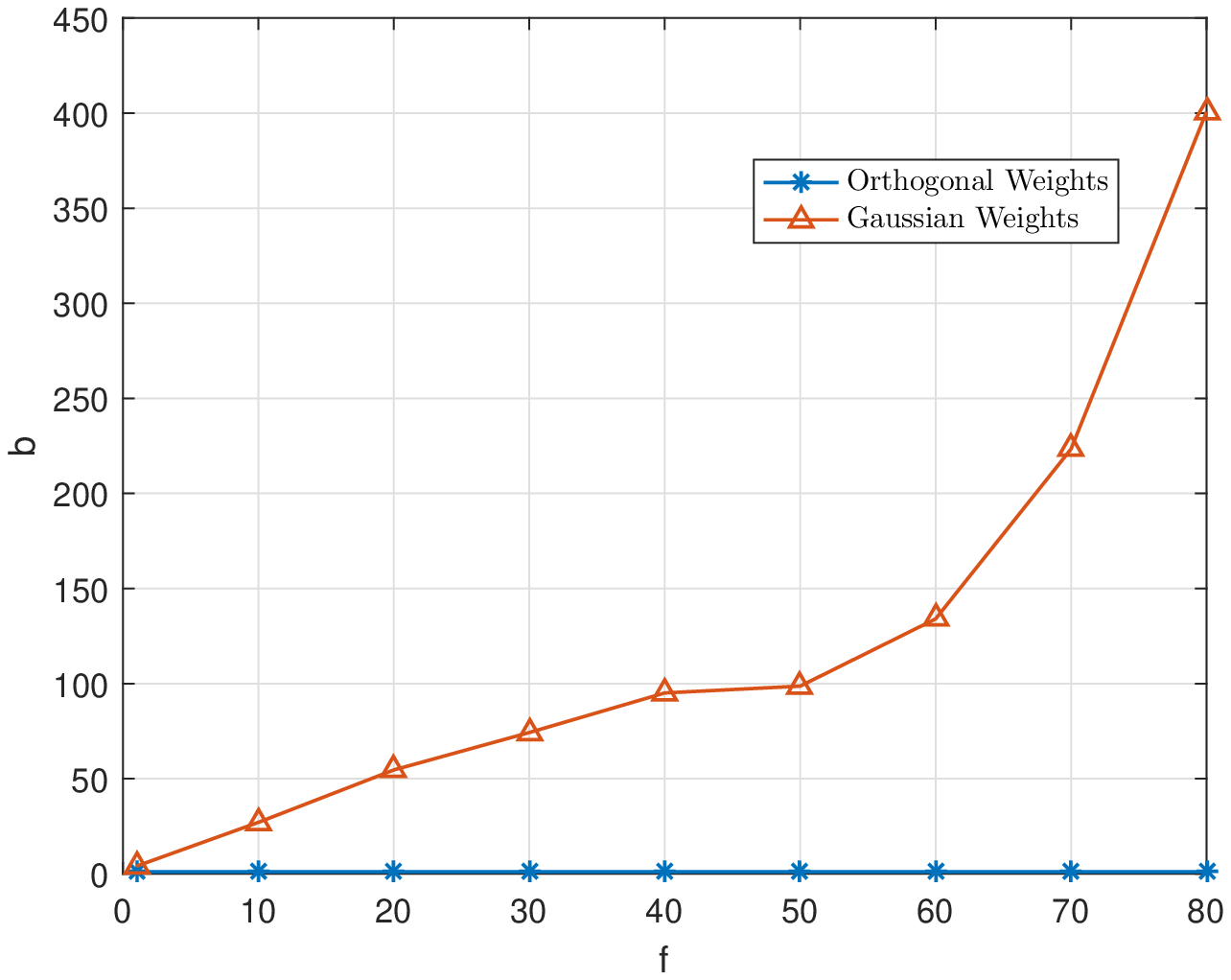}}\quad
\subfigure[$\lambda_{max}$, \rm tanh]{%
\epsfxsize=0.31\textwidth \leavevmode
\epsffile{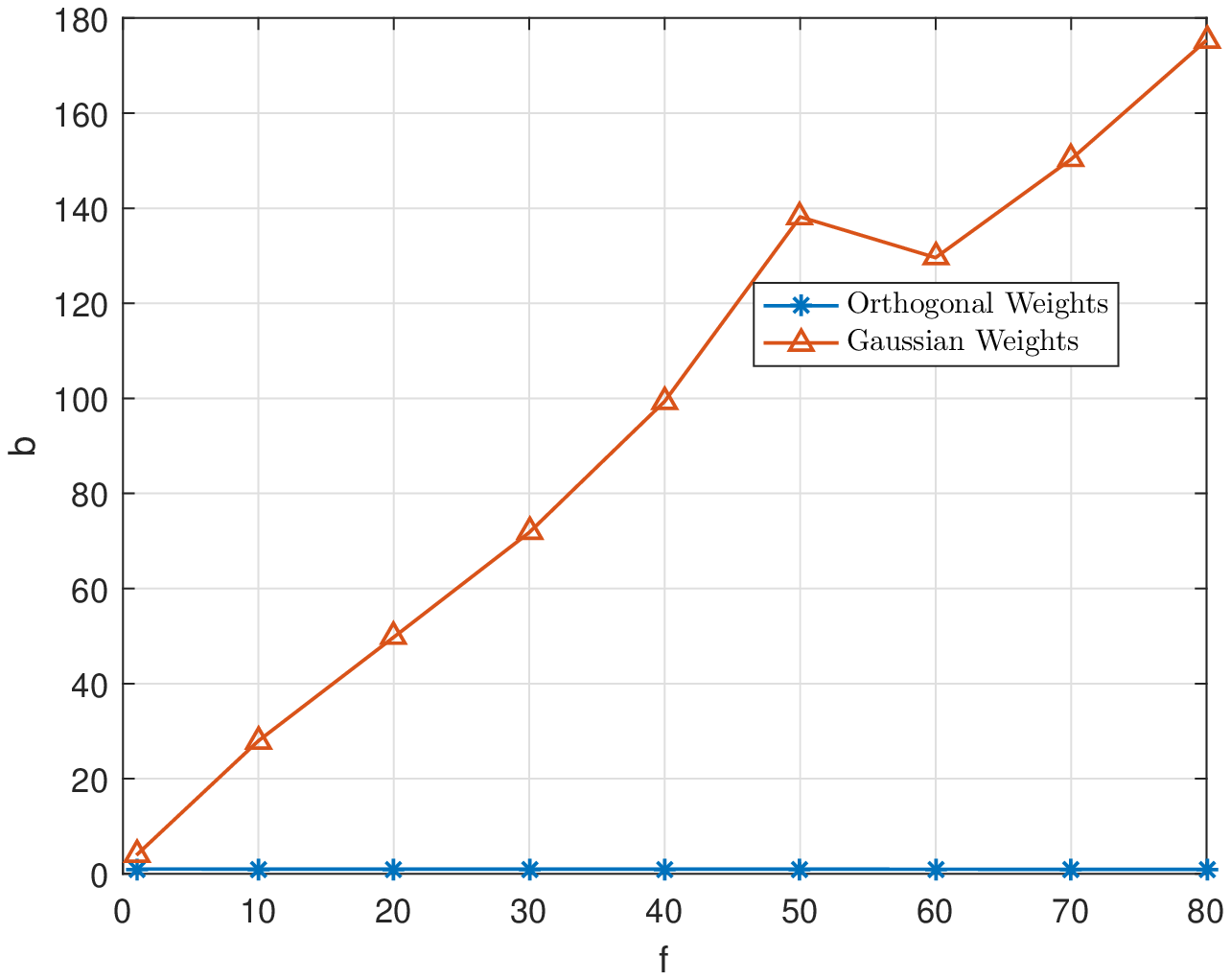}}\\

\subfigure[$\sigma_{\mathbf{J}\mathbf{J}^T}^2$, \rm ReLU]{%
\epsfxsize=0.31\textwidth \leavevmode
\epsffile{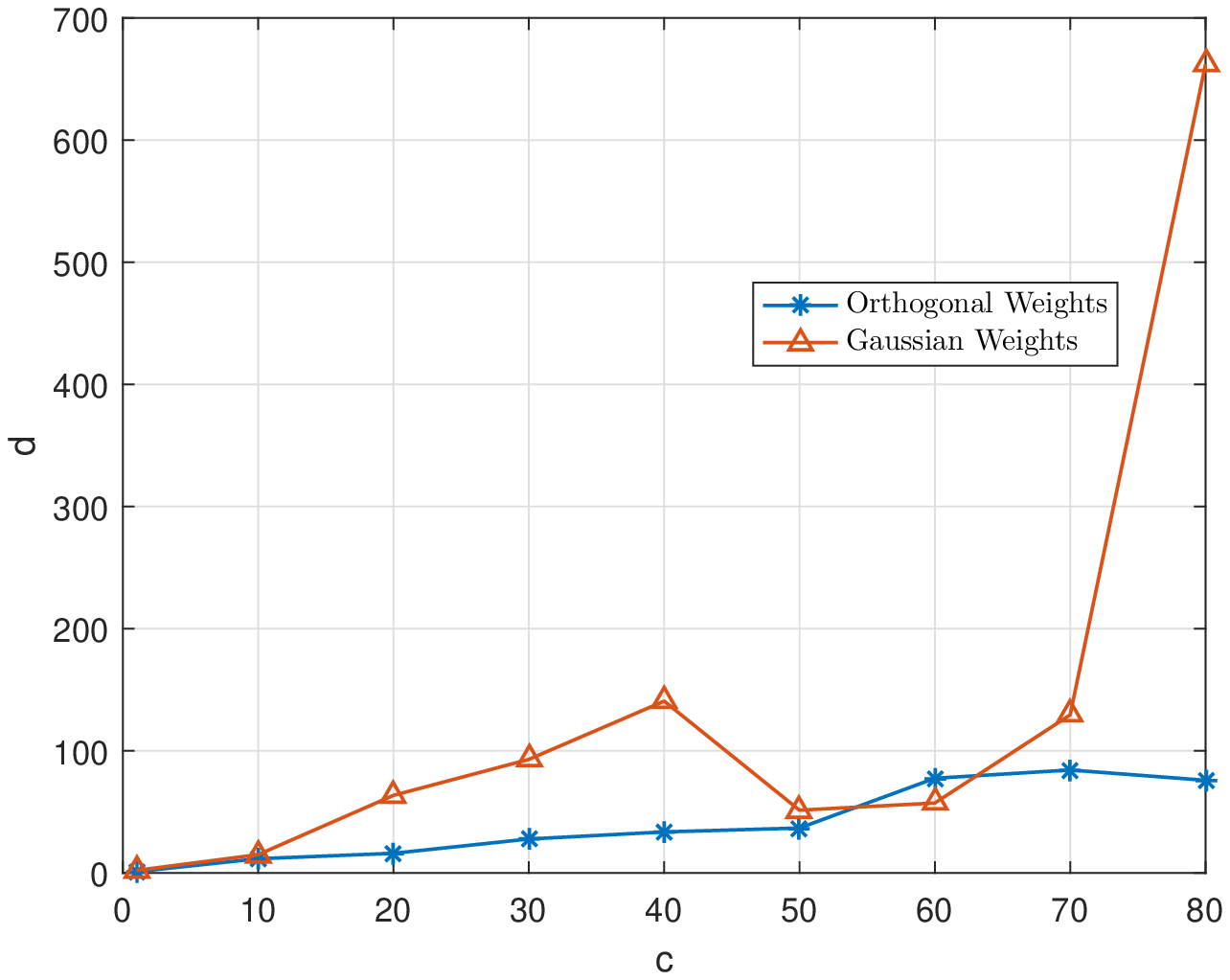}}\quad
\subfigure[$\sigma_{\mathbf{J}\mathbf{J}^T}^2$, \rm hard-tanh]{%
\epsfxsize=0.31\textwidth \leavevmode
\epsffile{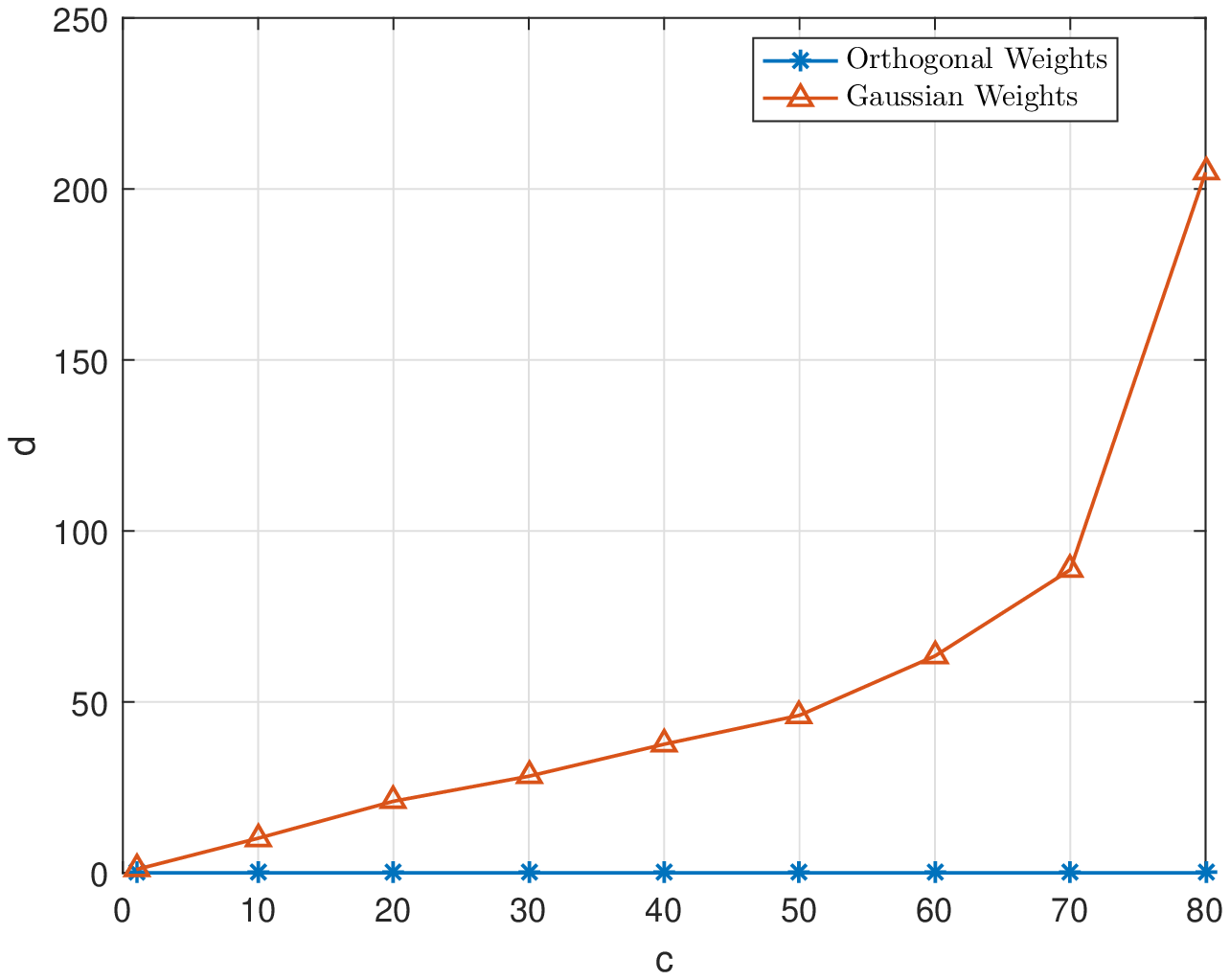}}\quad
\subfigure[$\sigma_{\mathbf{J}\mathbf{J}^T}^2$, \rm tanh]{%
\epsfxsize=0.31\textwidth \leavevmode
\epsffile{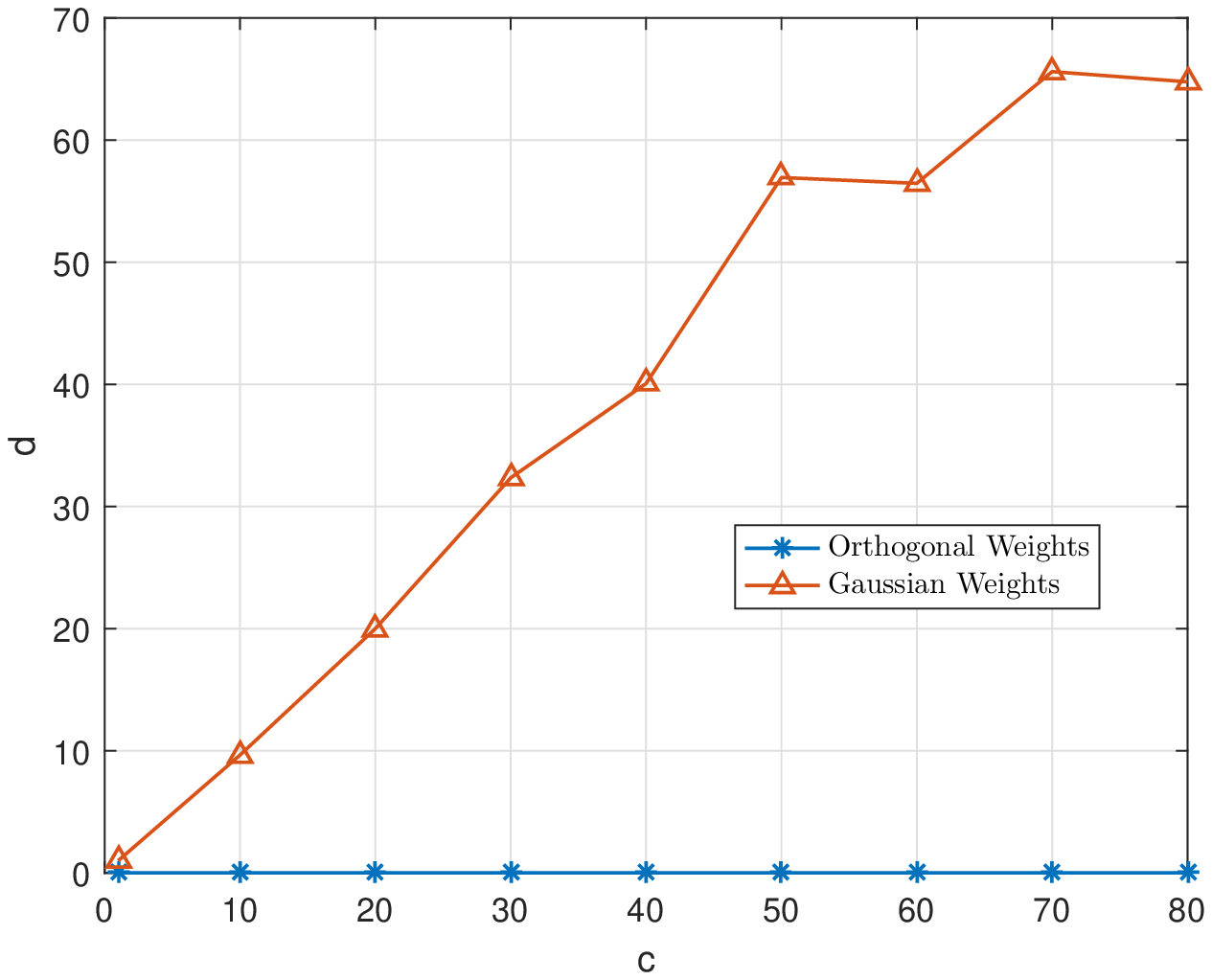}}

\caption{Variations of $\lambda_{max}$ and $\sigma_{\mathbf{J}\mathbf{J}^T}^2$ with respect to $L$ for different combinations of nonlinear activation functions and weight initializations. The width of the neural network is set to $1000$. \re{The results are obtained on a single realization.}}\label{fig:valsta_rth}
\end{center}
\end{figure*}

For the nonlinear networks, the random Gaussian weights and the random orthogonal weights are also respectively studied. When the Gaussian weights are adopted, we have
\begin{equation}\label{eq:largestevgaussian}
\lambda_{max}=s_{max}^2 = (\sigma_w^2p(q^*))^L\left(\frac{e}{p(q^*)}L+\mathcal{O}(1)\right),
\end{equation}
\begin{equation}\label{eq:variancegaussian}
\sigma_{\mathbf{J}\mathbf{J}^T}^2 = \frac{L}{p(q^*)},
\end{equation}
where $p(q^*)$ is the probability that a given neuron works in the linear regime with $\phi'(h)=1$, and it can be also explained as the fraction of neurons operating in the linear regime. For both the rectified-linear-unit (ReLU) and hard-tanh neural networks \footnote{\re{In this paper, the neural networks that only employ ReLU activation functions are referred to as ReLU neural networks. The other neural networks are defined in the same way.}}, we obviously always have $p(q^*)<1$, and this means that the Gaussian initializations can not realize dynamical isometry in the deep neural networks. Under the case where the random orthogonal weights are adopted, we have
\begin{equation}\label{eq:largestevorthogonal}
\lambda_{max}=s_{max}^2 = (\sigma_w^2p(q^*))^L\frac{1-p(q^*)}{p(q^*)}\frac{L^L}{(L-1)^{L-1}},
\end{equation}
\begin{equation}\label{eq:varianceorthogonal}
\sigma_{\mathbf{J}\mathbf{J}^T}^2 = \frac{1-p(q^*)}{p(q^*)}L.
\end{equation}
For ReLU networks, $p(q^*)=1/2$, and it can be seen that $\lambda_{max}$ and $\sigma_{\mathbf{J}\mathbf{J}^T}^2$ grow linearly with the depth. As a consequence, the dynamical isometry can not be realized in ReLU networks. However, in hard-tanh networks, $p(q^*)=erf(\frac{1}{\sqrt{2q^*}})$, thus we can tune $q^*$ to make $p(q^*)\approx 1-\frac{1}{L}$. In this way, the dynamical isometry is achievable in the orthogonal hard-tanh networks. In Fig. \ref{fig:valsta_rth}, with properly selected $q^*$ keeping $\chi$ around $1$,  the variations of $\lambda_{max}$ and $\sigma_{\mathbf{J}\mathbf{J}^T}^2$ with respect to $L$ are investigated for different combinations of nonlinear activation functions and weight initializations. It is shown that $\lambda_{max}$ and $\sigma_{\mathbf{J}\mathbf{J}^T}^2$ grow large as the depth of the neural network increases for all cases with random Gaussian weight initializations. However, for the cases with random orthogonal weight initializations, the hard-tanh and tanh neural networks with small $q^*$ can perform perfect dynamical isometry, i.e., all the eigenvalues concentrate at $1$ since  $\lambda_{max}$ is near to $1$ and $\sigma_{\mathbf{J}\mathbf{J}^T}^2$ is near to zero.

In \cite{pennington2018emergence}, the above methods to obtain the entire distribution of the singular values of the Jacobian are further developed to a calculational framework, which is useful in studying what combinations of nonlinear activation functions and weight initializations can yield the well conditioning that speed up the learning process. With the calculational framework, various combinations of weights initializations and nonlinear activation functions are analyzed. The results show that, beyond the hard-tanh activation function, a wide variety of nonlinear activation functions can realize dynamical isometry with random orthogonal weight initialization as the depth goes to infinity.

\re{However, the above results are more or less built on the free probability theory. In other words, the results only are true when the asymptotic freeness between every two matrix components in \eqref{eq:Jacobian} holds. Therefore, the applicability of the results in \cite{pennington2017resurrecting, pennington2018emergence} needs to be justified in practice \cite{pastur2020random, pastur2020random1}. As such, \cite{pastur2020random} provides a more complete proof for the results in \cite{pennington2017resurrecting, pennington2018emergence} under the Gaussian case, where the input data, the random weights and biases are assumed to be {\it i.i.d.} Gaussian variables. The proof in \cite{pastur2020random} is completed via rather standard techniques (e.g., Poincar\'{e}-Nash inequality \cite{pastur2011eigenvalue}) from RMT instead of directly applying conclusions from the free probability theory. Furthermore, in spirit of universality, \cite{pastur2020random1} extends the results in \cite{pastur2020random} to a more general framework analyzing the spectrum of input-output Jacobian under a more general {\it i.i.d.} case, where the random weights and biases are just {\it i.i.d.} variables with zero mean and finite fourth-order moment, but non-necessarily Gaussian. Thus, the line of works \cite{pastur2020random,pastur2020random1} actually give us a more general, more standard, and therefore more reliable analytical framework for the spectrum of the input-output Jacobian.
}

\re{
Beyond the conventional feed-forward neural networks, how to enable dynamical isometry in RNNs and CNNs is also studied in \cite{chen2018dynamical} and \cite{xiao2018dynamical}, respectively. Different from the conventional feed-forward neural networks, a mean field theory is introduced to analyze the signal propagations in RNNs. In particular, \cite{chen2018dynamical} develops the duality between the forward-propagation process of the signal and the back-propagation process of gradients in RNN. Overall, the input-output Jacobian spectra of RNN can be analyzed via RMT and the additional mean field theory, therefore the methods to achieve dynamical isometry can be developed. The simulation results in \cite{chen2018dynamical} show that a variety of RNNs with proper initializations achieving dynamical isometry are significantly easier to train. Analogously, mean field theory can be also utilized to analyze the signal propagation in CNNs \cite{xiao2018dynamical}. Furthermore, \cite{xiao2018dynamical} identifies an efficient construction approach for the convolution operators to facilitate random orthogonal initialization, therefore enables dynamical isometry in CNNs. As shown in the experimental results, the proposed construction method can speed up the training process of CNNs.
}

\subsection{Looking into the Loss Surface via the Hessian of the Weight Matrix}\label{subsec:Hessian}

In deep learning, training the neural network is actually optimizing a non-convex loss function, i.e., finding the global minimum of the loss surface, which is a geometric representation of the loss function \cite{granziol2019towards}. It is shown that even training a very simple neural network yields an intractable NP-complete problem \cite{blum1992training}. Thus, in the early stage, the neural networks were not favored compared to the classical machine learning methods that require only convex optimization. However, we all can see that nowadays the neural networks have achieved great practical successes in various fields. Despite some empirical or theoretical results which suggest that the local minimum is rarely an issue in large networks \cite{dauphin2014identifying, choromanska2015loss}, it is still hard to totally understand how the stochastic-gradient-descent (SGD) optimizer and simulated annealing methods make non-convex optimization problem tractable in the deep networks. Since the dimensions of the neural network and the input data are extremely large, RMT is considered as a powerful tool to explain the inner mechanism of deep learning. In this section, we will show the recent efforts made in understanding the loss surface of neural networks via RMT.

There are a few prior works that focus on the loss surface of the neural networks. Both \cite{choromanska2015loss} and \cite{dauphin2014identifying} show the prevalence of the saddle points as dominant critical points that plague the training process. In \cite{choromanska2015loss}, the authors propose to approximate the loss function with the Hamiltonian of the spherical spin-glass model, which originates from condensed matter physics. Therefore, the existence of the local minima at low loss values and saddle points at high loss values can be predicted via the knowledge of spherical spin-glass model from statistical physics. In addition, the existences of numerous local minima at low loss values are also highlighted. The related ideas are further investigated in \cite{kawaguchi2016deep, freeman2016topology, safran2016quality}. In \cite{dauphin2014identifying}, it is found that the {\it l.s.d.} of the Hessian at a critical point is a function of the loss value. Moreover, the shape of the spectrum of the Hessian at a critical point is similar to that of the semicircular law \cite{bray2007statistics}. In particular, the spectrum of the Hessian at the local minima is shifted right so much that all the eigenvalues of the Hessian are positive. On the contrary, the eigenvalues of the Hessian at the saddle points distribute around $0$, this means more negative eigenvalues exist in the spectrum of the Hessian. Therefore, the saddle points can be distinguished out via the faction of the negative eigenvalues of the Hessian. Besides, the Hessian contains more information about the loss surface. For example, the condition number of the Hessian determines the convergence rates of the first-order optimization methods on convex objectives \cite{nesterov2013introductory}. The existence of the negative eigenvalues of the Hessian indicates the non-convexity even at a local scale. Hessian analysis has been becoming a promising approach to study the geometric properties of the loss surface. In the following, we will introduce an RMT-based analytical framework for studying the spectra of the Hessian of the neural networks, which is proposed in \cite{pennington2017geometry}.

Considering a single-hidden-layer neural network without bias for simplicity, we denote the weight matrices by $\mathbf{W}^1\in\mathbb{R}^{n_1\times n_0}$ and $\mathbf{W}^2\in\mathbb{R}^{n_2\times n_1}$. Besides, the input data and output targets are denoted by $\mathbf{X}\in\mathbb{R}^{n_0\times m}$ and $\mathbf{Y}\in\mathbb{R}^{n_2\times m}$, where $n_0$, $n_1$, $n_2$, $m$ denote the input dimension, the number of neurons in the single layer, the output dimension, the number of data samples, respectively. In addition, the ReLU nonlinear activation function is employed, i.e., $\phi(z)=[z]_{+}=\max(z, 0)$. Therefore, the network output is given by
\begin{equation}\label{eq:network_model_Hessian}
\hat{\mathbf{Y}} = \mathbf{W}^2\phi(\mathbf{W}^1\mathbf{X}).
\end{equation}
The errors between the network output and the targets (a.k.a. the labels) are $e_{i\mu}=\hat{Y}_{i\mu}-Y_{i\mu}$, where $\mu$ is to index the samples. Considering the mean squared error, the loss value is given by
\begin{equation}\label{eq:loss_func}
\mathcal{L}=n_2\epsilon=\frac{1}{2m}\sum\limits_{i,\mu=1}^{n_2,m}e_{i\mu}^2,
\end{equation}
where $\epsilon$ is defined as the energy in the context and it actually characterizes the variance of the errors. The Hessian, denoted by $\mathbf{H}$, is defined as the matrix of second derivatives of the loss function with respect to the weights, namely, $H_{\alpha\beta}=\frac{\partial^2\mathcal{L}}{\partial\theta_\alpha\partial\theta_\beta}$, where $\theta_\alpha, \theta_\beta\in\{\mathbf{W}^1,\mathbf{W}^2\}$. $\mathbf{H}$ can be decomposed into two parts, $\mathbf{H}=\mathbf{H}_0+\mathbf{H}_1$, where $\mathbf{H}_0$ is a positive semi-definite matrix; $\mathbf{H}_1$ comes from the second derivatives and is therefore a symmetric matrix. More specifically, $\mathbf{H}_0$ and $\mathbf{H}_1$ are respectively given by
\begin{equation}\label{eq:Hessian_H0}
[H_0]_{\alpha\beta}\equiv\frac{1}{m}\sum\limits_{i,\mu=1}^{n_2,m}\frac{\partial \hat{Y}_{i\mu}}{\partial\theta_\alpha}\frac{\partial \hat{Y}_{i\mu}}{\partial\theta_\beta}
\equiv\frac{1}{m}[\mathbf{J}\mathbf{J}^T]_{\alpha\beta}
\end{equation}
and
\begin{equation}\label{eq:Hessian_H1}
[H_1]_{\alpha\beta}\equiv\frac{1}{m}\sum\limits_{i,\mu=1}^{n_2,m}e_{i\mu}\left(\frac{\partial^2 \hat{Y}_{i\mu}}{\partial\theta_\alpha\partial\theta_\beta}\right).
\end{equation}
It is worth noting that $\mathbf{J}$ in \eqref{eq:Hessian_H0} is the weight-output Jacobian, which is totally different from the input-output Jacobian in Section \ref{subsec:Jacobian}. The square neural networks where $n\equiv n_0=n_1=n_2$ are considered. In addition, we are interested in the asymptotic regime where both the network size and the data sets are very large. Besides, the limit ratio of the number of parameters to the effective number of samples, i.e., $c\triangleq2n^2/mn=2n/m$, is defined to characterize the network capacity. As we will see, $c$ is a important parameter that governs the shape of Hessian spectrum. From \eqref{eq:Hessian_H0}, it can be observed that $c$ also governs the rank of $\mathbf{H}_0$ since it determines the rank of $\mathbf{J}$.

To begin with, we make the following assumptions on the random neural network for the later derivation.
\begin{description}
  \item[$\mathcal{AS} 1$:] $\mathbf{H}_0$ and $\mathbf{H}_1$ are {\it freely independent}.
  \item[$\mathcal{AS} 2$:] The errors are {\it i.i.d.} Gaussian random variable $e_{i\mu}\sim \mathcal{N}(0, 2\epsilon)$. This assumption makes the gradients vanish in the large $m$ regime, specifying the analysis to critical points.
  \item[$\mathcal{AS} 3$:] Both the input data and the weights are {\it i.i.d.} Gaussian random variables.
\end{description}
The assumptions are quite mild in the random neural networks, and the reasonability of them is particularly discussed in \cite{pennington2017geometry}.

\re{
Under $\mathcal{AS} 1$, the Hessian $\mathbf{H}$ becomes a summation of two {\it freely independent} matrices, i.e., $\mathbf{H}_0$ and $\mathbf{H}_1$, and the spectrum of $\mathbf{H}$ can therefore be derived using the free probability theory.
With R-transform and the free probability theory, we get a general framework to compute the spectrum of the Hessian in steps: i) compute the Stieltjes transform of the {\it l.s.d.} of $\mathbf{H}_0$ and $\mathbf{H}_1$; ii) derive the corresponding R-transforms, i.e., $R_{\mathbf{H}_0}$ and $R_{\mathbf{H}_1}$, according to \eqref{eq:R-transform}; iii) obtain $R_{\mathbf{H}}$ via \eqref{eq:R-transform_sum} and further the Stieltjes transform of the {\it l.s.d.} of $\mathbf{H}$; iv) calculate the {\it l.s.d.} of $\mathbf{H}$ using the inverse Stieltjes transform.
}

Similar to quantum physics, we first simplify the Hessian by approximating $\mathbf{H}_0$ and $\mathbf{H}_1$ with random matrices. With the structural features of $\mathbf{H}_0 =\frac{1}{m}\mathbf{J}\mathbf{J}^T$ and $\mathbf{H}_1$, $\mathbf{H}_0$ and $\mathbf{H}_1$ are approximated with Wishart matrices and Wigner matrices, respectively. Therefore, the Hessian can be approximated with the Wishart-plus-Wigner model. Specifically, we assume that the elements of both $\mathbf{J}$ and $\mathbf{H}_1$ are {\it i.i.d.} Gaussian random variables. Hence, the spectra of $\mathbf{H}_0$ and $\mathbf{H}_1$ can be described with the general forms of the Mar\v{c}enko-Pastur distribution and the semi-circular distribution, respectively. Taking $\sigma^{\mathbf{H}_0}=1$ and $\sigma^{\mathbf{H}_1}=\sqrt{2\epsilon}$, the {\it l.s.d.} of $\mathbf{H}_0$ and $\mathbf{H}_1$ can be obtained as follows via \eqref{eq:SC_law_scaled} and \eqref{eq:MP_law_scaled}:
\begin{equation}\label{eq:lsd_H0}
f^{\mathbf{H}_0} = f_{MP}(\lambda; c, 1),
\end{equation}
and
\begin{equation}\label{eq:lsd_H1}
f^{\mathbf{H}_1} = f_{SC}(\lambda; \sqrt{2\epsilon}).
\end{equation}
According to \eqref{eq:intuition_Stieltjes} and \eqref{eq:R-transform}, we have
\begin{equation}\label{eq:R-transform_H0}
R_{\mathbf{H}_0} = \frac{1}{1-zc},
\end{equation}
and
\begin{equation}\label{eq:R-transform_H1}
R_{\mathbf{H}_1} = 2\epsilon z.
\end{equation}
Obviously, the R-transform of $f^{\mathbf{H}}$ can be derived as
\begin{equation}\label{eq:R-transform_H}
R_{\mathbf{H}} =\frac{1}{1-zc} + 2\epsilon z.
\end{equation}
The Stieltjes transform of $f^{\mathbf{H}}$ can be obtained through solving the following cubic equation,
\begin{equation}\label{eq:Stieltjes_H}
2\epsilon cm_{F^{\mathbf{H}}}^3-(2\epsilon+zc)m_{F^{\mathbf{H}}}^2+(z+c-1)m_{F^{\mathbf{H}}}-1=0.
\end{equation}

\begin{figure*}[!t]
\begin{center}
\psfrag{x}[cc][cc][.5][0]{$\lambda(\mathbf{H})$}
\psfrag{y}[cc][cc][.5][0]{$f^{\mathbf{H}}(\lambda)$}

\subfigure[$c=\frac{1}{3}$]{%
\epsfxsize=0.45\textwidth \leavevmode
\epsffile{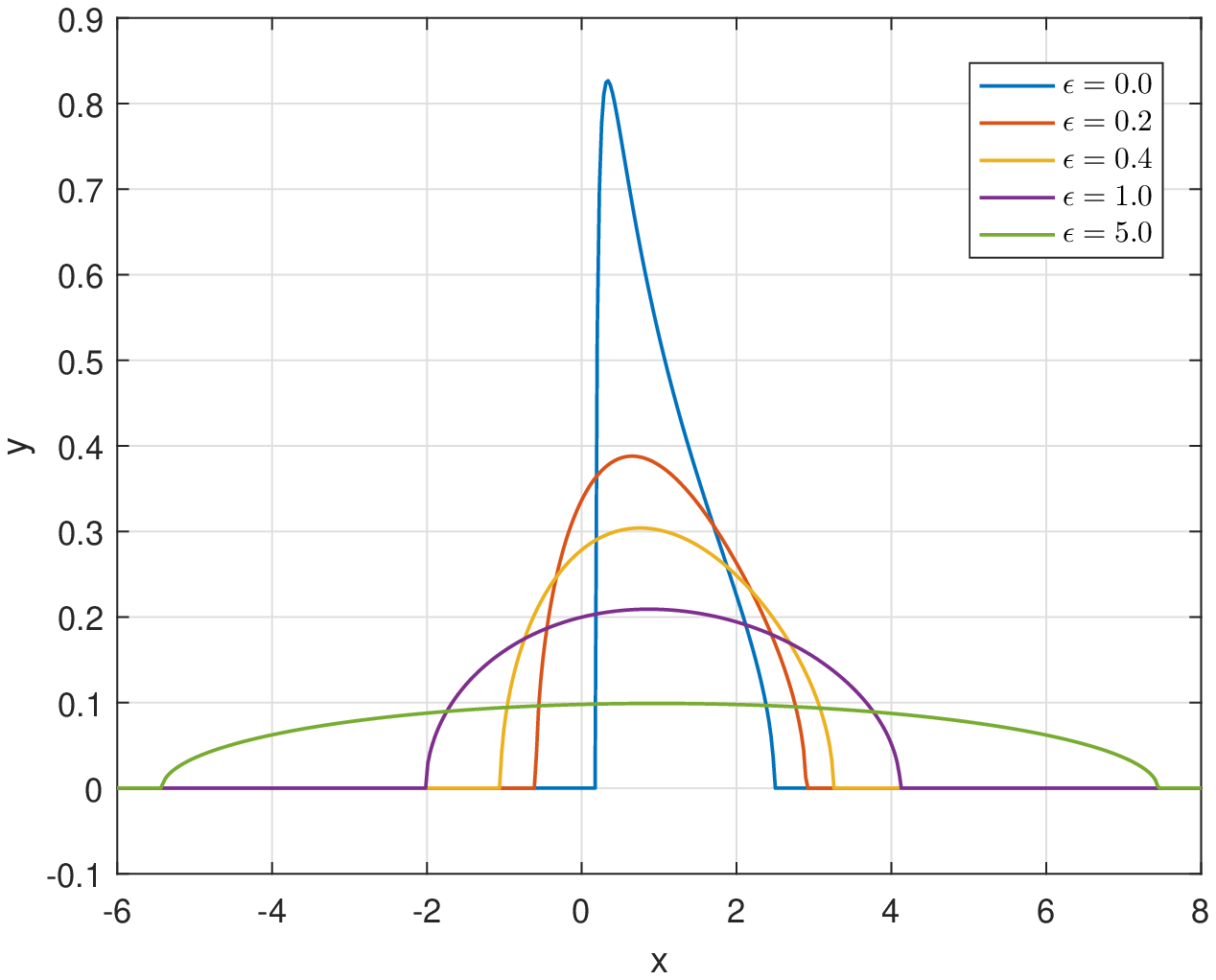}}\quad
\subfigure[$c=\frac{2}{3}$]{%
\epsfxsize=0.45\textwidth \leavevmode
\epsffile{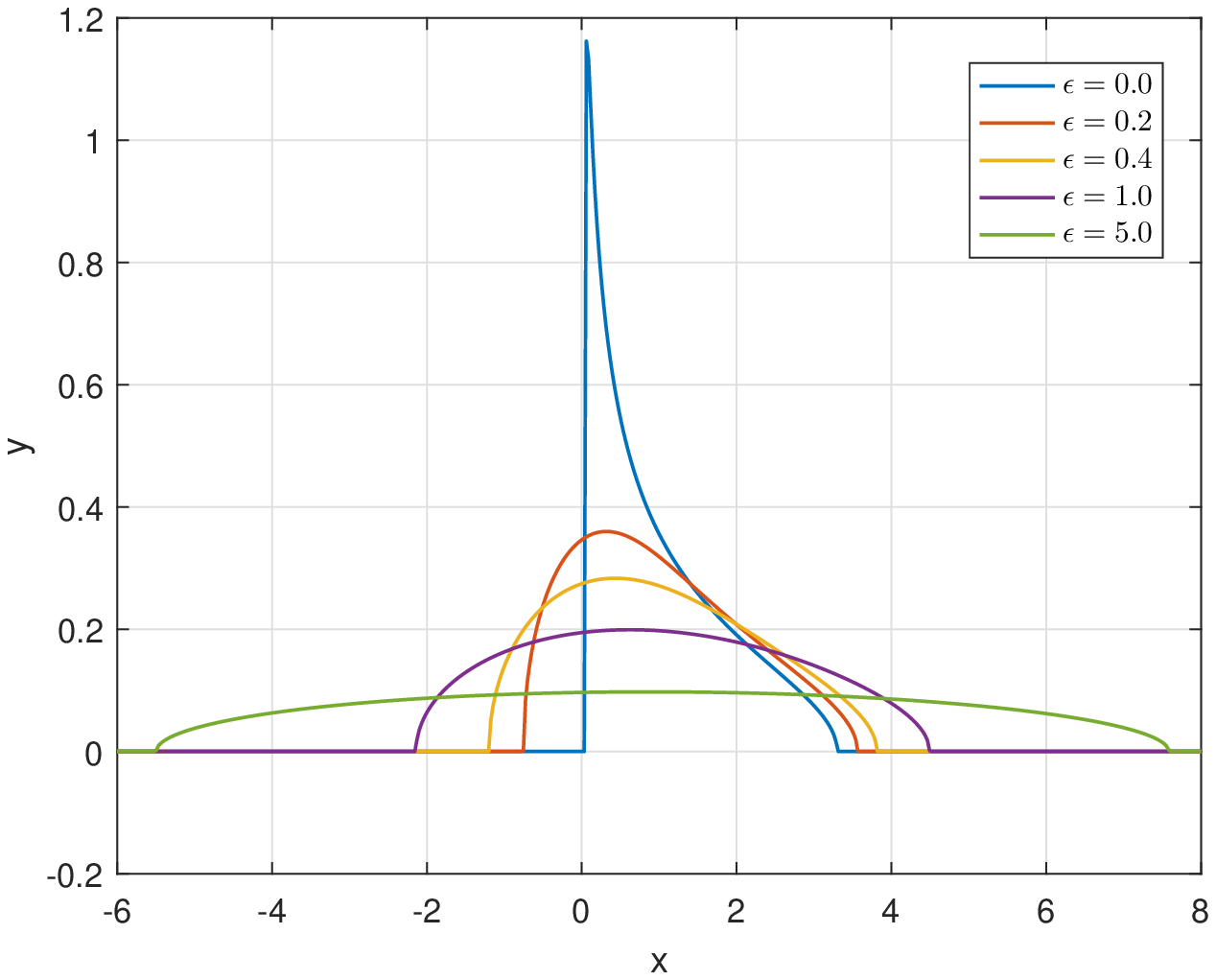}}\quad

\caption{Theoretical limit spectrum density of the Hessian at the critical points with different $\epsilon$'s and $c$'s.}\label{fig:WishartPlusWigner}
\end{center}
\end{figure*}

Finally, we can obtain $f^{\mathbf{H}}$ via the inverse Stieltjes transform and the Hessian spectra with different $c$'s and $\epsilon$'s are shown in Fig. \ref{fig:WishartPlusWigner}. Intriguingly, it can be observed that the shape of spectrum density of the Hessian at the critical point approaches the Mar\v{c}enko-Pastur distribution when $\epsilon$ is small enough. However, as $\epsilon$ grows large, $f^{\mathbf{H}}$ behaves more and more similar to the semi-circular distribution. Noting that $\epsilon$ is proportional to the loss value, we can therefore distinguish the saddle points at high loss values by observing the spectrum of the Hessian. Based on this, a more advanced quantity, namely, the {\it normalized index}, is induced to identify the critical points.

Obviously, $f^{\mathbf{H}}$ is a function parameterized by $\epsilon$ and $c$. The normalized index, or the fraction of the negative eigenvalues of the Hessian, is defined as \cite{bray2007statistics}
\begin{equation}\label{eq:normalizedindex}
\alpha(\epsilon, c)\triangleq \int_{-\infty}^{0}f^{\mathbf{H}}(\lambda;\epsilon,c)\dif\lambda = 1-\int_{0}^{\infty}f^{\mathbf{H}}(\lambda;\epsilon,c)\dif\lambda.
\end{equation}
It is observed that the normalized index of the critical points grows rapidly with $\epsilon$ in \cite{dauphin2014identifying, choromanska2015loss}, so that the critical points with many descent directions have large loss values. In addition, it is found that for small $\alpha$,
\begin{equation}\label{eq:ni_smallalpha}
\alpha(\epsilon, c)\approx \alpha_0(c)\left|\frac{\epsilon-\epsilon_c}{\epsilon_c}\right|^{\frac{3}{2}},
\end{equation}
where
\begin{equation}\label{eq:improvedcriticalvalue}
\epsilon_c=\frac{1}{16}\left(1-20c-8c^2+(1+8c)^{\frac{3}{2}}\right),
\end{equation}
is the critical value of $\epsilon$ below which all the critical points are minimizers. Therefore, we can determine whether a critical point is a saddle point analytically by comparing the energy at a critical point with $\epsilon_c$.

\begin{figure*}[!t]
\begin{center}
\psfrag{a}[cc][cc][.5][0]{$\lambda(\mathbf{H})$}
\psfrag{b}[cc][cc][.5][0]{$f^{\mathbf{H}}(\lambda)$}
\psfrag{c}[cc][cc][.5][0]{$\lambda(\mathbf{H}_0)$}
\psfrag{d}[cc][cc][.5][0]{$f^{\mathbf{H}_0}(\lambda)$}
\psfrag{e}[cc][cc][.5][0]{$\lambda(\mathbf{H}_1)$}
\psfrag{f}[cc][cc][.5][0]{$f^{\mathbf{H}_1}(\lambda)$}

\subfigure[$\epsilon=5$]{%
\epsfxsize=0.31\textwidth \leavevmode
\epsffile{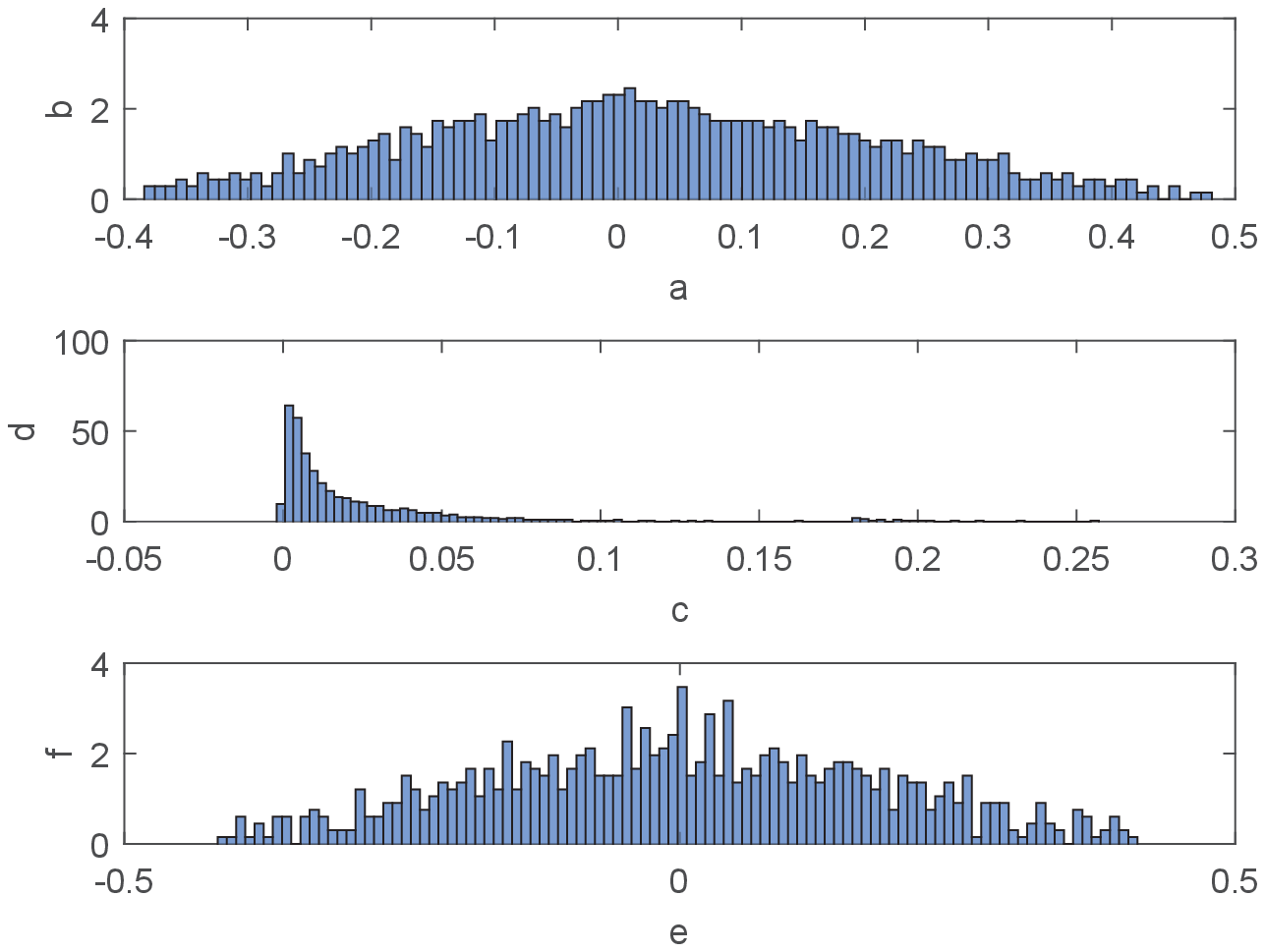}}\quad
\subfigure[$\epsilon=1$]{%
\epsfxsize=0.31\textwidth \leavevmode
\epsffile{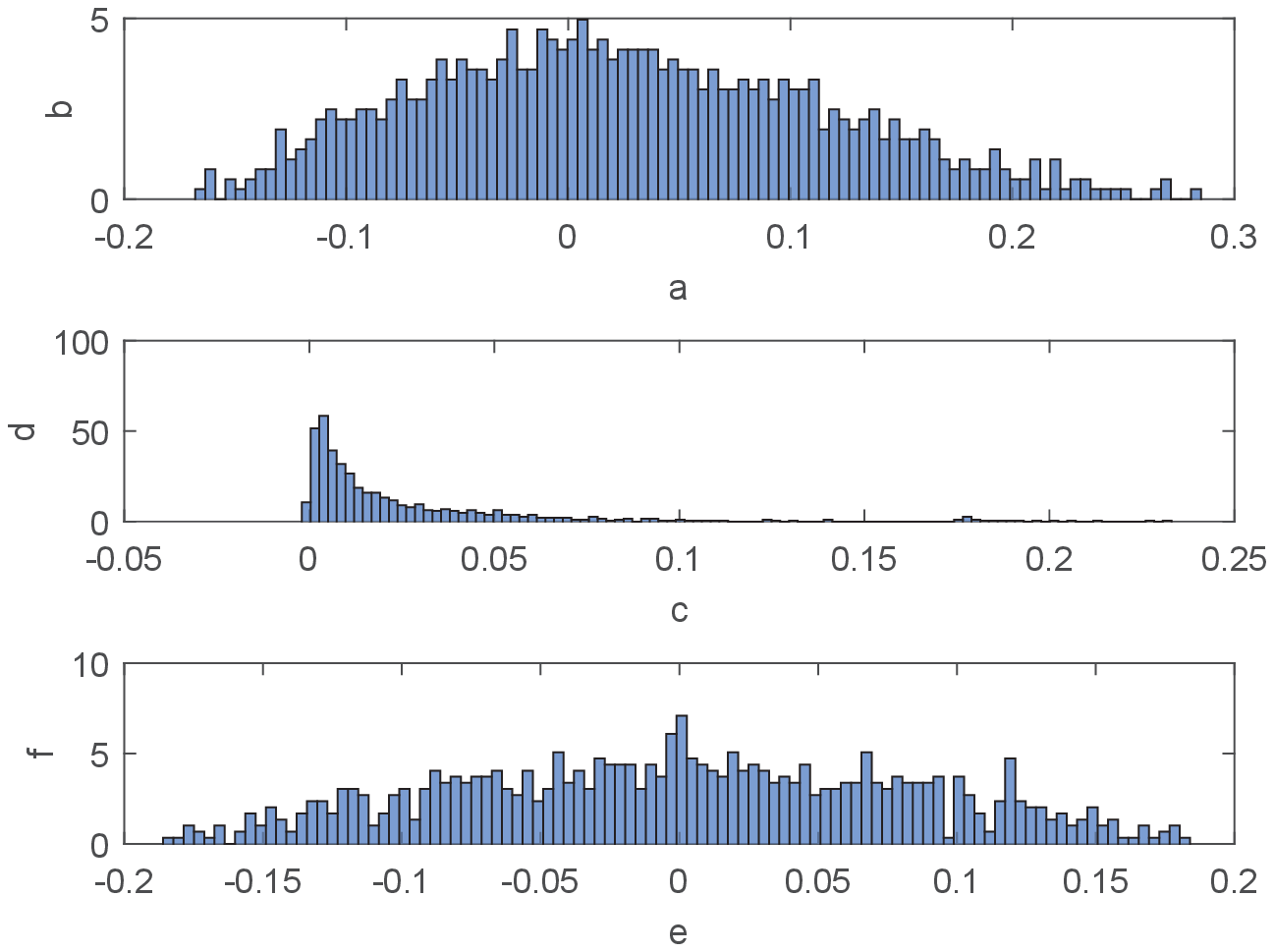}}\quad
\subfigure[$\epsilon=10^{-3}$]{%
\epsfxsize=0.31\textwidth \leavevmode
\epsffile{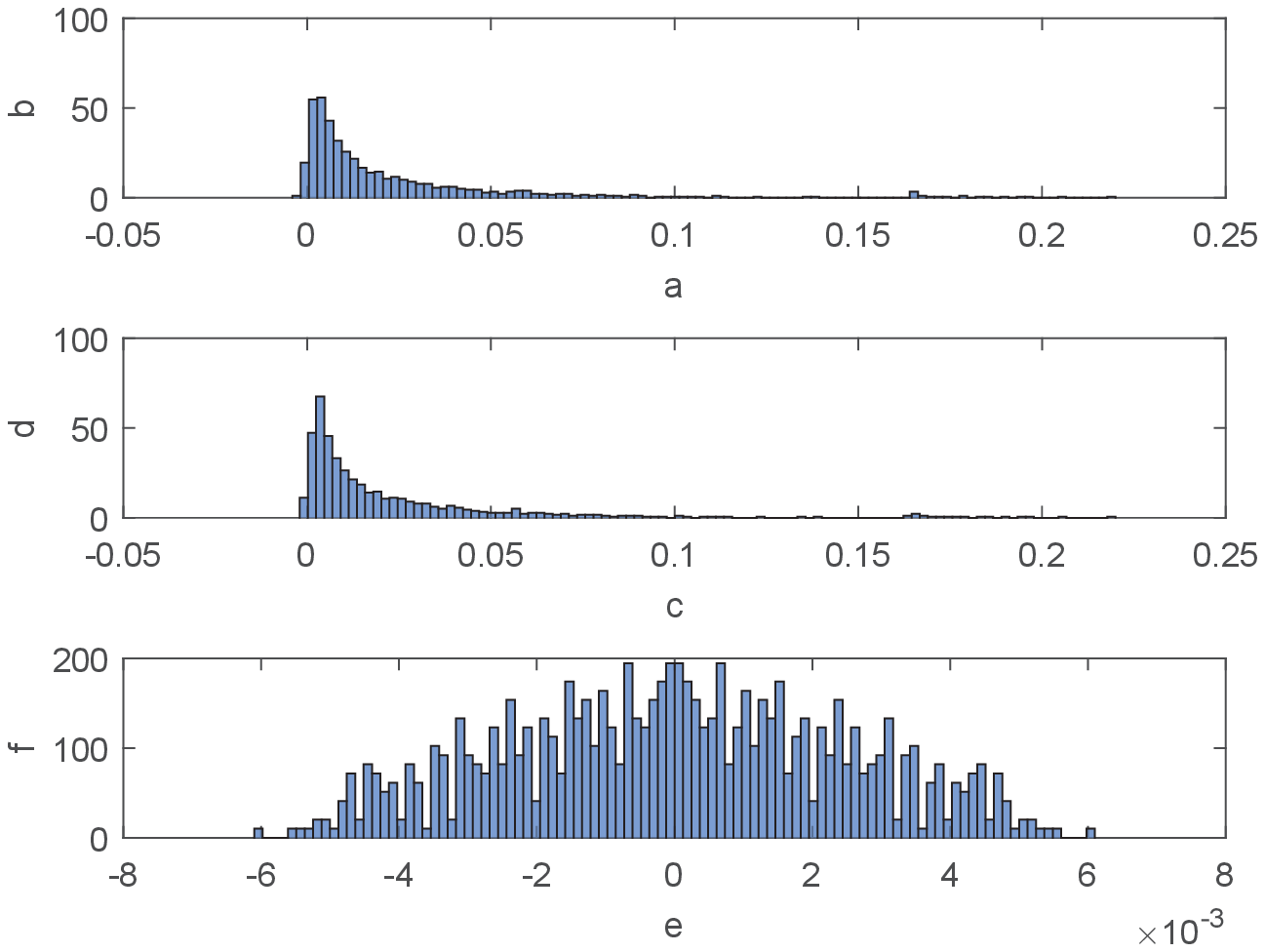}}

\caption{Empirical spectrum density of $\mathbf{H}$, $\mathbf{H}_0$, and $\mathbf{H}_1$ for the critical points with different levels of loss values in random neural networks.}\label{fig:Hessian_e_vals}
\end{center}
\end{figure*}

The above results mainly depend on $\mathcal{AS} 1$ -- $\mathcal{AS} 3$ and the additional assumption (denoted by $\mathcal{AS} 4$ for simplicity) that approximates $\mathbf{J}$ and $\mathbf{H}_1$ with {\it i.i.d.} Gaussian random variables. It is necessary to relax some unrealistic assumptions to acquire a deeper insight to the practical networks. In \cite{pennington2017geometry}, $\mathcal{AS} 1$ -- $\mathcal{AS} 3$ have been discussed in details and shown to be fairly mild. To validate $\mathcal{AS} 4$, we plot the empirical spectra of $\mathbf{H}$, $\mathbf{H}_0$ and $\mathbf{H}_1$ at the critical points with different levels of loss values in Fig. \ref{fig:Hessian_e_vals}.
\re{To be specific, the results in Fig. \ref{fig:Hessian_e_vals} are obtained in a single-layer random neural network as shown in \eqref{eq:network_model_Hessian} with $n_0=n_1=n_2=20$ and $m=160$. Besides, with the fact that $\epsilon$ is directly related to the loss values via \eqref{eq:loss_func}, we choose a set of parameters $\epsilon$'s with large gaps to show the difference of the Hessian spectra at critical points with different loss values more obviously.}
It is shown that the both the spectra of $\mathbf{H}_0$ and $\mathbf{H}_1$ deviate a little bit from the Mar\v{c}enko-Pastur distribution and the semicircular distribution.
Hence, more advanced and precise models are proposed to approximate the practical spectra of $\mathbf{H}_0$ and $\mathbf{H}_1$ and validated via numerical results \cite{pennington2017geometry}. The R-transforms of the spectra of $\mathbf{H}_1$ and $\mathbf{H}_0$ can be better approximated by
\begin{equation}\label{eq:improvedR-transform_H1}
R_{\mathbf{H}_1}(z) = \frac{\epsilon cz}{2-\epsilon c^2z^2},
\end{equation}
and
\begin{equation}\label{eq:improvedR-transform_H0}
R_{\mathbf{H}_0}(z) = \frac{\sigma}{1-\sigma zc},
\end{equation}
where $\sigma$ is an additional parameter to modify the Mar\v{c}enko-Pastur distribution so that it can better fit the spectrum of $\mathbf{H}_0$. Again, we can calculate the normalized index of critical points with energy $\epsilon$. Using the same techniques in obtaining \eqref{eq:ni_smallalpha}, we have
\begin{equation}\label{eq:improvednormalizedindex}
\alpha(\epsilon, c)\approx \tilde{\alpha}_0(c)\left|\frac{\epsilon-\epsilon_c}{\epsilon_c}\right|^{\frac{3}{2}},
\end{equation}
where $\tilde{\alpha}_0$ is used to show the difference with $\alpha_0$ in \eqref{eq:ni_smallalpha} and the critical value of $\epsilon$ is given by
\begin{equation}\label{eq:improvedcriticalvalue}
\epsilon_c=\frac{\sigma^2(27-18\chi-\chi^2+8\chi^{\frac{3}{2}})}{32c(1-c)^3}
\end{equation}
with $\chi=1+16c-8c^2$.

The most important step in the aforementioned computation framework is decomposing the Hessian as a summation of two freely independent matrices. This is further investigated for the practical deep neural networks in \cite{granziol2020beyond}. However, it is shown that the observed spectral shapes strongly deviate from the theoretical predictions even allowing for some outliers. With the numerical results obtained from the practical neural networks and data sets, they find that the spectra can be better approximated with the spectra of two new matrix ensembles, i.e., random Wigner/Wishart ensemble products and percolated Wigner/Wishart ensembles. One can see that, although RMT provides many useful tools to characterize the spectra of the Hessian of random neural networks, we still have a long way to go before totally understanding the loss surface of the practical deep networks.

\subsection{Designing the Nonlinearities to Preserve the Spectrum of the Data Covariance Matrix}\label{subsec:dataCM}

In deep learning, highly skewed spectra of data covariance matrices means strong anisotropy in the embedded feature space, which is regarded as an indicator of poor conditioning to impede the learning process \cite{pennington2017nonlinear}. The conventional solution is to introduce the batch normalization layer to rescale the variance of individual activations of the batch. However, the covariance is usually ignored. As a consequence, this may result in a large imbalance in singular values as the signal propagates through the neural networks. Hence, how to preserve the complete spectra of the data covariance matrices in the neural networks becomes an attractive question. Intriguingly, the following analysis of the data covariance matrix provides us another more efficient way to solve this problem from RMT.

The data covariance matrix, is actually the sample covariance matrix of the post-activations. For simplicity, we start from a single-layer neural network without bias. Here, we concatenate the random input vectors as a random data matrix $\mathbf{X}\in\mathbb{R}^{n_0\times m}$ with {\it i.i.d.} Gaussian elements $\mathbf{X}_{ij}\sim\mathcal{N}(0, \sigma_x^2)$, therefore the post-activation matrix of the neural network can be written as
\begin{equation}\label{eq:outputmatrix}
\mathbf{Y}=\phi(\mathbf{W}\mathbf{X}),
\end{equation}
where $\mathbf{W}\in\mathbb{R}^{n_1\times n_0}$ is the random weight matrix with {\it i.i.d.} Gaussian elements $\mathbf{W}_{ij}\sim\mathcal{N}(0, \sigma_w^2/n_0)$, and $\phi(\cdot)$ is the component-wise nonlinear activation function. In particular, $n_0$, $n_1$ denotes the input dimension and output dimension of the neural network, respectively; $m$ is the number of data samples in the data set. Besides, the asymptotic regime where $n_0$, $n_1$, and $m$ go to infinity with a constant rate is considered and we have some additional definitions as follows
\begin{equation}\label{eq:constantratio_nn}
\xi\triangleq \frac{n_0}{m}, \psi = \frac{n_0}{n_1},\ {\rm as}\ n_0, n_1, m \to\infty.
\end{equation}
In addition, a further assumption is needed for the nonlinear activation function. Denoting the pre-activation matrix as $\mathbf{Z}\triangleq\mathbf{W}\mathbf{X}$, let $\phi(\cdot)$ denote the activation function with zero mean and finite moments, i.e., $\phi(\cdot)$ satisfies
\begin{equation}\label{eq:assump_activation_func_zm}
\int\frac{\dif z}{\sqrt{2\pi}}e^{-\frac{z^2}{2}}\phi(\sigma_w\sigma_xz) = 0,
\end{equation}
and
\begin{equation}\label{eq:assump_activation_func_fm}
\left|\int\frac{\dif z}{\sqrt{2\pi}}e^{-\frac{z^2}{2}}\phi(\sigma_w\sigma_xz)^k\right|<\infty, \forall k>1.
\end{equation}
In the context, the Gram matrix $\mathbf{Y}\mathbf{Y}^T$ and output covariance matrix $\mathbf{F} =\frac{1}{m}\mathbf{Y}\mathbf{Y}^T$ are of our special interest. To be more specific, the literatures focus on the eigenvalues or the spectrum density of $\mathbf{F}$. We recall that the spectrum density function can be derived by calculating the corresponding Stieltjes transform. Noting the resolvent of $\mathbf{F}$ is defined as  $\mathbf{G}(z)=\left(\mathbf{F}-z\mathbf{I}_{n_1}\right)^{-1}$, according to \eqref{eq:intuition_Stieltjes}, the computation of the Stieltjes transform reduces to computing the trace of the resolvent, i.e.,  $$m_{\mathbf{F}}(z)=\frac{1}{n_1} tr\left(\mathbf{F}-z\mathbf{I}_{n_1}\right)^{-1}=\frac{1}{n_1}tr\mathbf{G}(z).$$

With the moment method in RMT \cite{tao2012topics}, $m_{\mathbf{F}}(z)$ can be computed and we can therefore obtain the spectrum of the output covariance matrix via the inverse Stieltjes transform. The results unfold as the following theorem \cite{pennington2017nonlinear}.
\begin{theorem}\label{th:mfz_nn}
Defining two constants $\eta$ and $\zeta$ as
\begin{equation}\label{eq:eta}
\eta = \int\frac{\dif ze^{-\frac{z^2}{2}}}{\sqrt{2\pi}}\phi(\sigma_w\sigma_xz)^2,
\end{equation}
\begin{equation}\label{eq:zeta}
\zeta = \left[\sigma_w\sigma_x\int\frac{\dif ze^{-\frac{z^2}{2}}}{\sqrt{2\pi}}\phi'(\sigma_w\sigma_xz)\right]^2,
\end{equation}
the Stieltjes transform of the spectrum density of $\mathbf{F}$ can be calculated by solving the following quart
\begin{equation}\label{eq:mfz_nn}
m_{\mathbf{F}}(z) = \frac{\psi}{z}P\left(\frac{1}{z\psi}\right) + \frac{1-\psi}{z},
\end{equation}
where
\begin{equation}\label{eq:P_nn}
P = 1 + (\eta-\zeta)tP_{\xi}P_{\psi} + \frac{P_{\xi}P_{\psi}t\zeta}{1-P_{\xi}P_{\psi}t\zeta},
\end{equation}
and
\begin{equation}\label{eq:P_xipsinn}
P_{\xi} = 1+(P-1)\xi, P_{\psi}=1+(P-1)\psi.
\end{equation}
\end{theorem}

In particular, we are interested in two special cases of \eqref{eq:mfz_nn}: $\eta=\zeta$ and $\zeta=0$. It is proved that $\eta=\zeta$ if and only if $\phi(\cdot)$ is a linear function, i.e., $\phi(z) = z$. In this case, $\mathbf{F}$ reduces to $\frac{1}{m}\mathbf{Z}\mathbf{Z}^T$, where $\mathbf{Z}=\mathbf{W}\mathbf{X}$ is a product of two Gaussian random matrices and the Stieltjes transform $m_{\mathbf{F}}(z)$ can be computed using the methods in \cite{dupic2014spectral}. Next, we will show that the other case, namely, $\zeta=0$, is more useful in designing the nonlinear activation functions. Without loss of generality, $\eta$ is set to $1$ while the general case can be recovered via a rescaling factor. When $\zeta=0$, \eqref{eq:mfz_nn} reduces to
\begin{equation}\label{eq:eta0}
z[m_{\mathbf{F}}(z)]^2 = \left((1-\frac{\psi}{\xi})z-1\right)m_{\mathbf{F}}(z)+\frac{\psi}{\xi} = 0,
\end{equation}
which is the exactly the Stieltjes transform of Mar\v{c}enko-Pastur distribution with parameter $c=\frac{\psi}{\xi}$. Noting that the input elements are assumed to be {\it i.i.d.} Gaussian random variables, the spectrum of the input covariance matrix also satisfies Mar\v{c}enko-Pastur law while the shape is governed by $\xi$. When $\psi=1$, we can observe that $\frac{1}{m}\mathbf{Y}\mathbf{Y}^T$ and $\frac{1}{m}\mathbf{X}\mathbf{X}^T$ have the same limit spectrum distribution, i.e., Mar\v{c}enko-Pastur distribution parameterized by $\xi$. So far, we have identified a novel type of nonlinear activation functions that can preserve the full spectra of the data covariance matrices as the signal propagates through the neural networks. Now we look back to the multi-layer neural networks, where the post-activation matrix of $l$-th layer is given by
\begin{equation}\label{eq:postam_llayer}
\mathbf{Y}^{l} = \phi(\mathbf{W}^{l}\mathbf{Y}^{l-1}), \mathbf{Y}^{0}=\mathbf{X}.
\end{equation}
Using the results in \eqref{eq:eta0}, we can design an activation function that satisfies $\zeta=0$ to approximately preserve the full singular value spectrum as the signal propagates through the neural networks, at least in the early training phase. With this observation, a lot of nonlinear activation functions can be designed to satisfy the condition $\zeta\approx 0$. This suggests that the design of the non-linear activation functions deserves further investigations to improve the learning speed of the training stage.

\begin{figure}[!t]
\begin{center}
\psfrag{x}[cc][cc][.7][0]{$x$}
\psfrag{y}[cc][cc][.7][0]{\rm $f_{\alpha}(x)$}
\epsfxsize=0.5\textwidth \leavevmode
\epsffile{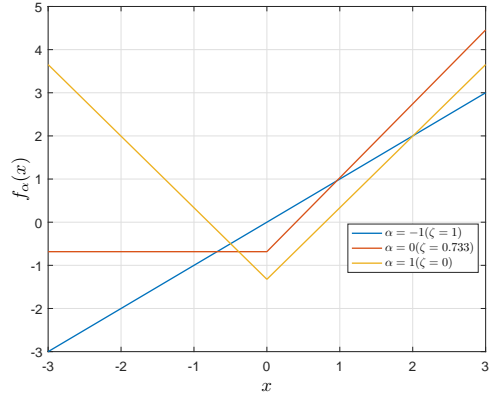}
\caption{The designed activation function $f_{\alpha}(x)$ for different $\alpha$.}\label{fig:variant_relu}
\end{center}
\end{figure}

\begin{figure*}[!t]
\begin{center}
\psfrag{x}[cc][cc][.55][0]{$x$}
\psfrag{y}[cc][cc][.45][0]{\rm probability density}
\subfigure[$\alpha=-1$ ($\zeta=1$)]{%
\epsfxsize=0.31\textwidth \leavevmode
\epsffile{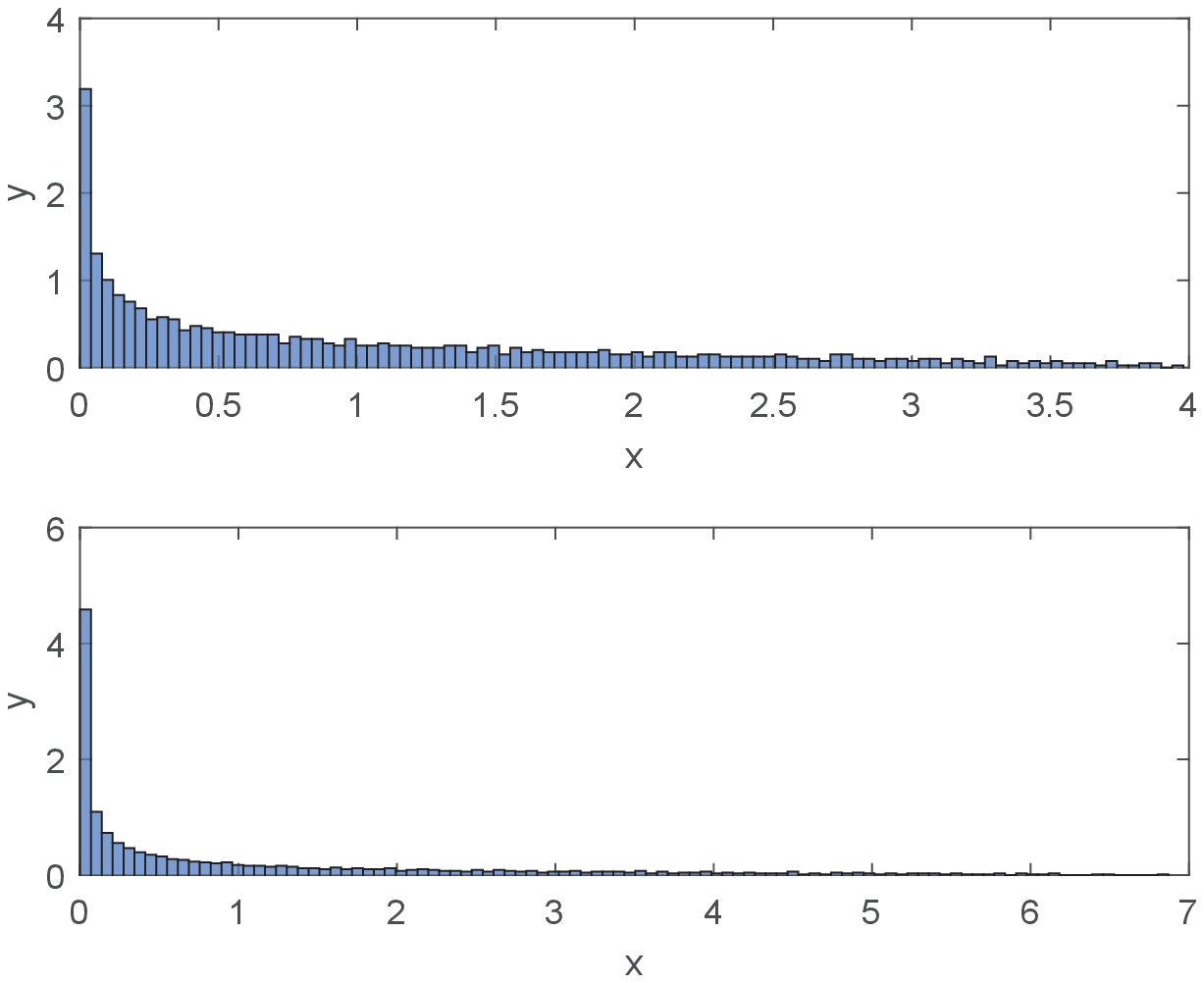}}\quad
\subfigure[$\alpha=0$ ($\zeta=0.733$)]{%
\epsfxsize=0.31\textwidth \leavevmode
\epsffile{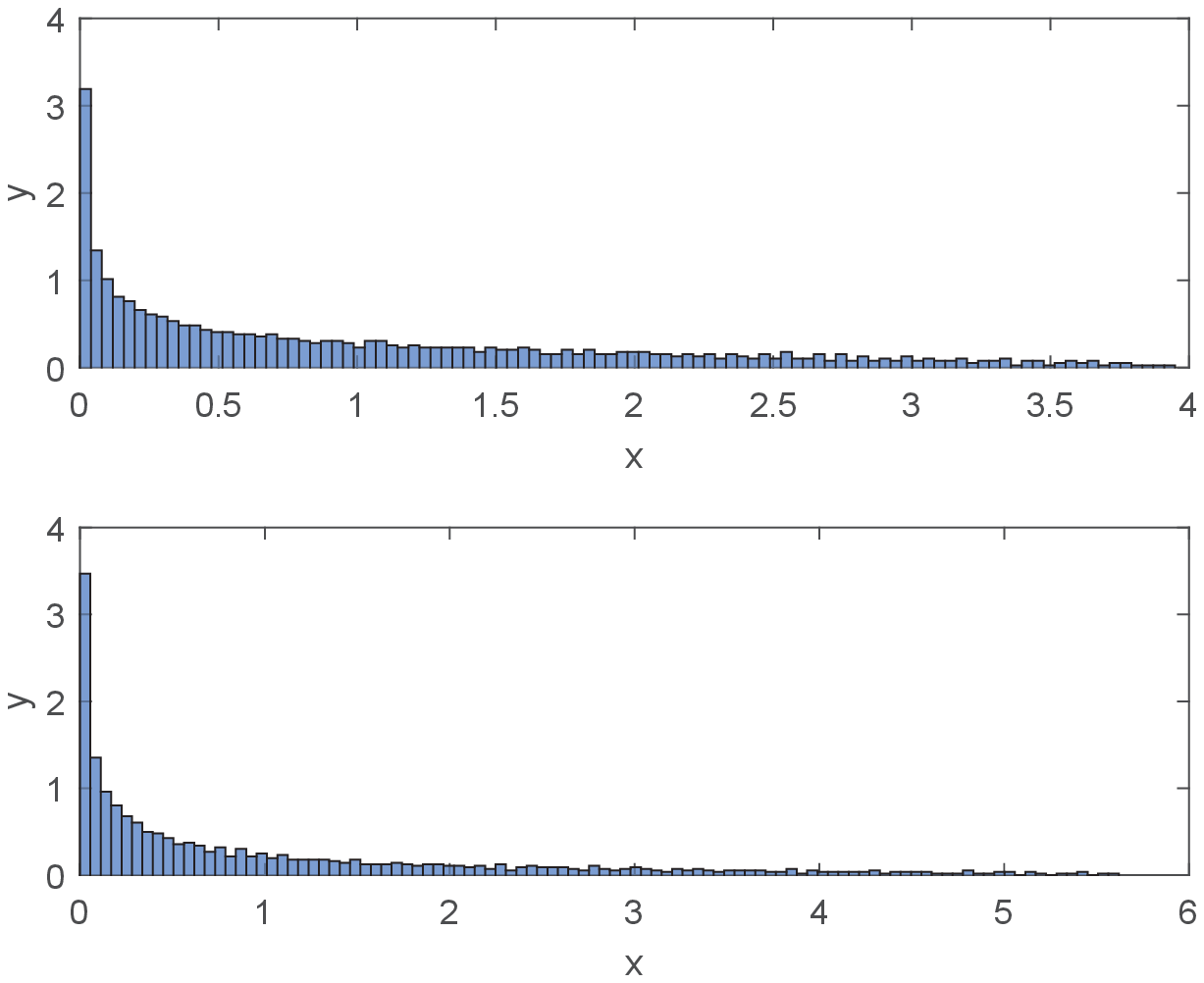}}\quad
\subfigure[$\alpha=1$ ($\zeta=0$)]{%
\epsfxsize=0.31\textwidth \leavevmode
\epsffile{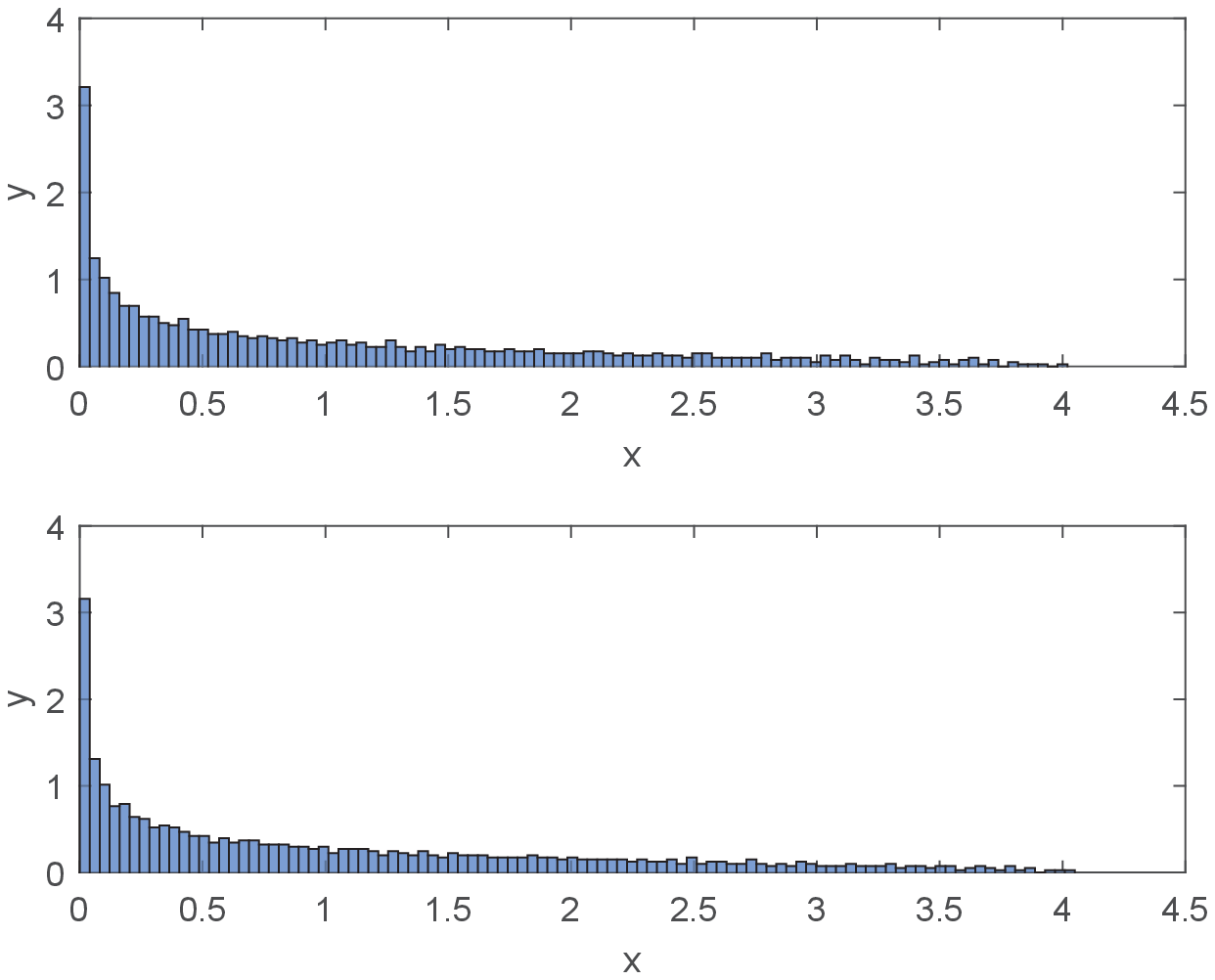}}\quad

\caption{Empirical spectrum density of the input covariance matrix and the output covariance matrix for different $\alpha$ in a single-layer neural network. The upper part and the bottom part of each subgraph show the spectrum of the input covariance matrix and that of the output covariance matrix, respectively.}\label{fig:io_spec_alpha_SL}
\end{center}
\end{figure*}

\begin{figure*}[!t]
\begin{center}
\psfrag{x}[cc][cc][.55][0]{$x$}
\psfrag{y}[cc][cc][.45][0]{\rm probability density}
\subfigure[$\alpha=-1$ ($\zeta=1$)]{%
\epsfxsize=0.31\textwidth \leavevmode
\epsffile{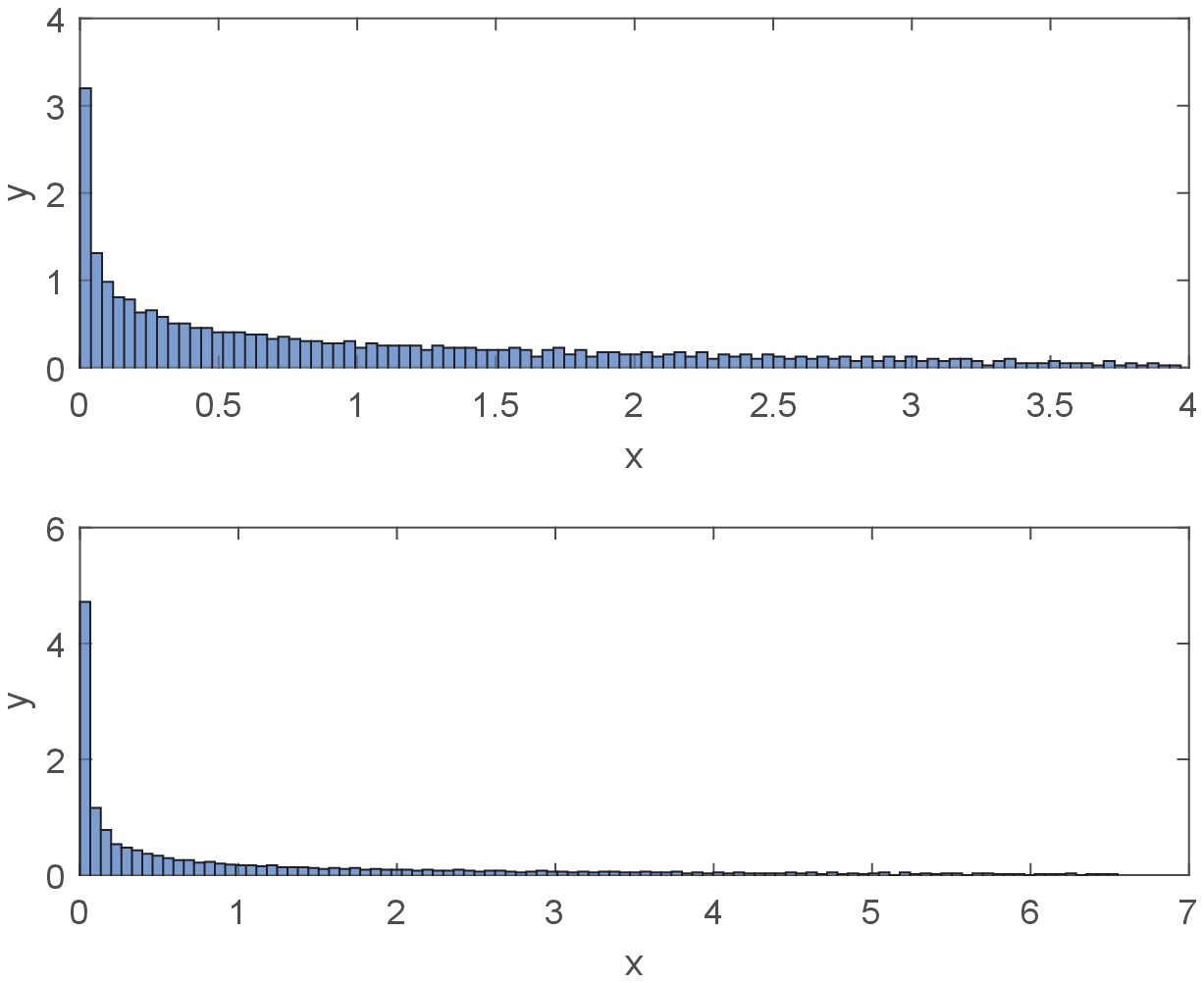}}\quad
\subfigure[$\alpha=0$ ($\zeta=0.733$)]{%
\epsfxsize=0.31\textwidth \leavevmode
\epsffile{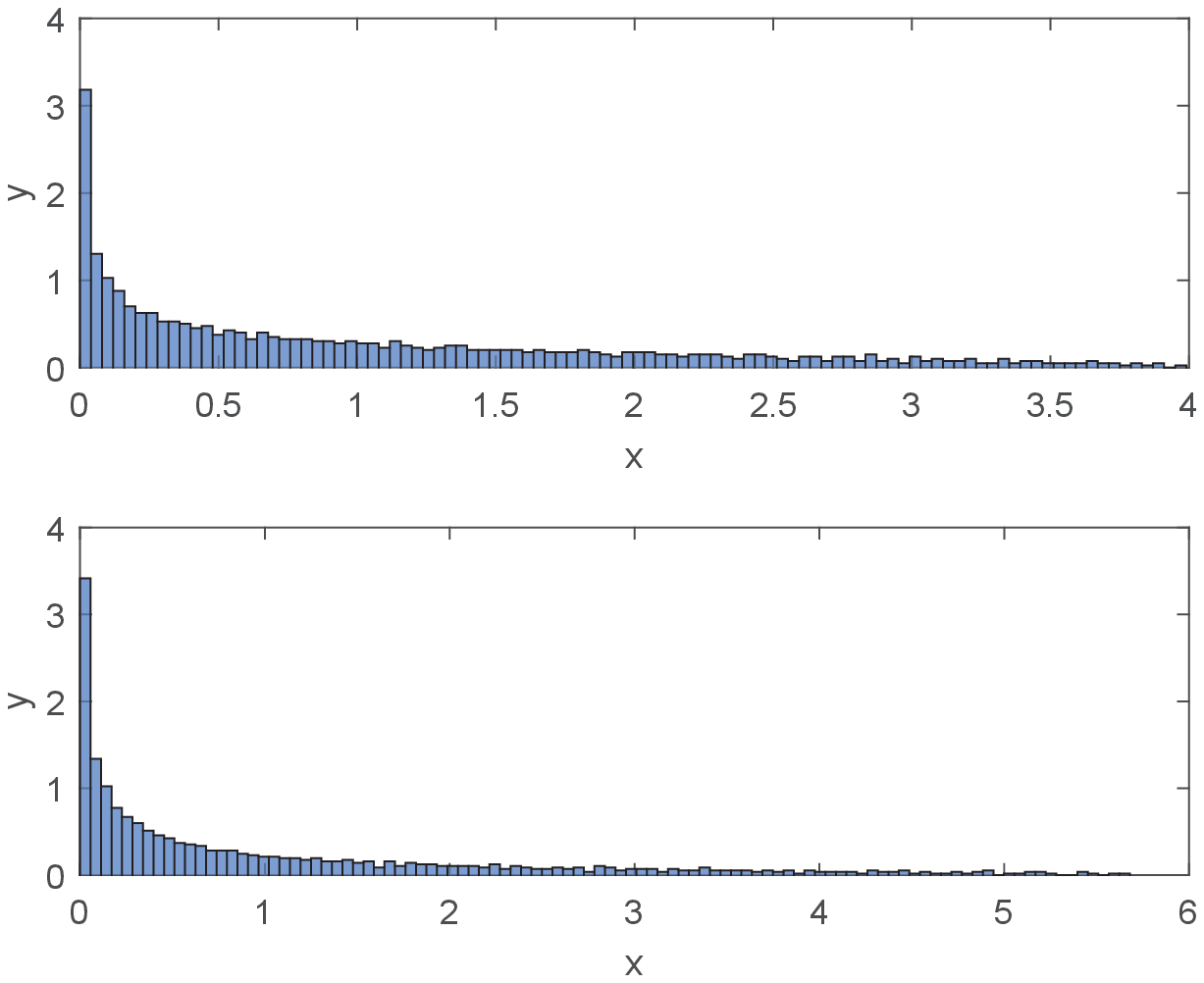}}\quad
\subfigure[$\alpha=1$ ($\zeta=0$)]{%
\epsfxsize=0.31\textwidth \leavevmode
\epsffile{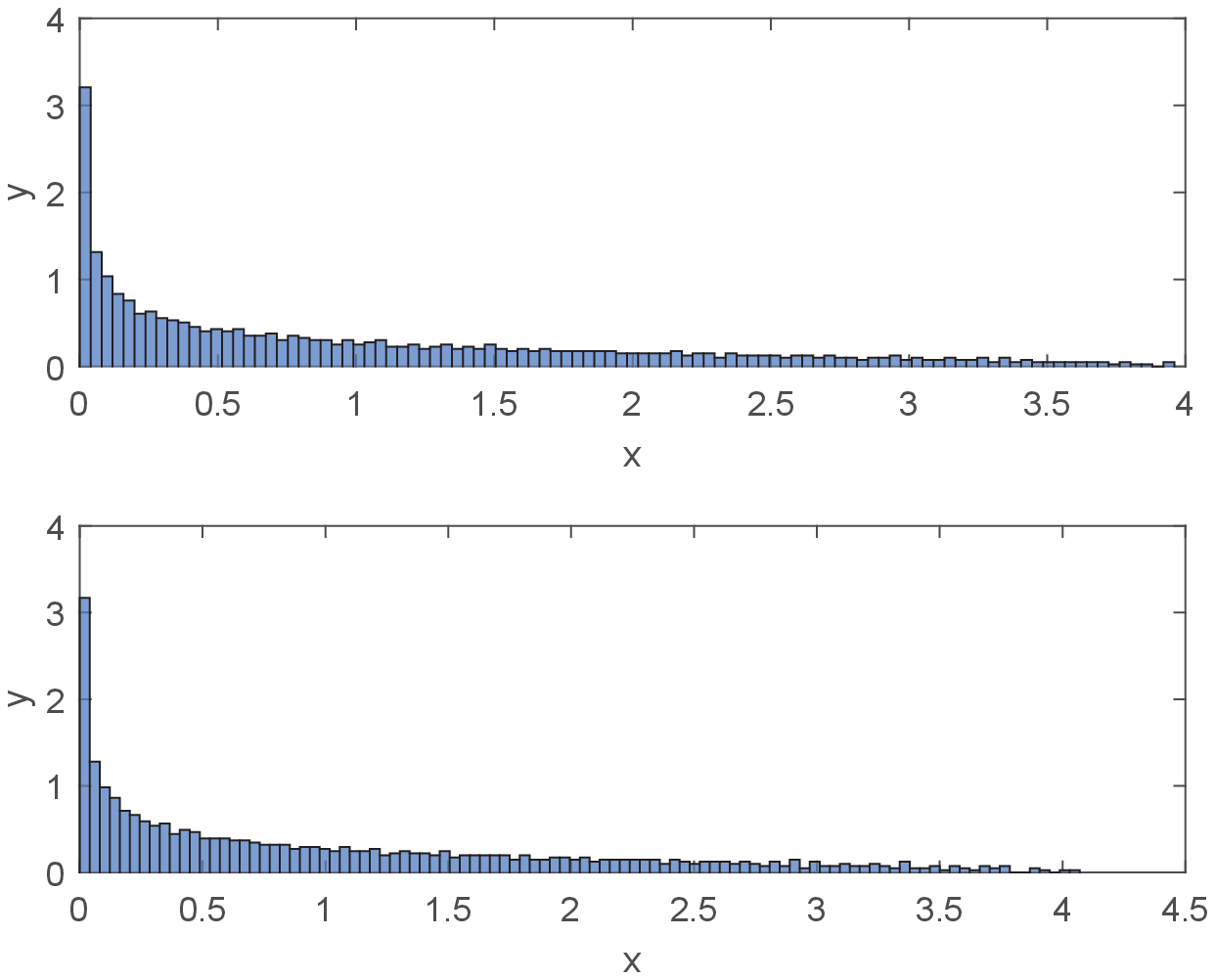}}\quad

\caption{Empirical spectrum density of the input covariance matrix and the output covariance matrix for different $\alpha$ in a neural network with $10$ layers. The upper part and the bottom part of each subgraph show the spectrum of the input covariance matrix and that of the output covariance matrix, respectively.}\label{fig:io_spec_alpha_ML}
\end{center}
\end{figure*}

In \cite{pennington2017nonlinear}, a variant of the ReLU activation function shown as follows is employed to study the impact of $\zeta$,
\begin{equation}\label{eq:vatiant_relu}
f_{\alpha}(x) = \frac{[x]_{+}+\alpha[-x]_{+}-\frac{1+\alpha}{\sqrt{2\pi}}}{\sqrt{\frac{1}{2}(1+\alpha^2)-\frac{1}{2\pi}(1+\alpha)^2}},
\end{equation}
where $\alpha$ is a parameter governing the shape of the activation function, and $\zeta$ can be adjusted by setting $\alpha$. Specifically, $f_{\alpha}(x)$ is the linear activation function and $\zeta=1$ when $\alpha=-1$; $f_{\alpha}(x)$ is the shifted ReLU activation function and $\zeta=0.733$ when $\alpha=0$; $f_{\alpha}(x)$ is the shifted absolute activation function and $\zeta=0$ when $\alpha=1$. The spectra of the input covariance matrix and the output covariance matrix for different activation functions in a single-layer neural network are shown in Fig. \ref{fig:io_spec_alpha_SL}. The corresponding results in a $10$-layer neural network are also shown in Fig. \ref{fig:io_spec_alpha_ML}. Obviously, the spectra of the data covariance matrices are skewed in the neural networks where $\zeta=1$ and $\zeta=0.733$. On the contrary, the spectra are perfectly preserved with $\zeta=0$. It should be highlighted that the spectra can be better preserved with smaller $\zeta$.

In \cite{adlam2019random}, a more general model for random neural networks is considered, i.e., the random biases are considered based on model in \eqref{eq:outputmatrix}. The post-activation matrix of a single-layer neural network is thus given by
\begin{equation}\label{eq:outputmatrixbias}
\mathbf{Y}=\phi(\mathbf{W}\mathbf{X}+\mathbf{B}),
\end{equation}
where $\mathbf{W}$ and $\mathbf{X}$ are as the same as that defined before; $\mathbf{B}=\mathbf{b}\mathbf{1}_m^T\in\mathbb{R}^{n_1\times m}$ (for $\mathbf{b}\in \mathbb{R}^{n_1}$) is the additive random bias matrix. The spectrum of the output covariance matrix is studied under the non-Gaussian data distributions and the non-zero bias distributions. The results in {\it Theorem} \ref{th:mfz_nn} are thus extended into a more general case. In addition, the bias is interpreted as a distribution induced to the activation function parameterized by $\mathbf{B}$, i.e., $\phi(\mathbf{Z}; \mathbf{B}):=\phi(\mathbf{Z}+\mathbf{B})$. Moreover, the analysis can be extended to an arbitrary distribution of activation functions $\phi(:; \mathbf{B})$ parameterized by $\mathbf{B}$. The results are obtained with the similar mathematical tools in \cite{pennington2017nonlinear} but more complex due to the consideration of the random biases. Hence, we do not present the details in this paper.
A quite significant discovery in \cite{adlam2019random} is that, for a specific noisy auto-encoding task, a non-trivial distribution over activation functions can outperform the existing possibly best single activation function. This indicates that the mixtures of nonlinearities might be more useful for approximating the kernel methods or the neural network architecture design. Besides, studying the relations between the spectrum properties of the data covariance matrices and the non-linearities in neural network may give us some inspirations about how to improve the learning speed by designing the nonlinear activation functions.

\re{
It should be noted that the results in \cite{adlam2019random,pennington2017nonlinear}, are obtained under the {\it i.i.d.} Gaussian assumptions on input data and random weights of neural networks. In spirit of the universality widely studied in RMT, \cite{benigni2019eigenvalue} extends the results to sub-Gaussian cases, where both the inputs and the random weights are not necessarily Gaussian. Besides, to further understand the effects of the nonlinear activation functions on the spectra of the data covariance matrices, \cite{benigni2019eigenvalue} derives the results under the cases where the activation functions are polynomial. Thus, \cite{benigni2019eigenvalue} actually extends the results to a more general class of activation functions. On the other hand, \cite{liao2018spectrum} extends the researches into a general case where the input data samples follow a Gaussian mixture model, which is more realistic in practice. Besides, \cite{liao2018spectrum} considers the average kernel matrix, which is the expectation of the output data covariance matrix with respect to the random weights. The mutual influence of different nonlinear activation functions and statistics of input data on the average kernel matrix is quantitatively described. The results reveal that, for different input data statistics, different activation functions have distinct performance on the classification learning task.
}

\subsection{\re{Understanding the Training and Performance of Neural Networks}}\label{subsec:doubledescent}

\re{
Deep neural networks with millions or sometimes even billions of parameters are so powerful that they can fit almost all the possible functional relations between the inputs and outputs. More generally, not only the neural networks, but also the other machine learning algorithms, e.g., support vector machine (SVM), the kernel methods or even more simpler linear regressors, are aimed to fit the training data. In general, the learning models with a large number of parameters can fit the train data very well. However, as the complexity of the learning models increases, the {\it overfitting} phenomenon usually appears. The trained model performs well on the train data set but shows poor performance on the test set. As a consequence, the curve of the prediction error with respect to the model complexity is usually U-shaped. Many techniques, e.g., regularization and dropout, are developed to avoid overfitting. However, recent researches show that deep neural networks and the kernel methods can generalize well even if they interpolate all the train data \cite{zhang2016understanding, belkin2018understand}. The learning models that achieve zero training error, a.k.a. the interpolators, have attracted a lot of attention recently in machine learning because state-of-the-art deep neural networks are belong to the models of this category \cite{hastie2019surprises}. The surprising generalization performance of the interpolators can be well explained by the {\it double descent theory} \cite{belkin2019reconciling}. It suggests that the prediction error decreases first and then increases as the complexity of the model increases under the so-called {\it interpolation threshold}. This corresponds to the conventional overfitting phenomenon. When the complexity of the model continues increasing and exceeds the interpolation threshold, the prediction error decreases again and often converges to the global minimum as the complexity of the model go to infinity \cite{mei2019generalization}.
}

\re{
The double descent phenomenon of the prediction error is first discussed generally in \cite{belkin2019reconciling} and is also observed in \cite{advani2020high, geiger2020scaling}. Here, we emphasize that the double descent phenomenon appears in the extremely complicated learning models, i.e., in the overparametrized regime \cite{hastie2019surprises, mei2019generalization}. In particular, the prediction error of the linear regression learning models is analytically derived in the asymptotic regime, where both the dimension of the learning model and the number of samples go to infinity \cite{hastie2019surprises, belkin2020two}. To be specific, \cite{hastie2019surprises} derives the asymptotic prediction error for a general model with correlated covariates and \cite{belkin2019reconciling} obtains the exact formula of the prediction error for {\it i.i.d.} Gaussian covariates. Besides, the asymptotic generalization error of the random features regression model is analyzed in \cite{mei2019generalization} and the results provide the first analytically tractable model capturing the double descent phenomenon without the misspecication structures assumption. Moreover, these works show that the double descent phenomenon of the generalization error can be theoretically analyzed via RMT \cite{hastie2019surprises, mei2019generalization}. In the following, we take the results from \cite{mei2019generalization} as an example to explain why overparametrized learning models perform so well in practice.
}

\re{
We first consider a specific problem of learning a function $f_d\in L^2(\mathbb{S}^{d-1}(\sqrt{d})$ on the $d$ dimensional sphere. Here, $\mathbb{S}^{d-1}(r)$ denotes the sphere of radius $r$ in $d$ dimensions and $r$ can be set to $\sqrt{d}$ without loss of generality. Besides, the {\it i.i.d.} training data samples $\{(\mathbf{x}_i, y_i)\}$ $(i=1, \cdots, n)$ satisfy $\mathbf{x}_i\sim Unif(\mathbb{S}^{d-1}(\sqrt{d})$ and $y_i = f_d(\mathbf{x}_i)+\epsilon_i$, with the {\it i.i.d.} $\epsilon_i$ independent of $\mathbf{x}_i$. The noise distribution is assumed to satisfy $\mathbb{E}_{\epsilon}(\epsilon_i)=0$, $\mathbb{E}_{\epsilon}(\epsilon_i^2)=\tau^2$ and $\mathbb{E}_{\epsilon}(\epsilon_i^4)<\infty$. Moreover, we consider the case where the data samples are fitted with the random features (RF) model, which is equivalent to the following function class
\begin{equation}\label{eq:RFmodel}
\mathcal{F}_{RF}(\mathbf{\Theta}) = \left\{f(\mathbf{x};\mathbf{a}, \mathbf{\Theta})\triangleq\sum_{i=1}^{N}a_i\phi\left(\langle\mathbf{\theta}_i, \mathbf{x}\rangle/\sqrt{d}\right)\right\},
\end{equation}
where $\langle\cdot,\cdot\rangle$ is the inner product operation.
Here, the random features model can be regarded as a single layer neural network where the weights $\mathbf{\Theta}\in\mathbb{R}^{N\times d}$ between the inputs and the pre-activations of the hidden layer are randomly chosen. $\mathbf{\theta}_i$, satisfying $\|\mathbf{\theta}_i\|_2=\sqrt{d}$, denotes the $i$-th row of $\mathbf{\Theta}\in\mathbb{R}^{N\times d}$, $\phi(\cdot)$ is the element-wise activation function and $\mathbf{a}=[a_i, \cdots, a_N]^T\in\mathbb{R}^{N}$ denote the weights between the post-activations of the hidden layer and the output. The training of the random features model is quite different from that of neural networks since only $\mathbf{a}$ needs to be trained. In general, $\mathbf{a}$ can be learnt by performing ridge regression
\begin{equation}\label{eq:ridgeregression}
\hat{\mathbf{a}}(\lambda) = \mathop{\arg\min}_{\mathbf{a}\in\mathbb{R}^{N}}\left\{\frac{1}{n}\sum_{j=1}^{n}\left(y_j-\sum_{i=1}^{N}a_i\phi\left(\langle\mathbf{\theta}_i, \mathbf{x}\rangle/\sqrt{d}\right)\right)^2+\frac{N\lambda}{d}\|\mathbf{a}\|_2^2\right\},
\end{equation}
where $\lambda$ is the regularization factor of the ridge regression. In addition, the ridge regularization path is shown to be closely related to the path of gradient flow when the mean square error (MSE) $\sum_{j=1}^{n}(y_j-f(\mathbf{x}_i;\mathbf{a}, \mathbf{\Theta}))^2$ is adopted. Particularly, the convergence point of the gradient flow is exactly the ridgeless limit of $\hat{\mathbf{a}}(\lambda)$, i.e., $\lim_{\lambda\to0}\hat{\mathbf{a}}(\lambda)$ and a positive $\lambda$ corresponds to an early stopping of the gradient descent procedure \cite{yao2007early}.
}

\re{
The prediction error (a.k.a. test error, generalization error or risk) is the expectation of the MSE with respect to the test data $\mathbf{x}\sim Unif(\mathbb{S}^{d-1}(\sqrt{d})$, which is independent of the train data. Denoting the train data samples with $\mathbf{X}=[\mathbf{x}_1, \cdots, \mathbf{x}_n]$, the prediction error of the random features model, $R_{RF}(f_d, \mathbf{X}, \mathbf{\Theta}, \lambda)$, is given by
\begin{equation}\label{eq:predictionerror}
R_{RF}(f_d, \mathbf{X}, \mathbf{\Theta}, \lambda) = \mathbb{E}_{\mathbf{x}}\left[\left(f_d(\mathbf{x})-f(\mathbf{x}; \hat{\mathbf{a}}(\lambda), \mathbf{\Theta})\right)^2\right].
\end{equation}
Note that we only take expectation with respect to $\mathbf{x}$. It is not important since $R_{RF}(f_d, \mathbf{X}, \mathbf{\Theta}, \lambda)$ concentrates around $\bar{R}_{RF}(f_d,\lambda)\triangleq\mathbb{E}_{\mathbf{X}, \mathbf{\Theta}, \mathbf{\epsilon}}R_{RF}(f_d, \mathbf{X}, \mathbf{\Theta}, \lambda)$ \cite{mei2019generalization}.
}

\re{
With the above analysis, the accurate approximation for the prediction error in the asymptotic regime ($d, n, N \to\infty$) can be derived via RMT. The derivations in \cite{mei2019generalization} are quite complicated, thus we here only present an informal overview of the results. With the following two ratios
\begin{equation}\label{eq:psi12}
\psi_1=\frac{N}{d}, \phi_2=\frac{n}{d}, \ {\rm as} \ d, n, \to\infty,
\end{equation}
the overparametrization ratio is defined as $\gamma=\psi_1/\psi_2=N/n$ \cite{hastie2019surprises}. $\gamma<1$ means the underparametrized regime while $\gamma>1$ means the overparametrized regime. The prediction error depends on $f_d(\cdot)$ (the characteristics of the function to be learnt), $\phi(\cdot)$ (the activation funcion), $\psi_1$, $\psi_2$, and $\tau^2$ (the noise variance). From the results in \cite{mei2019generalization}, the asymptotic ridgeless (the case where $\lambda\to0$) prediction error goes through decreasing-increasing-decreasing process as $\gamma$ increases. In addition, the global minimum of the prediction error is achieved in the highly overparametrized regime. This is exactly the double descent phenomenon (see Figure $3$ in \cite{mei2019generalization}). Moreover, for the specific regression problems with random feature kernels, the double descent phenomenon can be eliminated via optimal regularization and the prediction error monotonically decreases as $\gamma$ increases. This exactly justifies the effect of the regularization in avoiding overfitting.
}

\re{
Besides, there exists another structure of random neural networks, i.e., {\it extreme learning machine} (ELM) \cite{huang2011extreme}, which is quite similar to the random features model. The ELM can be described as
\begin{equation}\label{eq:ELM}
\hat{\mathbf{y}} = \mathbf{\beta}^T\phi(\mathbf{W}\mathbf{x}),
\end{equation}
where $\mathbf{x}\in\mathbb{R}^p$ is the input data, $\mathbf{W}\in\mathbb{R}^{n\times p}$ is a random weight matrix, $\mathbf{\beta}\in\mathbb{R}^{n\times d}$ is the coefficient matrix that maps the random feature $\phi(\mathbf{W}\mathbf{x})$ to the output $\hat{\mathbf{y}}\in\mathbb{R}^d$, and $\phi(\cdot)$ is the element-wise activation function. With the train data, the only trainable $\mathbf{\beta}$ can be trained quickly via ridge regression. Obviously, the ELM is almost the same with the random features model except the vectorial output. In \cite{louart2018random}, the asymptotic training error and generalization error are derived via RMT and are shown to depend on the hyper-parameters of the ELM. The results provide useful insights into the underlying mechanism of ELM and also give practical ways to tune the hyper-parameters. Beyond the feed-forward neural networks, the limiting training error and generalization performance of {\it linear echo state neural networks}, which are actually a class of RNNs, are analytically derived in the asymptotic regime \cite{couillet2016asymptotic}. The asymptotic results provide further new insights into the performance of more advanced neural networks.
}

\re{
Actually, the random features model \cite{rahimi2007random} can be regarded as not only a single-hidden layer neural network with random first layer weights, but also a random approximation of a kernel regression. Intuitively, the training of random features model can be divided into two parts: i) obtain the representation of the input $\mathbf{x}$ in the random feature space, namely, $[\phi(\langle\mathbf{\theta}_1, \mathbf{x}\rangle/\sqrt{d}), \cdots, \phi(\langle\mathbf{\theta}_N, \mathbf{x}\rangle/\sqrt{d})]$, via the random feature kernel. ii) perform a ridge regression between the kernel representation of $\mathbf{x}$ and the labels $y$ to learn the regression coefficients. \cite{mei2019generalization} also points that $\mathcal{F}_{RF}(\mathbf{\Theta})$ is indeed a reproducing kernel Hilbert space (RKHS) defined by the finite-rank approximation of the following kernel
\begin{equation}\label{eq:RKHS}
\mathcal{H}_N(\mathbf{x}, \mathbf{x}')=\frac{1}{N}\sum_{i=1}^{N}\phi\left(\langle\mathbf{\theta}_i, \mathbf{x}\rangle/\sqrt{d}\right)\phi\left(\langle\mathbf{\theta}_i, \mathbf{x}'\rangle/\sqrt{d}\right).
\end{equation}
}

\re{
Indeed, the neural networks is closely related to the kernel methods in machine learning. This is quite intuitive in the so-called {\it lazy training} regime, where the parameters of the neural networks change not much in the training process \cite{hastie2019surprises}. Consider a neural network with parameters $\mathbf{\theta}$ whose function is $f(\cdot; \mathbf{\theta}):\mathbb{R}^{d}\to\mathbb{R}$, $\mathbf{x}\mapsto f(\mathbf{x}; \mathbf{\theta})$. Assuming a random initialization for $\mathbf{\theta}$, say $\mathbf{\theta}_0$, makes $f(\mathbf{x}; \mathbf{\theta}_0)\approx0$, denote the parameters after training by $\mathbf{\theta}=\mathbf{\theta}_0+\mathbf{\beta}$ and $\mathbf{\beta}$ is small in the lazy training regime, we have the following approximate results with Taylor expansion:
\begin{equation}\label{eq:lazytraining}
\mathbf{x}\mapsto \nabla_{\mathbf{\theta}}f(\mathbf{x}; \mathbf{\theta}_0)^T\mathbf{\beta}.
\end{equation}
As we can see, the model in \eqref{eq:lazytraining} is linear in $\mathbf{\beta}$. We can now training the neural network by performing simple ridge regression with the known $\nabla_{\mathbf{\theta}}f(\mathbf{x}; \mathbf{\theta}_0)$. Thus, the asymptotic results derived for the ridge regression are closely related to the neural networks with lazy training.
}

\re{
We stress that the intuitive observation only holds in the lazy training regime. But this still provides us a train of thought to study the training dynamics of neural networks via kernel methods in a more general way. A recent line of researches \cite{jacot2018neural, li2018learning,oymak2020toward} show that the training dynamics of neural networks can be studied via the {\it Neural Tangent Kernel} (NTK), i.e., the training of neural networks can also be divided into two parts: i) learn the NTK which maps the input $\mathbf{x}$ to learning representations in another feature space. ii) perform ridge regression to learn the `regression coefficients' between the learning representations and the labels. Another kernel of interest is the {\it Conjugate Kernel} (CK), which also governs the training process and the generalization performance of neural networks.
}

\re{
In particular, the spectral properties of the two kernel matrices are closely related to training and generalization of neural networks \cite{fan2020spectra}. For example, the gradient descent process can be accelerated along the eigenvectors of the largest eigenvalues \cite{advani2020high}. Besides, the spectral distributions indicate the trainability and the extent of implicit bias towards simpler functions \cite{yang2019fine, xiao2020disentangling}. To introduce the two kernels, we consider the case where we use a neural network with $L$ hidden layers to fit the train data samples $\{(\mathbf{x}_i, y_i)\}$ ($i=1, \cdots, m$), the network outputs of the $m$ data samples $\hat{\mathbf{y}}=[\hat{y}_1, \cdots, \hat{y}_m]^T$ are given by
\begin{equation}\label{eq:network_model_Llayers}
\hat{\mathbf{y}} = \mathbf{w}^T\mathbf{X}^{L}, \mathbf{X}^{l} =\phi(\mathbf{W}^{l}\mathbf{X}^{l-1}), \mathbf{X}^0=\mathbf{X}, l=1,\cdots, L,
\end{equation}
where $\phi(\cdot)$ is the activation function, $\mathbf{w}\in\mathbb{R}^{n_L}$ is the coefficients that map $L$-th layer post-activations to the network ouput, $\mathbf{X}^{l}$ denote the post-activation matrix of $l$-th layer, $\mathbf{X}^0=\mathbf{X}=[\mathbf{x}_1, \cdots, \mathbf{x}_m]\in\mathbb{R}^{n_0\times m}$ is the matrix composed of $m$ input data in the train set, $\mathbf{W}^{l}\in\mathbb{R}^{n_l\times n_{l-1}}$ is the $l$-th layer weight matrix with $n_l$ the number of neurons of $l$-th layer, and $n_0$ denotes the input dimension.
The conjugate kernel is defined as the gram matrix of the post-activations of the final hidden layer, i.e.,
\begin{equation}\label{eq:CK}
\mathbf{K}^{CK}\triangleq(\mathbf{X}^{L})^T\mathbf{X}^{L}\in\mathbb{R}^{m\times m}.
\end{equation}
Besides, we use a weight vector $\mathbf{\theta}=[vec(\mathbf{W}^{1}), \cdots, vec(\mathbf{W}^{L}), \mathbf{w}]$ to denote the all the weights in the neural network. The network outputs can be rewritten as a function of the input data samples, i.e.,
\begin{equation}\label{eq:inputoutputfunction}
\hat{\mathbf{y}} = f_{\mathbf{\theta}}(\mathbf{X}) \ {\rm or} \ \hat{y}_i=f_{\mathbf{\theta}}(\mathbf{x}_i).
\end{equation}
Then we can obtain the Jacobian matrix of the network outputs with respect to the weight vector by
\begin{equation}\label{eq:weightJacobian}
\mathbf{J}=\nabla_{\mathbf{\theta}}f_{\mathbf{\theta}}(\mathbf{X})=[\nabla_{\mathbf{\theta}}f_{\mathbf{\theta}}(\mathbf{x}_1), \cdots, \nabla_{\mathbf{\theta}}f_{\mathbf{\theta}}(\mathbf{x}_m)]\in\mathbb{R}^{dim(\mathbf{\theta})\times m}.
\end{equation}
The neural tangent kernel is defined as
\begin{equation}\label{eq:NTK}
\mathbf{K}^{NTK}\triangleq\mathbf{J}^T\mathbf{J}\in\mathbb{R}^{m\times m}.
\end{equation}
}

\re{
With the Stieltjes transform, the limiting spectra of CK and NTK of the $L$-layer neural network are derived in the asymptotic regime \cite{fan2020spectra}, where both network width and the sample size grow to infinity with a constant ratio. The results in \cite{fan2020spectra} are actually great extensions of the researches about the training and generalization errors in linear regressions \cite{hastie2019surprises} and random features model \cite{mei2019generalization}. In addition, the studies on the spectra of CK and NTK of neural networks enable the analysis for the feature learning, training and generalization of neural networks in some more general scenarios, instead of only the lazy training regime. Last but not least, the experimental results in \cite{fan2020spectra} show that the spectra of CK and NTK appear interesting evolutions during training, the researches for random weight neural networks may shed some light on studying the interesting evolutions during training.
}

\section{Challenges and Opportunities}\label{sec:challenges}

One can see that RMT is a powerful mathematical tool to deal with the extremely large dimensional data and to analyze the large complex systems. However, there are still some critical challenges and opportunities that should be addressed.

\subsection{Complex Statistics of the Random Matrices}\label{subsec:complex_statistics}

As we can see, the major results about the specific eigenvalues in RMT mainly focus on the extreme eigenvalues. For some eigenvalue-based spectrum sensing algorithms, this will prevent us to determine the detection threshold and to evaluate the detection performance analytically. As an example, the detection threshold and the detection probability of the AGM method, which exploits the arithmetic mean to geometric mean of the eigenvalues of the sample covariance matrix, are hard to compute due to the complex test statistic. Hence, more advanced results about the complex statistics are expected to be derived so that more complex problems can be analyzed.

\subsection{Imperfect Randomness in Practice}\label{subsec:inperfect_randomness}

Most results in RMT can be regarded as the analogies with the {\it concentration of measure} phenomenon \cite{tao2012topics} in probability theory, and their validity relies on the independence for the entries of the random matrices. However, the independence of the entries of random matrices in practical scenarios may not hold perfectly. For example, the discrete noise samples may be not {\it i.i.d.} due to the non-ideal sample filter design \cite{zeng2010review}, and this will cause many spectrum sensing methods out of gear. Fortunately, the noise prewhitening technique can be used to solve this problem \cite{zeng2009eigenvalue}. On the other hand, with the fact that massive antennas and higher carrier frequencies will be employed in the 5G and beyond communications systems, the channel statistics become quite different. The {\it i.i.d} Rayleigh fading channels should be modified with the Rician fading channel models, in which the constant line-of-sight (LOS) components exist. This can also make the assumptions about the independence not hold true. Therefore, the asymptotic analysis for the future large complex communication systems is quite challenging. It is quite interesting to study how much impact will the imperfect randomness has on the results obtained under the perfect {\it i.i.d.} assumption.


\subsection{Demand for New Technical Tools}\label{subsec:new_models}

\re{
Deep neural networks are extremely powerful in exploiting nonlinear features from the data. Intuitively, we have to develop new technical tools to analyze the nonlinear random matrix models.
}
Also, many theoretical analysis for neural networks are based on quite simple neural networks with equally wide layers or without biases. We expect to get some inspirations from these simplified neural networks, but these simplifications also make the conclusions deviate from the practical results. For example, a recent work, i.e., \cite{granziol2020beyond}, shows that the observed spectral shapes of practical neural networks and datasets strongly deviate from the theoretical results. In addition, the products of random Wigner/Wishart matrices and the percolated Wigner/Wishart matrices are found to be better in approximating the practical spectra. This indicates that new tools are needed to make the theoretical analysis more practical.

\subsection{\re{Wish for Universal Theories}}\label{subsec:universality}
\re{
In the theoretical analysis of deep learning/general machine learning techniques, the input data or the network weights are usually assume to be {\it i.i.d.} Gaussian. These assumptions are quite strong and may diverge a lot from the practical scenarios. Hence, one wishes to build theories on a statistical model capturing the practical domain-specific data (e.g., the shift- and rotation-invariant property of images), beyond the simple {\it i.i.d.} Gaussian modeling of the data. For example, \cite{seddik2020random} shows that the deep learning representations of GAN-data behaves as Gaussian mixture model (GMM) via a concentration of measure approach. Moreover, the impacts of nonlinearities on the classification performance are studied under the Gaussian mixture data model in \cite{liao2018spectrum}. Besides, to analyze the input-output Jacobian of neural networks, \cite{pastur2020random1} proposed a general analytical framework which accounts for {\it i.i.d.} random weights but non-necessarily Gaussian. On the other hand, as discussed in Section \ref{subsec:dataCM}, the spectral behaviors of the nonlinear matrices produced by deep neural networks depends on the nonlinearities via a few parameters. Many kinds of nonlinearities have been well studied in the line of works \cite{benigni2019eigenvalue, pennington2017resurrecting, pennington2018emergence}. In spirit of the {\it universality} widely studied in RMT, one may wish universal theories for more general data/weight distributions, and nonlinearities.
}

\section{Conclusions}\label{sec:conclusions}
In this paper, we have investigated the applications of RMT in wireless communications and deep learning. First, we have reviewed the basic concepts and the well-known results in RMT. Then, we have introduced some typical applications in wireless communications: designing the spectrum sensing algorithms for the cognitive radio systems and analyzing the asymptotic performance of the multiuser receivers for the large communication systems. Afterwards, we have provided an overview of the applications in understanding and improving the emerging deep neural networks. In particular, we have respectively introduced the RMT-based analysis methods for studying the spectra of the Hessian, Jacobian and data covariance matrix the of the neural networks. We also have presented the works devoting to understanding the training and generalization performance by analyzing the limit training error, generalization error and the related kernel matrices of neural networks. Finally, we have highlighted the challenges and opportunities in applying RMT to the practical large complex systems. We hope this article can establish a connection between engineering applications and mathematical field in which RMT will keep to be powerful.

\bibliographystyle{ws-jktr}
\bibliography{IEEEabrv, ref_RMT_in_WCDL}

\end{document}